\documentclass[preprint,11pt]{elsarticle}
\biboptions{compress} 

\usepackage[a4paper]{geometry}
\usepackage{lmodern} 

\usepackage{amssymb}
\usepackage{amsfonts}
\usepackage{amsmath}
\usepackage{mathtools} 
\usepackage{xfrac} 
\usepackage{gensymb} 
\usepackage{booktabs} 
\usepackage{siunitx} 

\usepackage{xcolor} 
\usepackage{graphicx}

\usepackage{caption}

\usepackage{floatrow}
\usepackage{pgf}
\usepackage{subfig}

\usepackage{multirow}
\usepackage{dashrule}

\usepackage{acro} 
\usepackage{comment} 
\usepackage{lineno} 
\usepackage[ruled, vlined, english, onelanguage, algo2e]{algorithm2e}

\usepackage{mdframed}

\definecolor{julia-purple}{HTML}{9558B2}
\definecolor{julia-green}{HTML}{389826}
\definecolor{julia-blue}{HTML}{4063D8}
\definecolor{julia-red}{HTML}{CB3C33}

\definecolor{gray03}{gray}{0.3}
\definecolor{gray04}{gray}{0.4}
\definecolor{gray08}{gray}{0.8}
\definecolor{gray09}{gray}{0.9}

\definecolor{minted-red}{HTML}{b00040}

\definecolor{gray08}{gray}{0.8}

\usepackage{makecell} 
\usepackage{hyperref}
\usepackage[capitalise]{cleveref}
\Crefname{equation}{}{}

\newcommand{\orcid}[1]{\href{https://orcid.org/#1}{\includegraphics[width=10pt]{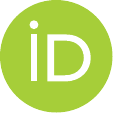}}}

\numberwithin{equation}{section}

\usepackage{fancyhdr}
\usepackage{arydshln} 

\newcommand*\nid{\text{d} }   
\newcommand{\pd}[2]{\frac{\partial #1}{\partial #2} } 
\newcolumntype{?}[1]{!{\vrule width #1pt}} 

\usepackage{tikz}

\usetikzlibrary{patterns}
\usetikzlibrary{backgrounds}
\usetikzlibrary{calc}
\usetikzlibrary{arrows.meta}
\usetikzlibrary{decorations.pathreplacing}

\Crefname{algocf}{Algorithm}{Algorithms}

\DeclareAcronym{ODE}{
	short = ODE,
	long  = ordinary differential equation
}

\DeclareAcronym{IVP}{
	short = IVP,
	long  = initial value problem
}

\DeclareAcronym{RHS}{
	short = RHS,
	long = right-hand side
}

\DeclareAcronym{PDE}{
	short = PDE,
	long  = partial differential equation
}

\DeclareAcronym{MoL}{
	short = MoL,
	long  = method of lines
}

\DeclareAcronym{SSP}{
	short = SSP,
	long  = strong stability preserving
}

\DeclareAcronym{TVD}{
	short = TVD,
	long  = total variation diminishing
}

\DeclareAcronym{TVB}{
	short = TVB,
	long  = total variation bounded
}

\DeclareAcronym{DG}{
	short = DG,
	long  = Discontinuous Galerkin
}

\DeclareAcronym{RKDG}{
	short = RKDG,
	long  = Runge-Kutta Discontinuous Galerkin
}

\DeclareAcronym{DGSEM}{
	short = DGSEM,
	long  = discontinuous Galerkin spectral element method
}

\DeclareAcronym{HLLC}{
	short = HLLC,
	long  = Harten-Lax-Van Leer-Contact
}

\DeclareAcronym{HLLE}{
	short = HLLE,
	long  = Harten-Lax-Van Leer-Einfeldt
}

\DeclareAcronym{HLL}{
	short = HLL,
	long  = Harten-Lax-Van Leer
}

\DeclareAcronym{MHD}{
	short = MHD,
	long  = magnetohydrodynamics
}

\DeclareAcronym{vrMHD}{
	short = VRMHD,
	long  = visco-resistive Magnetohydrodynamics
}

\DeclareAcronym{GLM}{
	short = GLM,
	long  = generalized Lagrangian multiplier
}

\DeclareAcronym{glm-mhd}{
	short = GLM-MHD,
	long  = generalized Lagrangian multiplier magnetohydrodynamics
}

\DeclareAcronym{PRKM}{
	short = PRKM,
	long = Partitioned Runge-Kutta method
}
\DeclareAcronym{ARKMs}{
	short = ARKMs,
	long = Additive Runge-Kutta methods
}
\DeclareAcronym{PERK}{
	short = P-ERK,
	long = Paired Explicit Runge-Kutta
}

\DeclareAcronym{PERRK}{
	short = P-ERRK,
	long = Paired Explicit Relaxation Runge-Kutta
}

\DeclareAcronym{RRKM}{
	short = RRKM,
	long = Relaxation Runge-Kutta method
}

\DeclareAcronym{IMEX}{
	short = IMEX,
	long = implicit-explicit
}

\DeclareAcronym{AMR}{
	short = AMR,
	long = adaptive mesh refinement
}

\DeclareAcronym{CFL}{
	short = CFL,
	long = Courant-Friedrichs-Lewy
}

\DeclareAcronym{DoF}{
	short = DoF,
	long = degree of freedom
}

\DeclareAcronym{IDT}{
	short = IDT,
	long = incremental direction technique
}

\DeclareAcronym{SBP}{
	short = SBP,
	long = summation-by-parts
}

\DeclareAcronym{BR1}{
	short = BR1,
	long = Bassi-Rebay 1
}

\DeclareAcronym{EC/ES}{
	short = EC/ES,
	long = entropy-conservative/entropy-stable
}

\DeclareAcronym{CRM}{
	short = CRM,
	long = Common Research Model
}

\journal{??}

\begin{document}

\mdfsetup{
	backgroundcolor=gray08,
	roundcorner=10pt,
	topline=false,
	rightline=false,
	bottomline=false,
	leftline=false}
	
	\begin{frontmatter}
		
		
		
		\title{Paired Explicit Relaxation Runge-Kutta Methods: \\
		Entropy-Conservative and Entropy-Stable \\High-Order Optimized Multirate Time Integration}
		
		
		\cortext[cor1]{Corresponding author}
		
		\author[1]{Daniel Doehring\corref{cor1} \orcid{0009-0005-4260-0332}}
		\author[2]{Hendrik Ranocha \orcid{0000-0002-3456-2277}}
		\author[1]{and Manuel Torrilhon \orcid{0000-0003-0008-2061}}
		
		\affiliation[1]{
			organization={Applied~and~Computational~Mathematics,~RWTH~Aachen~University},
			country={Germany}.
		}

		\affiliation[2]{
			organization={Institute~of~Mathematics,~Johannes~Gutenberg~University~Mainz},
			country={Germany.}
		}

		\begin{abstract}
			We present novel entropy-conservative and entropy-stable multirate Runge-Kutta methods based on Paired Explicit Runge-Kutta (P-ERK) schemes with relaxation for conservation laws and related systems of partial differential equations.
			Optimized schemes up to fourth-order of accuracy are derived and validated in terms of order of consistency, conservation of linear invariants, and entropy conservation/stability.
			We demonstrate the effectiveness of these P-ERRK methods when combined with a high-order, entropy-conservative/stable discontinuous Galerkin spectral element method on unstructured meshes.
			The Paired Explicit \textit{Relaxation} Runge-Kutta methods (P-ERRK) are readily implemented for partitioned semidiscretizations arising from problems with equation-based scale separation such as non-uniform meshes.
			We highlight that the relaxation approach acts as a time-limiting technique which improves the nonlinear stability and thus robustness of the multirate schemes.
			The P-ERRK methods are applied to a range of problems, ranging from compressible Euler over compressible Navier-Stokes to the visco-resistive magnetohydrodynamics equations in two and three spatial dimensions.
			For each test case, we compare computational load and runtime to standalone relaxed Runge-Kutta methods which are outperformed by factors up to four.
			All results can be reproduced using a publicly available repository.
		\end{abstract}
		
		%
		%
		
		%
		%
		
		
	\end{frontmatter}
	
	
	
	
	%
	\section{Introduction}
	The concept of entropy is of pivotal importance in both analysis and numerical discretization of hyperbolic conservation laws and related systems of \acp{PDE}.
	This insight dates back to the seminal works by Oleĭnik \cite{oleinik1959uniqueness}, Lax 
	\cite{lax1957hyperbolic2}, 
	and Kružkov \cite{kruvzkov1970first} who realized that entropy provides the sought guidance to resolve the non-uniqueness issue caused by shocks.
	Employing entropy principles allowed the aforementioned authors to answer questions originally raised by Riemann, Rankine, and Hugoniot which remained open for more than 70 years.
	Besides being an indispensable tool for the analysis of conservation laws \cite{lax1973hyperbolic, dafermos2005hyperbolic}, entropy is at its heart 
	a physical, thermodynamic quantity.
	All physical processes are irreversible, i.e., carry an inherent direction and thus lead to an increase of the physical entropy.
	A prominent example from hydrodynamics is the entropy increase across shocks \cite{courant1999supersonic}.
	Another notoriously complicated hydrodynamic feature are scale-spanning turbulent flows which continue to challenge various branches of science and engineering.
	Recent studies on \ac{EC/ES} numerical schemes \cite{gassner2016split, SJOGREEN2018153, MANZANERO2020109241, rojas2021robustness, chan2022entropy} have shown how \ac{EC/ES} schemes can be used to obtain more robust simulations for under-resolved turbulent flows.
	Thus, to accurately simulate such and even more intricate phenomena, entropy can and should be employed as a guiding principle in addition to the traditional balance laws governing mass, momentum, and energy.

	The construction of \ac{EC/ES} finite volume schemes has been pioneered by Tadmor \cite{tadmor1986entropy, tadmor1987numerical} and was since then extended to a variety of high-order accurate discretization schemes, see for instance \cite{fjordholm2012arbitrarily, FISHER2013518, gassner2013skew, carpenter2016towards} and references therein.
	In this work we focus on the \ac{DGSEM} method, i.e., a nodal \ac{DG} method \cite{hesthaven2007nodal} with spatially collocated interpolation/quadrature points \cite{black1999conservative, kopriva2009implementing}.
	These in principle arbitrarily high-order accurate schemes achieve entropy conservation/stability by computing the volume integral with \ac{SBP} operators and flux-differencing techniques \cite{gassner2016split, FISHER2013518, gassner2013skew, CHEN2017427, chan2018discretely}.
	In addition, to conserve entropy, inviscid surface flux functions \cite{Chandrashekar_2013, winters2016affordable, ranocha2018comparison, Ranocha2020Entropy} must be used.
	Entropy stability however, i.e., monotone decay of the mathematical entropy, is guaranteed for instance by standard dissipative Lax-Friedrichs type surface fluxes \cite{tadmor2003entropy}.

	In the previous paragraph we briefly outlined the building blocks of an \ac{EC/ES} \textit{spatial} discretization.
	To obtain a (fully) discrete \ac{EC/ES} scheme for an unsteady problem, it remains to discretize the \textit{temporal} derivative in a way that preserves this property.
	Examples of \ac{EC/ES} space-time methods can be found in \cite{hiltebrand2014entropy, friedrich2019entropy, gaburro2023high}.
	If the \ac{MoL} technique \cite{cockburn2001runge} is employed, different techniques to ensure the conservation of nonlinear functionals of the \ac{ODE} system (such as the mathematical entropy) need to be employed.
	Options for this are projection-based methods \cite[Chapter IV.4]{hairer2006geometric}, \cite{calvo2010projection, kojima2016invariants, najafian2025quasi}, the \ac{IDT} \cite{CalvoIDTPreservationInvariants}, or the relaxation methodology developed in \cite{RelaxationKetcheson, RelaxationRanocha}.
	As pointed out in \cite{kang2022entropy, ranocha2020fully, ranocha2020general}, a relaxation approach to ensure conservation of energy and other quadratic invariants was already proposed in \cite{sanz1982explicit, sanz1983method}.
	Therein, the authors consider a modified Leapfrog scheme and apply it to the Korteweg--de~Vries and nonlinear Schrödinger equation.
	The relaxation idea was also presented in \cite{dekker1984stability}.
	
	In this work we combine the relaxation methodology \cite{RelaxationKetcheson, RelaxationRanocha, ranocha2020general} and the \ac{PERK} methods \cite{vermeire2019paired, nasab2022third, doehring2024fourth} to obtain an optimized multirate time integration scheme that is \ac{EC/ES}, which we accordingly label \ac{PERRK}.
	The closest related work to the present study is \cite{kang2022entropy}, where the authors extend the second-order accurate multirate Runge-Kutta method from \cite{constantinescu2007multirate} with the relaxation methodology.
	Therein, however, the authors consider besides pure \acp{ODE} only the Burgers equation and do not quantify speed-up of the relaxed multirate methods over the baseline relaxed time integration methods.
	In contrast, we demonstrate in this work speed-up and reduced computational effort for up to fourth-order discretizations of the compressible Navier-Stokes-Fourier, \ac{vrMHD} and inviscid Euler equations.

	\vspace{0.5em}
	The paper is organized as follows:
	We begin by briefly sketching out the considered \ac{PDE} systems, \ac{PERK} methods, the notion of the (global) solution entropy, and the relaxation approach in \cref{sec:Preliminaries}.
	Then, we proceed to the construction of \ac{PERRK} methods in \cref{sec:PERRKs} which includes a detailed discussion on how the relaxation approach can act as a time-limiting technique to improve the robustness of the multirate scheme.
	The obtained methods are thoroughly validated in \cref{sec:Validation} where conservation of entropy and linear invariants is demonstrated without impairing the accuracy, i.e., order of consistency, of the schemes.
	In \cref{sec:Applications}, the \ac{PERRK} schemes are applied to a laminar, effectively incompressible flow around the SD7003 airfoil, a Mach $0.5$ flow of a visco-resistive magnetized fluid around a cylinder, an inviscid, transonic flows around the NACA0012 airfoil and ONERA M6 wing, and a viscous, transonic flow across the NASA \ac{CRM}.
	For each test case, we demonstrate improved performance of the \ac{PERRK} schemes over single-rate relaxed Runge-Kutta methods.
	Furthermore, we showcase increased robustness of the relaxed multirate \ac{PERRK} schemes over the standard \ac{PERK} methods, especially for under-resolved simulations.
	We conclude the paper in \cref{sec:Conclusions}.
	\section{Preliminaries}
	\label{sec:Preliminaries}
	%
	\subsection{Scope}
	In this work we consider \ac{PDE} systems of the form
	\begin{subequations}
		\label{eq:BalanceLaw}
		\begin{align}
			\boldsymbol u(t_0) &= \boldsymbol u_0 \\
			\label{eq:BalanceLawPDE}
			\partial_t \boldsymbol u(t, \boldsymbol x) + \nabla \cdot \boldsymbol f\big(\boldsymbol u, \nabla \boldsymbol u\big) 
			+ \boldsymbol g(\boldsymbol u, \nabla \boldsymbol u)
			&= 
			\boldsymbol s\big(t, \boldsymbol x; \boldsymbol u \big)
			, \quad t \in \big [t_0, t_f \big ], \: \boldsymbol x \in \Omega \subset \mathbb R^{d}
		\end{align}
	\end{subequations}
	where $\boldsymbol u \in \mathbb R^V$ denotes the vector of conserved variables and $\boldsymbol f \in \mathbb R^{V \times d}$ is a flux function which may contain diffusive terms.
	Nonconservative terms are collected in $\boldsymbol g$ and $\boldsymbol s$ is a source term.
	The classic theory of entropy conservation evolved around hypebolic conservation laws, i.e., system \eqref{eq:BalanceLaw} with $\boldsymbol f = \boldsymbol f(\boldsymbol u)$, $\boldsymbol g = \boldsymbol s = \boldsymbol 0$.
	Examples for entropy conservation problems with non-conservative terms, i.e., $\boldsymbol g \neq \boldsymbol 0$ are the shallow water equations \cite{GASSNER2016291} or the coupled Burgers equation \cite{castro2013entropy}.
	In the presence of viscous terms, i.e., $\boldsymbol f = \boldsymbol f(\boldsymbol u, \nabla \boldsymbol u)$ entropy is not conserved, but expected to decay monotonically/stable, assuming suitable boundary conditions.

	Here, we employ the \ac{DGSEM} \cite{hesthaven2007nodal, black1999conservative, kopriva2009implementing} to obtain the semidiscretization	of \eqref{eq:BalanceLaw}:
	\begin{subequations}
		\label{eq:Semidiscretization}
		\begin{align}
			\boldsymbol U(t_0) &= \boldsymbol U_0 \\
			\label{eq:SemidiscretizationODE}
			\boldsymbol U'(t) &= \boldsymbol F\big(t, \boldsymbol U(t) \big), \quad \boldsymbol U \in \mathbb R^M \: .
		\end{align}
	\end{subequations}
	The $M$ unknowns of \cref{eq:Semidiscretization} are the nodal \ac{DG} solution polynomial coefficients.
	Throughout this work we mainly consider convection-dominated problems, i.e., configurations for which the spectral radius $\rho_{\boldsymbol F}$ of the Jacobian of the fully discrete system \cref{eq:SemidiscretizationODE} $J_{\boldsymbol F}(\boldsymbol U) \coloneqq \pd{\boldsymbol F}{\boldsymbol U}$ scales as 
	\begin{equation}
		\label{eq:Convection_Diffusion_Scaling}
		\rho \sim \frac{\vert a \vert}{h} \gg \frac{\vert d \vert}{h^2} \: ,
	\end{equation}
	where $h$ denotes the smallest characteristic size of a grid cell, which is for the quadrilaterals and hexahedra considered in this work estimated based on the minimal distance between corner nodes.
	In \cref{eq:Convection_Diffusion_Scaling} $a, d$ are suitable quantifiers for the influence of convection and diffusion, respectively.
	%
	\subsection{Paired Explicit Runge-Kutta Methods}
	\label{subsec:PERKMs}
	\ac{PERK} methods have been developed in \cite{vermeire2019paired} and consecutively extended to third- \cite{nasab2022third} and fourth-order \cite{doehring2024fourth} of consistency.
	Furthermore, an embedded scheme based on the second-order method has been proposed in \cite{vermeire2023embedded} which enables error-based timestep control.
	
	The central idea of the \ac{PERK} methods is to pair Runge-Kutta schemes with different domains of stability and computational cost in an efficient way.
	These schemes belong to the class of \acp{PRKM} \cite{rentrop1985partitioned, bruder1988partitioned, gunther2001multirate} and are naturally applicable to semidiscretizations \eqref{eq:Semidiscretization} which may be partitioned according to different characteristic timescales:
	\begin{subequations}
		\label{eq:PartitionedODESys}
		\begin{align}
			\def\arraystretch{1.4}
			\boldsymbol U(t_0) &= \boldsymbol U_0, \\
			\label{eq:PartitionedODESys2}
			\boldsymbol U'(t) &= \begin{pmatrix} \boldsymbol U^{(1)}(t) \\ \vdots \\ \boldsymbol U^{(R)}(t) \end{pmatrix}' = \begin{pmatrix} \boldsymbol F^{(1)}\big(t, \boldsymbol U^{(1)}(t), \dots,  \boldsymbol U^{(R)}(t) \big) \\ \vdots \\ \boldsymbol F^{(R)}\big(t, \boldsymbol U^{(1)}(t), \dots, \boldsymbol U^{(R)}(t) \big) \end{pmatrix} = \boldsymbol F\big(t, \boldsymbol U(t) \big) \: .
		\end{align}
	\end{subequations}
	Here, $R$ denotes the number of \textit{partitions/levels}.
	Based on the local characteristic speed 
	\begin{equation}
		\label{eq:CharSpeed_DG}
		\nu(\boldsymbol x_m) \sim \max_{i = 1, \dots, d}\frac{(k+1) \cdot \rho_i(\boldsymbol x_m)}{h(\boldsymbol x_m)} > 0, \quad m = 1, \dots, M
	\end{equation}
	the unknowns $\left \{ U_m \right \}_{m=1, \dots, M}$ are assigned to different partitions $r = 1, \dots, R$ which are integrated with different schemes.
	In \eqref{eq:CharSpeed_DG} $d$ denotes the number of spatial dimensions, $h$ is a measure for the smallest characteristic cell size,
	$k$ is the \ac{DG} solution polynomial degree, and $\rho_i$ is the spectral radius of the Jacobian of the directional flux $\boldsymbol f_i$:
	\begin{equation}
		\label{eq:DirectionalFluxJacobian}
		J_i \coloneqq \pd{\boldsymbol f_i}{\boldsymbol u}
	\end{equation}
	cf. \eqref{eq:BalanceLawPDE}.

	In the following, the superscript $(\cdot)^{(r)}$ denotes quantities corresponding to the $r$'th partition.
	These may be a subset of unknowns $\boldsymbol U^{(r)} \in \mathbb R^{M_r}$, or a method used for solving the corresponding system with \ac{RHS} $\boldsymbol F^{(r)}$.
	Note that for the sake of readability we oftentimes pad the vectors $\boldsymbol U^{(r)}, \boldsymbol K^{(r)} \in \mathbb R^{M^{(r)}}$ with zeros to match the dimensionality $\mathbb R^M$ of the full system \eqref{eq:Semidiscretization}.
	Evaluations of the corresponding \ac{RHS} $\boldsymbol F^{(r)}$ are, however, only performed on the actual index set $M^{(r)}$ with property $\cup_{r=1}^R M^{(r)} = M$.
	\ac{ODE} systems of the form \eqref{eq:PartitionedODESys} encourage the use of \acp{PRKM}, which, in Butcher form, are given by \cite[Chapter~II.15]{HairerWanner1}
	\begin{subequations}
		\label{eq:PartitionedRKSystem}
		\begin{align}
			\label{eq:PartitionedRKFirstEq}
			\boldsymbol U_0 &= \boldsymbol U(t_0) \: , \\
			\label{eq:PartitionedRKSecondEq}
			\boldsymbol K_i^{(r) } &= \boldsymbol F^{(r) }\Biggl( \underbrace{t_n + c_i^{(r) } \Delta t}_{\eqqcolon t_{n,i}}, \underbrace{\boldsymbol U_n + \Delta t \sum_{j=1}^S \sum_{k=1}^R  a_{i,j}^{(k) } \boldsymbol K_j^{(k) }}_{\eqqcolon \boldsymbol U_{n,i}} \Biggr), \quad i = 1, \dots, S; \: \: r = 1, \dots, R \\
			\label{eq:PartitionedRKUpdateStep}
			\boldsymbol U_{n+1} & = \boldsymbol U_n + \Delta t  \sum_{i=1}^S  \sum_{r=1}^R   b_i^{(r) }  \boldsymbol K_i^{(r)} \: ,
		\end{align}
	\end{subequations}
	where $S$ denotes the number of stages $\boldsymbol K_i \in \mathbb R^M$.

\begin{figure}
		\centering
		\subfloat[{Schematic three-level ($R=3$) non-uniform mesh as obtained e.g. by quad/octree \ac{AMR}.
		The black dots resemble the Gauss-Lobatto-Legendre quadrature points of a $k=2$ solution polynomial.
		}]{
			\label{fig:ExampleMesh}
			\centering
			\includegraphics[width=.4\textwidth]{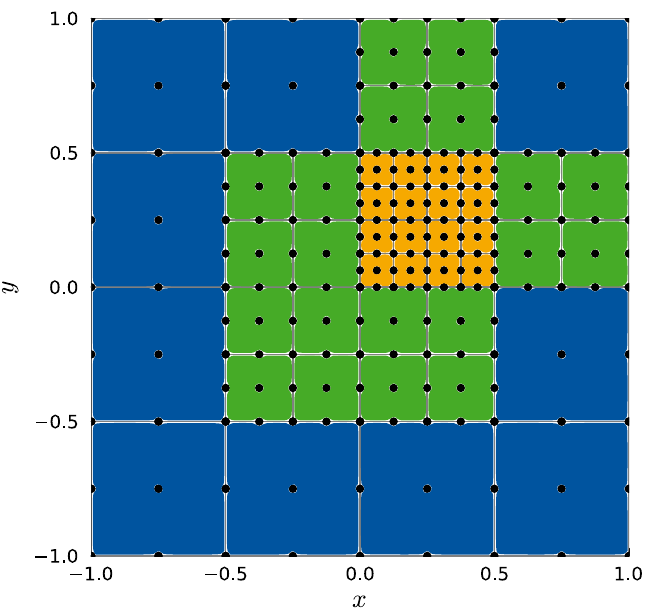}
		}
		\hfill
		\subfloat[{Domains of absolute stability of different-degree optimized fourth-order accurate stability polynomials.
		The spectrum corresponds to a \ac{DGSEM} discretization of the 2D advection equation with $\boldsymbol a = (1, 1)$ through $k=3$ solution polynomials and Rusanov/Local Lax-Friedrichs \cite{RUSANOV1962304} numerical flux.
		}]{
			\label{fig:StabDomains}
			\centering
			\resizebox{.5\textwidth}{!}{\includegraphics{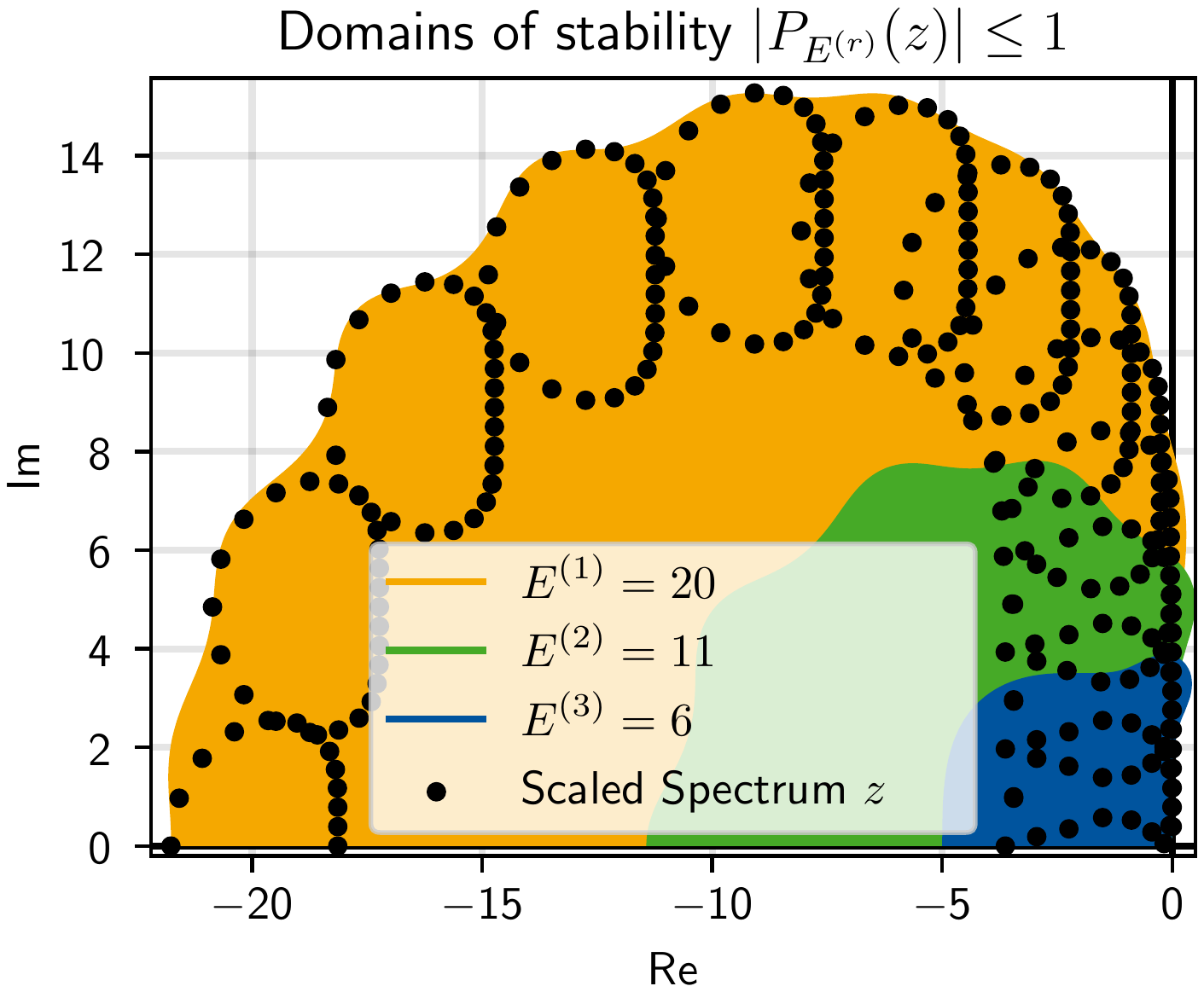}}
		}
		\caption[Illustration of the concept of multirate time integration with PERK methods.]
		{Illustration of the concept of multirate time integration with \ac{PERK} methods.
		In \cref{fig:ExampleMesh} a non-uniform mesh as obtained e.g. by classic quad/octree \ac{AMR} is shown.
		In \cref{fig:StabDomains} the domains of absolute stability of optimized stability polynomials are presented.
		These methods could be combined into a $R=3$ \ac{PERK} scheme to integrate a \ac{PDE} on meshes like the one shown in \cref{fig:ExampleMesh} efficiently.
		}
		\label{fig:ExampleMesh_StabDomains}
	\end{figure}

	When integrating \eqref{eq:Semidiscretization} with an explicit method, the maximum linearly stable timestep is restricted by some form of the \ac{CFL} condition:
	\begin{equation}
		\label{eq:TimestepConstraintGlobal}
		\Delta t \overset{!}{\leq} \text{CFL} \cdot \min_{\boldsymbol x \in \Omega} \frac{1}{\nu (\boldsymbol x)} \: .
	\end{equation}
	In contrast to local time-stepping methods \cite{berger1984adaptive, dawson2001high, lorcher2007discontinuous, dumbser2007arbitrary, krivodonova2010efficient, grote2015runge, luther2024adaptive},
	which reduce the timestep $\Delta t$ in regions with high characteristic speeds $\nu(\boldsymbol x)$ \eqref{eq:CharSpeed_DG}, \ac{PERK} schemes increase the $\text{CFL}$ stability bound locally to avoid reduction in the global timestep $\Delta t$:
	\begin{equation}
		\label{eq:TimestepConstraintLocal}
		\Delta t = \text{CFL}^{(r)} \cdot \frac{1}{\nu^{(r)}} \, , \quad \forall \: r = 1, \dots, R \: .
	\end{equation}
	This is achieved by optimizing the stability polynomials of the composing schemes $\left (A^{(r)}, b^{(r)} \right )$ to achieve larger $\text{CFL}^{(r)}$.
	In particular, equipped with the spectrum $\boldsymbol \sigma$ of the Jacobian of the fully-discrete system \eqref{eq:Semidiscretization} $J_{\boldsymbol F}$ the stability polynomials $P^{(r)}(z)$ can be optimized to yield the maximum timestep $\Delta t$ for this particular spectrum $\boldsymbol \sigma(J_{\boldsymbol F})$ \cite{ketcheson2013optimal, doehring2024manystage}.
	The key insight here is that the number of free coefficients of the $P^{(r)}(z)$ stability polynomial is equal to $E^{(r)} - p$, where $E^{(r)}$ denotes the number of stage \textit{evaluations} of the $r\text{-th}$ scheme and $p$ is the order of the method.
	Thus, computationally demanding methods with a higher number of stage evaluations $E^{(r)}$ \textit{and} larger domain of stability are applied in regions with higher characteristic speeds $\nu^{(r)}$, while less expensive methods with a smaller domain of stability are used in regions with lower characteristic speeds $\nu^{(r)}$.
	This concept is schematically illustrated in \cref{fig:ExampleMesh_StabDomains}.
	Standard, i.e., non-reducible \cite[Chapter~IV.12]{HairerWanner2} Runge-Kutta methods have $E = E^{(1)}= S$ stage evaluations, i.e., every stage present in e.g. the Butcher tableau is actually computed and utilized in the final update step or for the computation of subsequent stages.

	To foster the understanding of \ac{PERK} schemes, we present the Butcher tableaus of a second-order, $S=6$ stage, $E= \left \{ E^{(1)}, E^{(2)} \right \} =  \{6, 3\}$ stage evaluation \ac{PERK} method:
	\begin{equation}
			\label{eq:PERK_ButcherTableauClassic_p2_E36}
			\renewcommand\arraystretch{1.3}
			\begin{array}
					{c|c|c c c c c c|c c c c c c}
					i & \boldsymbol c & & & A^{(1)} & & &  & & & A^{(2)} & & &\\
					\hline
					1 & 0 & & & & & &             										 & & & & & &\\
					2 & \sfrac{1}{10} & \sfrac{1}{10} & & & & &        & \sfrac{1}{10} & & & & & \\
					3 & \sfrac{2}{10} & \sfrac{2}{10} - a_{3,2}^{(1)} & a_{3,2}^{(1)} & & & &      & \sfrac{2}{10} & 0 & & & &  \\ 
					4 & \sfrac{3}{10} & \sfrac{3}{10} -a_{4,3}^{(1)}& 0 & a_{4,3}^{(1)} & & &    & \sfrac{3}{10} & 0 & 0 & \\ 
					5 & \sfrac{4}{10} & \sfrac{4}{10} - a_{5,4}^{(1)} & 0 & 0 & a_{5,4}^{(1)} & & & \sfrac{4}{10} 								 & 0  & 0 & 0 & & \\
					6 & \sfrac{5}{10} & \sfrac{5}{10} - a_{6,4}^{(1)} & 0 & 0 & 0 & a_{6,4}^{(1)} & & \sfrac{5}{10} - a_{6,5}^{(2)}  & 0 & 0 & 0 & a_{6,5}^{(2)} & \\
					\hline
					& \boldsymbol b^T & 0 & 0 & 0 & 0 & 0 & 1          & 0 & 0 & 0 & 0             & 0 & 1
			\end{array}
	\end{equation}
	The first method $r= 1$ with Butcher matrix $A^{(1)}$ requires computation of all six stages $\boldsymbol K_i^{(1)}, \allowbreak i = 1, \dots, 6$ while the second method $r=2$ with $A^{(2)}$, however, requires only the computation of stages $\boldsymbol K_1^{(2)}, \allowbreak \boldsymbol K_5^{(2)}, \allowbreak \boldsymbol K_6^{(2)}$, i.e., three evaluations $E^{(2)}=3$ of $\boldsymbol F^{(2)}$.
	When looking at the second method in isolation, it resembles a reducible \cite[Chapter~IV.12]{HairerWanner2} method, i.e., the stages $\boldsymbol K_2^{(2)}, \allowbreak \boldsymbol K_3^{(2)}, \allowbreak \boldsymbol K_4^{(2)}$ do not influence the final solution $\boldsymbol U_{n+1}^{(2)}$ and the Butcher tableau could be truncated to a three-stage method.
	In the context of \acp{PRKM}, however, the second to fourth stage guarantee an internally consistent \cite{hundsdorfer2015error} update of the intermediate state $ \boldsymbol U_{n,i}$, cf. \eqref{eq:PartitionedRKSecondEq}.
	This is assured for identical abscissae $\boldsymbol c^{(r)} = \boldsymbol c \: \forall \: r = 1, \dots, R$ fulfilling the usual convention
	\begin{equation}
		\label{eq:InternalConsistency}
		\sum_{j=1}^{i-1} a^{(r)}_{i,j} \overset{!}{=} c_i \quad \forall \: i = 1, \dots , S; \: \: \forall \: r = 1, \dots , R \: .
	\end{equation}
	The intermediate, cheap update of the inactive partitions is essential for the overall accuracy of the method \cite{hundsdorfer2015error}.

	The particular form of the Butcher tableaus \eqref{eq:PERK_ButcherTableauClassic_p2_E36} comes with three benefits:
	\begin{enumerate}
		\item The \ac{PERK} methods are low-storage, requiring only storage of two Runge-Kutta stages $\boldsymbol K_i$ at the same time.
		\item The computational costs per stage are not increasing if higher stage (evaluation) methods are used, as always at most two stages are used for computation of the intermediate state $\boldsymbol U_{n,i}$.
		\item The sparse structure of both $A^{(r)}$ and $\boldsymbol b^{(r)}$ allows for a simplified computation of the coefficients $a_{i,j}^{(r)}$ from an optimized stability polynomial $P^{(r)}(z)$
	\end{enumerate}

	In previous works \cite{doehring2024fourth, doehring2024multirate} we discussed how results derived for partitioned Runge-Kutta methods (PRKMs) extend to \ac{PERK} methods.
	In particular, we show that criteria for linear stability, order of accuracy, internal consistency, and conservation of linear invariants are met by the \ac{PERK} methods proposed in \cite{vermeire2019paired, nasab2022third, doehring2024fourth}.
	The archetypes of these third- and fourth-order \ac{PERK} methods are given by \cref{eq:PERK_ButcherTableauOwn_p3} and \cref{eq:PERK_ButcherTableau_p4}, respectively.
	Thus, we omit this discussion here except that we remark that for a conservative \ac{PRKM} \eqref{eq:PartitionedRKSystem} the weight vector $\boldsymbol b^T$ must be identical across partitions \cite{hundsdorfer2015error} \cite[Chapter IV.1]{hairer2006geometric}, i.e., 
	\begin{equation}
		\label{eq:Conservation}
		b_i^{(r_1)}	 = b_i^{(r_2)} = b_i, \quad \forall \: i = 1, \dots, S; \: \: \forall \: r_1, r_2 = 1, \dots, R \: .
	\end{equation}
	Thus, in the following we only consider an identical weight vector $\boldsymbol b^T$.
	\subsection{The Notion of Entropy}
	In the context of this work, we are interested in \ac{ODE} systems \eqref{eq:Semidiscretization} with states $\boldsymbol U \in \mathbb R^M$ with standard Euclidian inner product $\langle \cdot, \cdot \rangle$ and norm.
	Additionally, we assume a smooth convex global entropy functional 
	\begin{equation}
		\label{eq:GlobalEntropyFunctional}	
		H: \mathbb R^M \rightarrow \mathbb R
	\end{equation}
	whose temporal evolution is tied to a physical meaningful solution of \eqref{eq:Semidiscretization} and thus also \eqref{eq:BalanceLaw}.
	By the chain rule, the evolution over time of $H(\boldsymbol U)$ is given by
	\begin{equation}
		\frac{\nid}{\nid t} H \big(\boldsymbol U(t) \big ) 
		= 
		\left \langle \frac{\partial H(\boldsymbol U)}{\partial \boldsymbol U}, \frac{\nid}{\nid t} \boldsymbol U(t) \right \rangle 
		\overset{\cref{eq:PartitionedODESys2}}{=}
		\left \langle \frac{\partial H(\boldsymbol U)}{\partial \boldsymbol U}, \boldsymbol F\big(t, \boldsymbol U(t) \big) \right \rangle \: .
	\end{equation}
	Entropy-conservative systems (e.g. appropriately discretized Euler equations) thus satisfy 
	\begin{equation}
		\label{eq:EntropyConservationContinuous}
		\forall \: \boldsymbol U \in \mathbb R^M, \: t \in [0, T]: \quad \left \langle \frac{\partial H(\boldsymbol U)}{\partial \boldsymbol U}, \boldsymbol F(t, \boldsymbol U) \right \rangle = 0
	\end{equation}
	with discrete equivalent
	\begin{equation}
		\label{eq:EntropyConservationDiscrete}
		H(\boldsymbol U_{n+1}) = H(\boldsymbol U_n) = H(\boldsymbol U_0)  \quad \: \forall \: n \in \mathbb N 
		\: .
	\end{equation}
	For entropy-dissipative systems (e.g. correspondingly discretized Navier-Stokes-Fourier equations) these equalities are replaced by inequalities, usually by the convention that the (mathematical) entropy decreases over time, i.e., 
	\begin{equation}
		\label{eq:EntropyDissipationDiscrete}
		\left \langle \frac{\partial H(\boldsymbol U)}{\partial \boldsymbol U}, \boldsymbol F(t, \boldsymbol U) \right \rangle \leq 0
		\quad \Rightarrow \quad 
		H(\boldsymbol U_{n+1}) \leq H(\boldsymbol U_n) \leq H(\boldsymbol U_0)
		\quad \: \forall \: n \in \mathbb N 
		\: .
	\end{equation}

	In this work we are concerned with the compressible flows of ideal gases described by the Euler, Navier-Stokes-Fourier, or \ac{MHD} equations.
	For these systems, the mathematical entropy $s$ as a function of the continuous conserved variables $\boldsymbol u$ is given by
	\begin{equation}
		\label{eq:EntropyMathematical}
		s(\boldsymbol u) \coloneqq - \underbrace{\rho}_{\equiv u_1} \cdot s_\text{therm}(\boldsymbol u) = - \rho \cdot \log \left( \frac{p(\boldsymbol u)}{\rho^\gamma} \right) \: ,
	\end{equation}
	which is then integrated over the spatial domain $\Omega$ to obtain the total entropy $\eta$:
	\begin{equation}
		\label{eq:TotalEntropyContinuosDefinition}
		\eta(t) \coloneqq \eta \big(\boldsymbol u(t, \boldsymbol x) \big)	= \int_{\Omega} s \big (\boldsymbol u(t, \boldsymbol x) \big ) \, \nid \boldsymbol x \: .
	\end{equation}
	Note that we denote with $\eta : \mathbb R \rightarrow \mathbb R$ the continuous equivalent of the discrete global entropy $H: \mathbb R^M \rightarrow \mathbb R$.
	\subsection{Relaxation Runge-Kutta Methods}
	\label{subsec:RelaxationRungeKuttaMethods}
	The term \ac{RRKM} was coined in \cite{RelaxationKetcheson} for Runge-Kutta methods whose final update step \eqref{eq:PartitionedRKUpdateStep} is scaled by the \textit{relaxation parameter} $\gamma_n \in \mathbb R$:
	\begin{equation}
		\label{eq:RRKM_Update}
		\boldsymbol U_{n+1}(\gamma_n) = \boldsymbol U_n + \Delta t \gamma_n \sum_{i=1}^Sb_i \boldsymbol K_i \: .
	\end{equation}
	Recall that $\boldsymbol K_i$ denotes the $i$-th stage of the Runge-Kutta method, i.e., 
	\begin{align}
		\label{eq:StageSingleEq}
		\boldsymbol K_i = \boldsymbol F \left( t_{n, i}, \boldsymbol U_{n,i} \right), \quad i = 1, \dots, S
	\end{align}
	with intermediate state $\boldsymbol U_{n,i} = \boldsymbol U_n + \Delta t \sum_{j=1}^S  a_{i,j} \boldsymbol K_j $ and time $t_{n,i} = t_n + c_i \Delta t$, cf. \cref{eq:PartitionedRKSecondEq}.

	The introduction of this relaxation parameter allows to satisfy exact entropy conservation, i.e., 
	\begin{equation}
		\label{eq:EntropyConservationRRKM}	
		H \big(\boldsymbol U_{n+1}(\gamma_n) \big) \overset{!}{=} H(\boldsymbol U_n) + \gamma_n \Delta t \sum_{i=1}^S b_i \left \langle \frac{\partial H(\boldsymbol U_{n, i})}{\partial \boldsymbol U_{n, i}}, \boldsymbol K_i \right \rangle \: .
	\end{equation}
	It is worthwhile to examine \eqref{eq:EntropyConservationRRKM} for a trivial function $\boldsymbol H$, for instance $\boldsymbol H (\boldsymbol U) = I \boldsymbol U$ where $I$ is the identity matrix $I \in \mathbb R^{M \times M}$.
	Thus, for a single component $k$ of $\boldsymbol U$ we have $H_k(\boldsymbol U) = U^{(k)} = \sum_{m=1}^M \delta_{km} U^{(m)}$ with Kronecker delta $\delta_{km}$.
	In this case \eqref{eq:EntropyConservationRRKM} corresponds to the well-known linear-invariant (i.e., mass) conservation property of Runge-Kutta methods, as for every component $k$ of $\boldsymbol U$ we have 
	\begin{align}
		H_k \big(\boldsymbol U_{n+1}(\gamma_n) \big) = U_{n+1}^{(k)}	(\gamma_n) &= U_n^{(k)} + \gamma_n \Delta t \sum_{i=1}^S b_i \left \langle \boldsymbol{e}_k, \boldsymbol K_i \right \rangle \\
		\overset{\eqref{eq:RRKM_Update}}{\Leftrightarrow} \sum_{m=1}^M \delta_{km}\left[ \boldsymbol U_n + \Delta t \gamma_n \sum_{i=1}^Sb_i \boldsymbol K_i	 \right] &= U_n^{(k)} + \gamma_n \Delta t \sum_{i=1}^S b_i K_i^{(k)}
	\end{align}
	where $\boldsymbol{e}_k \in \mathbb R^M$ is the $k$-th (Euclidian) unit vector.
	Additionally, we immediately see that arbitrary linear invariants are conserved for every choice of $\gamma_n$.
	The same holds of course also for linear functions resulting in a scalar entropy, such as $H(\boldsymbol U) = \boldsymbol q^T \boldsymbol U$ where $\boldsymbol q \in \mathbb R^M$ could resemble a vector corresponding to quadrature weights.

	For nonlinear functions $H(\boldsymbol U)$ such as the (mathematical) entropy \eqref{eq:EntropyConservationRRKM} does not trivially hold, however.
	In this scenario, we search for a relaxation parameter $\gamma_n$ such that the true change in entropy
	\begin{equation}
		\label{eq:EntropyChangeDefinition}
		\Delta H_n \coloneqq \Delta t \sum_{i=1}^S b_i \left \langle \frac{\partial H(\boldsymbol U_{n,i})}{\partial \boldsymbol U_{n,i}}, \boldsymbol K_i \right \rangle	
	\end{equation}
	matches the entropy of the relaxed update state $\boldsymbol U_{n+1}(\gamma_n)$:
	\begin{equation}
		\label{eq:EntropyConservationRRKM_Short}
		H \big(\boldsymbol U_{n+1}(\gamma_n) \big) = H \left( \boldsymbol U_n + \Delta t \gamma_n \sum_{i=1}^Sb_i \boldsymbol K_i  \right) \overset{!}{=} H(\boldsymbol U_n) + \gamma_n \Delta H_n (\boldsymbol U_n) \: .
	\end{equation}
	For entropy-conservative systems the the scalar product in \cref{eq:EntropyChangeDefinition} is thus exactly zero, i.e.,
	\begin{equation}
		\left \langle \frac{\partial H(\boldsymbol U_{n,i})}{\partial \boldsymbol U_{n,i}}, \boldsymbol K_i \right \rangle \overset{\cref{eq:StageSingleEq}}{=} \left \langle \frac{\partial H(\boldsymbol U_{n,i})}{\partial \boldsymbol U_{n,i}}, \boldsymbol F\big(t_{n,i}, \boldsymbol U_{n,i} \big) \right \rangle = 0
	\end{equation}
	Note that we demand the equality \cref{eq:EntropyConservationRRKM_Short} also for entropy diffusive systems, i.e., systems with $\langle \partial_{\boldsymbol U} H, \boldsymbol F \rangle \leq 0$ to suppress erronous additional changes in entropy due to the time integration method.
	
	We can rewrite \eqref{eq:EntropyChangeDefinition} by introducing the \textit{entropy variables} \cite{tadmor2003entropy}
	\begin{equation}
		\label{eq:EntropyVariables}
		\boldsymbol w(\boldsymbol u) \coloneqq \frac{\partial s(\boldsymbol u)}{\partial \boldsymbol u} \in \mathbb R^V
	\end{equation}
	with \eqref{eq:TotalEntropyContinuosDefinition} as
	\begin{equation}
		\label{eq:EntropyChange_EntropyVariables}
		\Delta H_n = \Delta t \sum_{i=1}^S b_i \left \langle \int_\Omega \boldsymbol W_{n,i} \, \nid \Omega \, , \boldsymbol K_i \right \rangle	
	\end{equation}
	which will be used in the remainder.
	The vector of discretized entropy variables $\boldsymbol W$ is obtained by evaluating the entropy variables at every node of the mesh and then re-assembling them into a vector of dimension $\mathbb R^M$.
	Note that integration of discrete integrands is denoted by the differential form $\nid \Omega$, while for integration of continuous quantities we employ $\nid \boldsymbol x$, cf. \eqref{eq:TotalEntropyContinuosDefinition}.

	It is well-known that the first order condition for a Runge-Kutta method is that $\sum_{i=1}^S b_i = 1$.
	Thus, we interpret $\gamma_n$ as acting on the timestep $\Delta t$ rather than on the weights $\boldsymbol b^T$.
	Consequently, the obtained update $\boldsymbol U_{n+1}(\gamma_n)$ is understood as an approximation at time $t_{n+1} = t_n + \gamma_n \Delta t$ \cite{RelaxationKetcheson, RelaxationRanocha}.

	Conversely, when interpreting the relaxed result $\boldsymbol U_{n+1}(\gamma_n)$ as an approximation to $t_{n+1} = t_n + \Delta t$ the \ac{RRKM} can be viewed as a projection method \cite[Chapter IV.4]{hairer2006geometric} \cite{RelaxationKetcheson}.
	In that case, the relaxation mechanism can be interpreted as an instance of the incremental direction technique (IDT) \cite{CalvoIDTPreservationInvariants, RelaxationKetcheson, RelaxationRanocha}.

	At this point, numerous questions arise, concerning for instance existence and uniqueness of $\gamma_n$, order of the method, and stability properties.
	Here, we briefly state the main results obtained in \cite{RelaxationKetcheson, RelaxationRanocha, ranocha2020general} which settle these concerns.
	\begin{itemize}
		\item \textbf{Existence and uniqueness of $\gamma_n$:} Clearly \eqref{eq:EntropyConservationRRKM_Short} is trivially fulfilled for $\gamma_n = 0$.
		A positive (i.e., nontrivial) unique solution which is close to unity exists for $\gamma_n$ for strictly convex entropy functionals $H$ coupled with at least second-order accurate schemes for sufficiently small $\Delta t$.
		\item \textbf{Accuracy:} Given a Runge-Kutta method of order $p$, the relaxed update $\boldsymbol U_{n+1}(\gamma_n)$ is an order $p$ approximation of the true solution at time $t_{n+1} = t_n + \gamma_n \Delta t$ for relaxation parameters $\gamma_n = 1 + \mathcal O \left(\Delta t^{p-1} \right)$ (which follows from the existence and uniqueness of $\gamma_n$).
		We validate this property in \cref{subsec:OrderOfConvergence}.
		\item \textbf{Linear Stability:} For $\gamma_n \leq 1$ the relaxed method is linearly stable for a given timestep $\Delta t$ if the standard method is linearly stable for this $\Delta t$.
		We provide a brief illustration on this manner by means of a concrete example in \cref{subsec:EntropyRelaxationTimeLimiting}.
		\item \textbf{Nonlinear Stability:} We remark that also e.g. \ac{SSP} methods \cite{shu1988efficient, gottlieb2001strong} can be relaxed while maintaining their nonlinear stability properties.
	\end{itemize}
	We are now equipped with the necessary tools to introduce the \ac{PERK} methods with relaxation.
	\section{Paired Explicit Relaxation Runge-Kutta Methods}
	\label{sec:PERRKs}
	The generalization of the relaxation methodology to \ac{PERK} schemes is straightforward.
	Due to the shared weight vector cf. \eqref{eq:Conservation} all stages $i$ with $b_i \neq 0$ are evaluated across all partitions $r = 1, \dots, R$.
	Thus, the \textit{direction}
	\begin{equation}
		\boldsymbol d_n \coloneqq \sum_{i=1}^S b_i \boldsymbol K_i
	\end{equation}
	and entropy change $\Delta H$ (cf. \eqref{eq:EntropyConservationRRKM_Short}) are computed completely analogous to standard Runge-Kutta methods.
	As we are only concerned with \ac{PERK} methods of order two and higher, the conditions for existence and uniqueness of a solution $\gamma_n$ are fulfilled.
	Furthermore, the \ac{PERK} schemes considered in \cite{vermeire2019paired, nasab2022third, doehring2024fourth, doehring2024multirate} have nonnegative weights $b_i \geq 0$ which ensures that the relaxed update is dissipation preserving for such problems \cite{RelaxationRanocha} and avoids stage downwinding \cite{gottlieb1998total}.

	The sparsity in the weight vector $\boldsymbol b^T$ renders \ac{PERK} methods especially attractive for the relaxation methodology.
	The vanilla second-order method as proposed in \cite{vermeire2019paired} has only $b_S = 1$ with all other entries being zero.
	The third-order method proposed in \cite{nasab2022third} has non-zero $b_{S-1} = 0.75, \allowbreak b_S = 0.25$.
	The third-order variant proposed in \cite{doehring2024multirate} adds $b_1 = \sfrac{1}{6}$ to $b_{S-1} = \sfrac{1}{6}, \allowbreak b_S = \sfrac{2}{3}$ which avoids a negative entry in Butcher array of the original scheme \cite{nasab2022third}, thus there is no downwinding of stages \cite{gottlieb1998total}.
	Finally, the fourth-order method developed in \cite{doehring2024fourth} has $b_{S-1} = b_S = 0.5$.
	Archetypes for the Butcher tableaus of these methods are provided in \ref{sec:P-EKR_ButcherTableauArchetypes}.
	
	This sparsity in $\boldsymbol b^T$ is a very attractive feature as this reduces the number of integrations involved in a sequential computation of $\Delta H$, cf. \eqref{eq:EntropyChangeDefinition}.
	Of course, with the usual computation-storage tradeoff one can always reduce \eqref{eq:EntropyChangeDefinition} to a single integration if all intermediate approximations $\boldsymbol U_{n, i}$ and stages $\boldsymbol K_i$ are stored, even though they are not required for the actual update \eqref{eq:PartitionedRKUpdateStep}.
	We remark that the transform from conserved $\boldsymbol U_{n,i}$ to entropy variables $\boldsymbol W_{n,i}$ would also be required for this approach.
	In our implementation of the \ac{PERRK} methods we opt for the sequential computation of \eqref{eq:EntropyChangeDefinition} to avoid unnecessary memory overhead.

	In practice, we compute \eqref{eq:EntropyChange_EntropyVariables} as
	\begin{equation}
		\label{eq:EntropyChange_EntropyVariables_ReArranged}
		\Delta H_n = \Delta t \sum_{i=1}^S b_i \sum_j \int_{\Omega_j} \sum_{l} \big \langle \boldsymbol w_{n, i, j, l}, \boldsymbol k_{i, j, l} \big \rangle	\, \nid \Omega_j 
	\end{equation}
	which follows from interchanging the scalar product and the integral (a quadrature rule on the discrete level).
	Thus, for every element $\Omega_j$ in the mesh we compute the dot product $\big \langle \boldsymbol w_{n, i, j, l}, \boldsymbol k_{i, j, l} \big \rangle$ for each node $l$ in the element.
	Note that on every node we have exactly $V$ variables, thus the inner product in \cref{eq:EntropyChange_EntropyVariables_ReArranged} is a scalar product on $\mathbb R^V$.
	\subsection{Solving the Relaxation Equation}
	To ensure an entropy-stable time update the scalar nonlinear equation 
	\begin{equation}
		\label{eq:RelaxationEquation}
		r(\gamma_n) \coloneqq H \left( \boldsymbol U_n + \gamma_n \Delta t \boldsymbol d_n\right) - H(\boldsymbol U_n) - \gamma_n \Delta H_n \overset{!}{=} 0
	\end{equation}
	needs to be solved for $\gamma_n$.
	In the preparation of this work we considered three main approaches to solve \eqref{eq:RelaxationEquation}.
	These are a standard bisection procedure which is quite robust and easy to implement.
	The drawback of this method is, however, that one needs to supply a bracketing interval $[\gamma_{\text{min}}, \gamma_{\text{max}}]$.
	As the bisection method converges only linearly, the choice of this bracketing interval heavily influences the number of iterations taken to solve \eqref{eq:RelaxationEquation} to some desired accuracy.
	Second, we considered the secant method which can be seen as a middle-ground between the bisection method and Newton-Raphson method.
	In particular, the secant method is of similar computational complexity as the bisection method, but offers a higher convergence rate.
	As for bisection, however, the secant method also requires two initial guesses $\gamma_0, \gamma_1$ which are used to compute the first update.
	The presence and significance of these hyperparameters render the secant method only a fallback option for our purposes.

	Thus, we mainly considered variants of the Newton-Raphson method which converge quadratically near the root.
	The required derivative of $\frac{\nid r}{\nid \gamma_n}$ is given by \cite{RelaxationRanocha}
	\begin{equation}
		r'(\gamma_n) = \int_\Omega \big \langle \boldsymbol w(\boldsymbol U_n + \gamma_n \Delta t \boldsymbol d_n), \Delta t \boldsymbol d_n \big \rangle \, \nid \Omega - \Delta H_n \: .
	\end{equation}
	To determine the stepsize in the Newton-Raphson procedure we experimented with different line-search methods.
	In practice (see \cref{sec:EntropyConservation}) we found that the vanilla Newton-Raphson procedure with trivial stepsize $1$ converges extremely fast ($1$ to $3$ iterations) and additional line-search methods are not necessary.
	In the very first timestep we set $\gamma_1 = 1$ while in subsequent timesteps we use the previously obtained value for the relaxation parameter as the initial guess for the Newton-Raphson iteration.
	As stopping criteria we employ the residual of the relaxation equation \eqref{eq:RelaxationEquation}, the magnitude of the update step and a maximum number of iterations to prevent iterating divergent sequences.
	\subsection{Entropy Relaxation as Time Limiting}
	\label{subsec:EntropyRelaxationTimeLimiting}
	The traditional discretization of a hyperbolic \ac{PDE} involves a spatial discretization that is \ac{TVD}/\ac{TVB} or offers at least some form of control over oscillations at discontinuities.
	This is then paired with an \ac{SSP} time integration scheme to prevent the introduction of spurious oscillations due to the time marching mechanism \cite{gottlieb2001strong, gottlieb1998total}.
	The latter scheme may be explicit or implicit, however, even for unconditionally linearly stable (A-stable) implicit methods of order $p > 1$, the nonlinearly stable, i.e., stability-preserving timestep is finite \cite{gottlieb2001strong}.
	In order to render implicit schemes more efficient than their explicit counterparts timesteps beyond the \ac{SSP} stability limit are often employed.
	Since in this case global non-oscillarity can no longer be guaranteed, localized \textit{time-limiting} strategies have been studied \cite{forth2001second, duraisamy2003concepts, duraisamy2007implicit} with some recent contributions \cite{arbogast2020third, puppo2023quinpi, puppo2023quinpi2}.
	The central idea in the referenced works is to blend a high-order implicit scheme with the first-order, unconditionally nonlinearly stable implicit/backward Euler method.

	In order to determine the limiting factor/blending coefficient different criteria based e.g. on temporal derivatives \cite{duraisamy2003concepts, puppo2023quinpi} have been employed.
	Recently, an indicator based on the local numerical entropy production 
	\begin{equation}
		\label{eq:EntropyProduction_CellLocal_Quinpi2}
		\Delta s_{n, j} \coloneqq 
		\frac{1}{\Delta t} \left[ 
			s \left( \boldsymbol U_{n+1, j} \right) - s \left( \boldsymbol U_{n, j} \right) - 
			\frac{\Delta t}{\Delta x} 
			\sum_{i=1}^S b_i \left( \Psi_{i, j-\frac{1}{2}}  - \Psi_{i, j+\frac{1}{2}} \right)
		\right]
	\end{equation}
	was proposed in \cite{puppo2023quinpi2}.
	Derived in the context of one-dimensional finite-volume schemes, \eqref{eq:EntropyProduction_CellLocal_Quinpi2} denotes the entropy production on cell $j$ due to numerical entropy fluxes $\Psi_{i, j \pm \frac{1}{2}}$ exchanged from time $t_n$ to $t_{n+1}$.
	The numerical entropy fluxes are required to be consistent with the entropy flux $\psi(\boldsymbol u)$ \cite{tadmor2003entropy} which corresponds to the entropy $s(\boldsymbol u)$ in the sense that $\left(\frac{\partial \psi}{\partial \boldsymbol u}\right)^T = 
	\boldsymbol w^T
	\frac{\partial \boldsymbol f}{\partial \boldsymbol u}$ \cite{tadmor2003entropy}.
	Equation \eqref{eq:EntropyProduction_CellLocal_Quinpi2} can be understood as a localized version of the relaxation equation \eqref{eq:EntropyConservationRRKM} as both determine the change in entropy due to the numerical scheme.
	In fact, by integrating \eqref{eq:EntropyProduction_CellLocal_Quinpi2} over all cells and time interval $[t_n, t_{n+1}]$ we recover \eqref{eq:EntropyConservationRRKM} for a suitable choice of the numerical entropy flux function $\Psi$.
	
	We would like to remark that the relaxation methodology can also be applied to local entropies, i.e., convex functionals that are defined on subdomains $\Omega_j \subset \Omega$ \cite{ranocha2020fully}.
	In the extreme, these subdomains could be reduced to single cells, resulting in local entropy changes $\Delta H_j$ which can be computed analogous to \cref{eq:EntropyProduction_CellLocal_Quinpi2}.

	In \cite{puppo2023quinpi2}, the authors employ \eqref{eq:EntropyProduction_CellLocal_Quinpi2} as an indicator if the timestep needs to be reduced to avoid oscillations.
	Similarly, the relaxation mechanism as discussed above can act as a limiter for monotonicity violations due to discontinuities in the partitioned Runge-Kutta scheme.
	This is particularly attractive as we have demonstrated in \cite{doehring2024multirate} that at the intersection of partitions spurious oscillations may occur.
	These are intensified for larger differences in active stage evaluations $E^{(1)} - E^{(2)}$, larger timesteps, and reduced grid resolution \cite{doehring2024multirate}.
	In the next section we follow up on this insights and demonstrate in what sense the relaxation methodology can act as a time-limiting mechanism for \ac{PERK} methods.
	\subsubsection{1D Advection Equation}
	\label{subsec:EntropyRelaxationTimeLimiting_1DAdvection}
\begin{figure}
		\centering
		\resizebox{.75\textwidth}{!}{\includegraphics{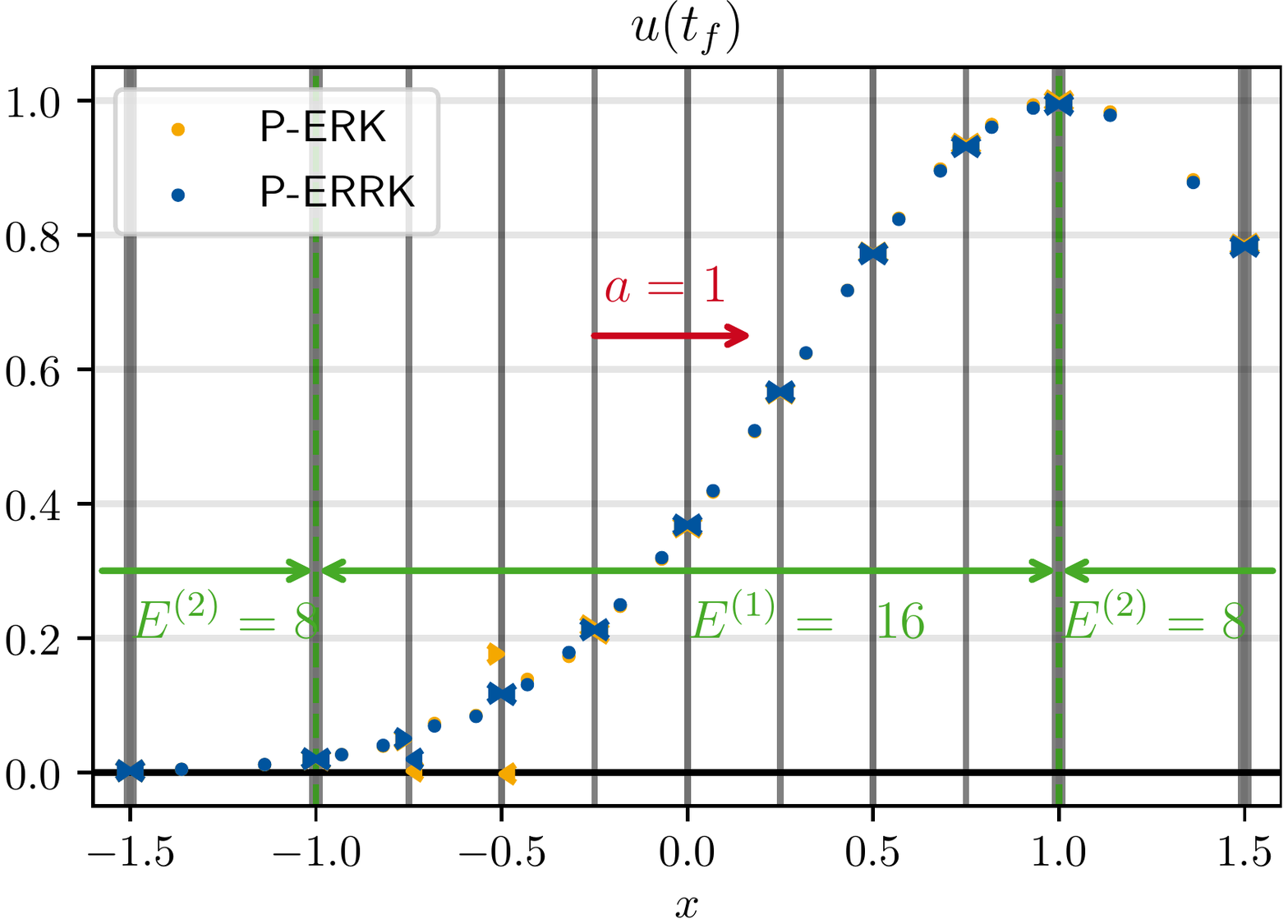}}
		\caption[Computational mesh for the VRMHD flow past a cylinder.]
		{Relaxation as time-limiting: Advection equation \cref{eq:AdvectionEq_Plain} on two-level, non-uniform mesh solved with standard (orange) and relaxed (blue) \ac{PERK} methods.
		The interpolated solution is displayed at the Gauss-Lobatto-Legendre nodes values at $t_f = 9$ on the interval 
		$[-1.6, 1.6]$.
		Cell boundaries are indicated by the vertical lines and contain $k+1 = 4$ nodes.
		To distinguish the nodes at cell boundaries, left boundary nodes are indicated by left-pointing triangles, while right boundary nodes are indicated by right-pointing triangles.
		Interior nodes are indicated by circles.}
		\label{fig:TimeLimiting_Solution}
	\end{figure}

	We consider the standard advection equation 
	\begin{equation}
		\label{eq:AdvectionEq_Plain}
		u_t + u_x = 0
	\end{equation}
	on a non-uniform grid with periodic boundaries.
	As a first motivation we present the advection of a Gaussian pulse $u_0(x) = \exp \left(-x^2 \right)$ on $\Omega = [-4, 4]$.
	The base mesh size is $\Delta x^{(2)} = 0.5$ while cells in $[-1, 1]$ are refined once, i.e., have length $\Delta x^{(1)} = 0.25$.
	We optimize $E = \{8, 16\}$ second-order accurate methods for a $k=3$ \ac{DGSEM} discretization with Upwind/Godunov surface flux.
	We run the simulation with $\Delta t = 0.2$ until $t_f = 9$, i.e., the Gaussian pulse has traversed the domain once and is now centered at $x=1$, with left tail in the refined region.
	For the standard \ac{PERK} method we observe severe oscillations in the second and third refined cell, \cref{fig:TimeLimiting_Solution}.
	In fact, these oscillations are strong enough to violate the positivity of the solution.
	Next, we consider the \ac{PERRK} method.
	For the standard scalar advection equation \cref{eq:AdvectionEq_Plain} a suitable mathematical entropy is given by $s(u) = u^2$ which coincides with the energy of the solution.
	We solve the relaxation equation \eqref{eq:RelaxationEquation} using the Newton-Raphson method which is limited to 10 iterations and also stops if the function residual is less or equal $10^{-14}$ or the update step is less than $10^{-15}$.
	If the Newton iteration does not converge to a root, a non-relaxed step is taken corresponding to $\gamma = 1$.
	This is, however, not observed for this simulation.
	Furthermore, on average only three Newton iterations per simulation step are performed.

	It is evident from \cref{fig:TimeLimiting_Solution} that the relaxation approach greatly reduces oscillations in the second cell in direction of the advection velocity at the coarse-fine intersection.
	In order to better visualize Gauss-Lobatto-Legendre nodes at the cell interfaces, we employ right-pointing triangles for the right boundary nodes and left-pointing triangles for the left boundary nodes, while interior nodes are indicated by circles.
	By closely examining \cref{fig:TimeLimiting_Solution} we observe that the entropy-stable relaxed scheme manages to preserve positivity of the solution, while the standard \ac{PERK} method shows negative nodal values.
	This comes at the expense, however, that we observe slight global diffusion of the solution, which is visible at the tip of the Gaussian pulse.

\begin{figure}
		\centering
		\hfill
		\subfloat[{Entropy evolution for standard (orange) and relaxed (blue) two-level multirate methods.}]{
			\label{fig:TimeLimiting_Entropy}
			\centering
			\resizebox{.47\textwidth}{!}{\includegraphics{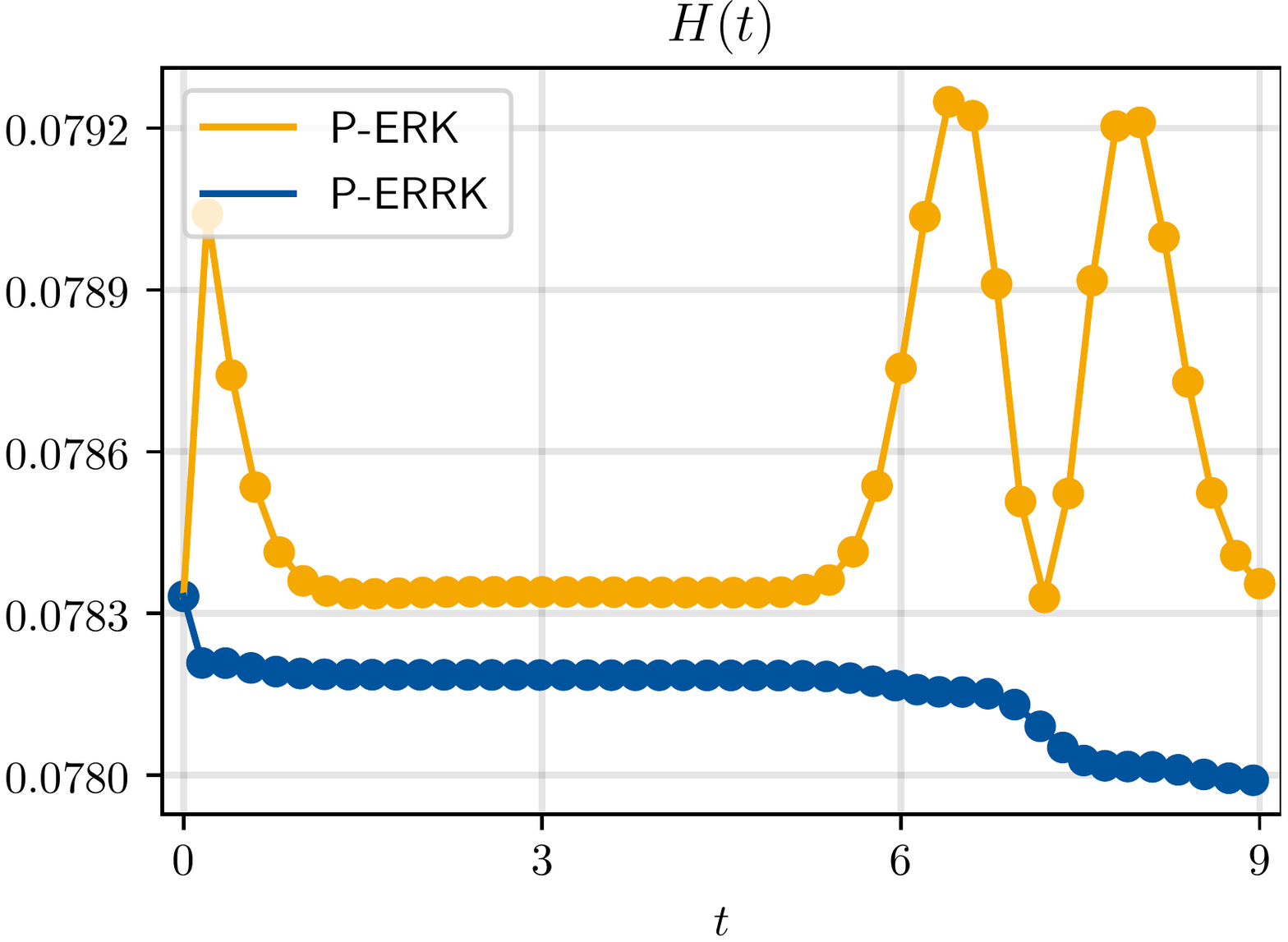}}
		}
		\hfill
		\subfloat[{Relaxation parameter $\gamma$ over the course of the simulation.}]{
			\label{fig:TimeLimiting_Gamma}
			\centering
			\resizebox{.47\textwidth}{!}{\includegraphics{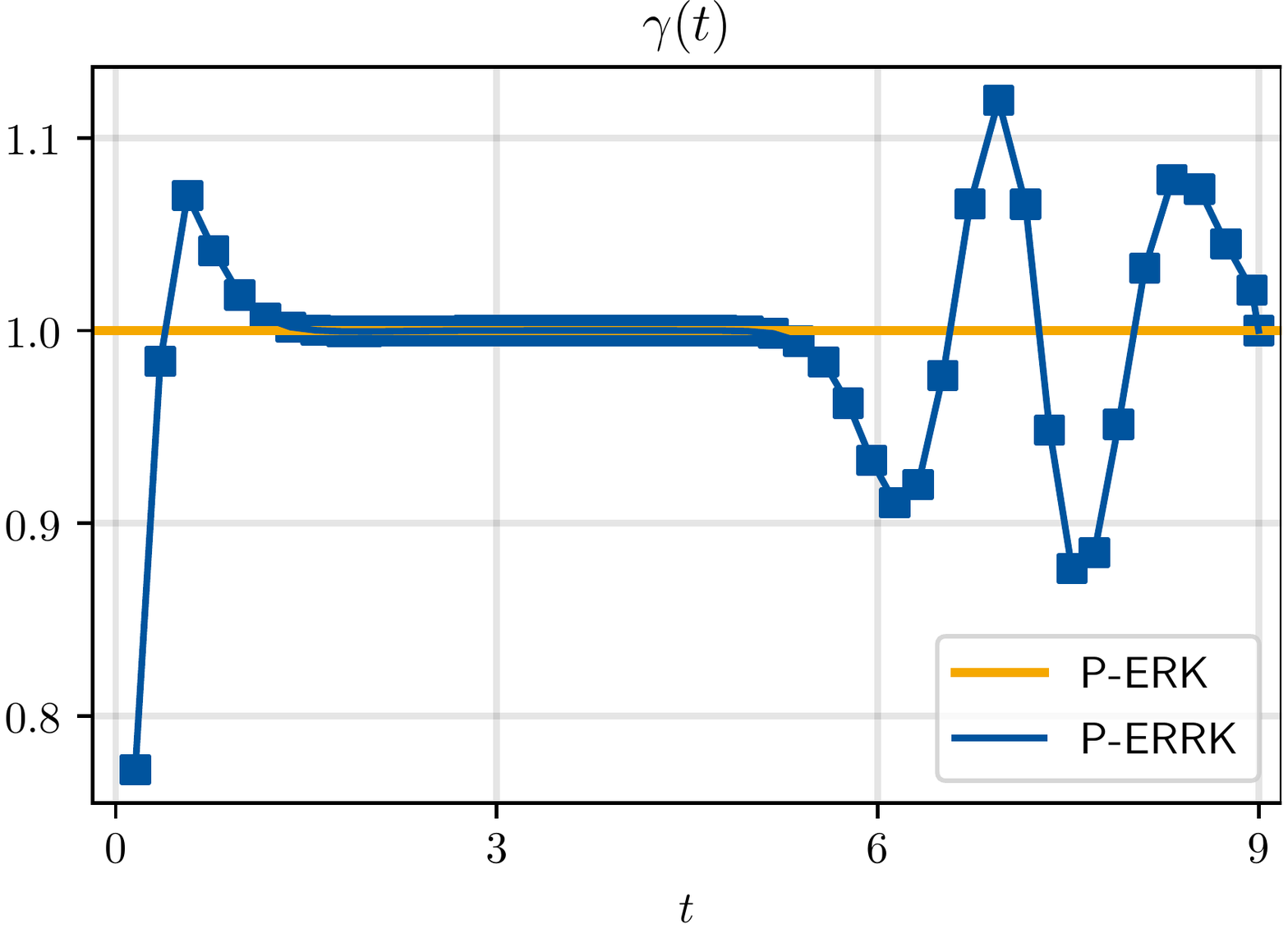}}
		}
		\caption[Relaxation as time-limiting: Entropy and relaxation parameter evolution.]
		{Relaxation as time-limiting: Advection equation \cref{eq:AdvectionEq_Plain} on two-level, non-uniform mesh.
		The temporal evolution of the entropy $H(t)$ is depicted in \cref{fig:TimeLimiting_Entropy} for the standard (orange) and relaxed (blue) \ac{PERK} method.
		The relaxation parameter $\gamma(t)$ is shown in \cref{fig:TimeLimiting_Gamma} for the relaxed method.}
		\label{fig:TimeLimiting_Entropy_Gamma}
	\end{figure}

	By examining the temporal evolution of the global entropy $H(t)$ \cref{fig:TimeLimiting_Entropy} for the two methods it is evident that the relaxation methodology adds diffusion only when required, i.e., when the entropy for the standard \ac{PERK} scheme increases.
	In particular, these are the first timestep and the traversal of the refined region by the Gaussian pulse happening roughly in the $t \in [5, 9]$ time interval.
	This is also reflected in the evolution of the relaxation parameter $\gamma_n$ \cref{fig:TimeLimiting_Gamma} which deviates from $1$ only when the entropy of the standard \ac{PERK} method changes.
	Noticeably, $\gamma$ is negatively correlated with the change in entropy $H$, i.e., for increasing entropy $\gamma_n$ is reduced and vice versa, which is related to the convexity of the entropy functional.
	To some extend, the relaxation methodology acts in this case also as a timestep adapting mechanism which is based on a nonlinear solution functional instead of the deviation between a higher- and a lower-order embedded method.
	The relaxation methodology thus ensures that the entropy dissipative nature of the spatial discretization is preserved.
	In particular, we cannot expect an entropy conservation scheme here as we employ the diffusive Upwind/Godunov flux and not e.g. the central flux.

	To shed further light on the stabilizing effect due to relaxation we consider the fully discrete system
	\begin{equation}
		\label{eq:FullyDiscreteSystem}	
		\def\arraystretch{1.3}
		\boldsymbol U_{n+1} = \begin{pmatrix} \boldsymbol U^{(1)}_{n+1} \\ \boldsymbol U^{(2)}_{n+1} \end{pmatrix} = \begin{pmatrix} D^{(1)}(\gamma_{n+1}) \\ D^{(2)}(\gamma_{n+1}) \end{pmatrix}
	\boldsymbol U_{n} = D(\gamma_{n+1}) \, \boldsymbol U_n
	\end{equation}
	which is obtained by applying the $p=2, E = \{8, 16 \}$ two-level \ac{PERK} method \cref{eq:PartitionedRKSystem} to 
	\begin{equation}
		\label{eq:LinearTwoLevelPartitionedODESys}
		\boldsymbol U'(t) = \begin{pmatrix} \boldsymbol U^{(1)}(t) \\ \boldsymbol U^{(2)}(t) \end{pmatrix}' = \begin{pmatrix} L^{(1)} \\ L^{(2)} \end{pmatrix} \boldsymbol U(t) = L \, \boldsymbol U(t)
	\end{equation}
	which is the linear equivalent of \cref{eq:PartitionedODESys2}.
	The update matrix $D(\gamma)$ can now be used to confirm linear stability by assessing whether all eigenvalues of $D(\gamma)$ are inside the unit circle.
	This is indeed the case for both the standard \ac{PERK} method, i.e., $\gamma = 1$ and the relaxed method with $\gamma \neq 1$, which confirms the discussion on linear stability from \cref{subsec:RelaxationRungeKuttaMethods}.
	Note that for the standard \ac{PERK} method $D(\gamma) \equiv D(1)$ for all times, i.e., the update matrix remains constant.
	For the \ac{PERRK} methods, however, $D(\gamma_n)$ is in general not identical among timesteps due to the relaxation parameter $\gamma_n$.
	We provide plots of the spectra of the constant $D(1)$ and $D(\gamma_1)$ matrices in the appendix, see \cref{fig:LinearStability}.
	Furthermore, one may employ $D$ to examine the nonlinear stability of the time integration methods.
	It is well-known that the fully-discrete method is monotonicity-preserving in the sense of Harten \cite{HARTEN1997260} if and only if
	\begin{equation}
		\label{eq:MonotonicityRequirement}
		D_{ij} \geq 0 \quad \: \forall \: i, j = 1, \dots, N \: .
	\end{equation}
	Since the underlying spatial discretization is a fourth-order \ac{DGSEM} scheme we cannot expect \eqref{eq:MonotonicityRequirement} to hold.
	Nevertheless, it is interesting to examine the minimum and maximum entries of $D(1)$ and $D(\gamma_1)$ to see if the relaxation mechanism has a stabilizing effect.
	For the standard \ac{PERK} method we have $\min_{i,j} D_{i,j}(1) \approx -0.370$ and $\max_{i,j} D_{i,j}(1) \approx 1.115$.
	For the relaxed method we have for $\gamma_1 \approx 0.772$ that $\min_{i,j} D_{i,j}(\gamma_1) \approx -0.286$ and $\max_{i,j} D_{i,j}(\gamma_1) \approx 0.86$.
	Clearly, the relaxation mechanism nudges the entries of $D$ closer to the $[0, 1]$ interval.
	Additionally, we would like to remark that for both cases the row-sums equal one, i.e., $\sum_j D_{ij} = \boldsymbol 1$, thereby reflecting the conservation property of the \ac{PERK}/\ac{PERRK} methods. 
	This is discussed in more detail in \cref{subsec:ConservationOfLinearInvariants}.
	\subsubsection{Kelvin-Helmholtz Instability with Adaptive Mesh Refinement}
	\label{subsec:KelvinHelmholtzInstability}
	Motivated by these insights we reconsider a numerical example conducted in \cite{doehring2024multirate}, wherein we studied the applicability of the standard \ac{PERK} schemes to the Kelvin-Helmholtz instability on an adaptively refined mesh.
	We employ the setup from \cite{rueda2021subcell} with initial condition at $t_0 = 0$ given by
	\begin{equation*}
		\begin{pmatrix}
			\rho \\ v_x \\ v_y \\ p
		\end{pmatrix}
		= 
		\begin{pmatrix}
			\frac{1}{2} + \frac{3}{4} b \\
			\frac{1}{2} (b - 1) \\
			\frac{1}{10} \sin(2 \pi x) \\
			1
		\end{pmatrix}, 
		\quad b \coloneqq \tanh \left( 15y + 7.5 \right) - \tanh \left(15y - 7.5 \right)
	\end{equation*}
	on periodic domain $\Omega = [-1, 1]^2$ which we simulate up to $t_f = 3.2$.
	The compressible Euler equations are discretized using the \ac{DGSEM} with $k=3$ local polynomials, \ac{HLLE} surface-flux \cite{einfeldt1988godunov}, subcell shock-capturing \cite{hennemann2021provably} to avoid spurious oscillations, and entropy-conserving volume flux \cite{Ranocha2020Entropy}.
	Initially, the domain is discretized with $2^4 \times 2^4$ cells which are allowed to be refined five times, i.e., the finest cells may have edge length $h = 2^{-8}$.
	The indicator for \ac{AMR} coincides in this case with the shock-capturing indicator \cite{hennemann2021provably} and is based on the density $\rho$.
	In \cite{doehring2024multirate} the $p=3, \allowbreak E = \{4, \allowbreak 6, \allowbreak 11\}$ scheme was constructed which we employ here as well.
	The only parameter changed compared to \cite{doehring2024multirate} is the \ac{AMR} interval which is decreased to $N_\text{AMR} = 7$, i.e., the mesh is allowed to change every $7^\text{th}$ timestep.
	The relaxation equation is solved with the Newton-Raphson method with a maximum of $5$ iterations, root-residual tolerance of $10^{-15}$ and stepsize tolerance of $10^{-13}$.

	We observe significantly improved robustness of the relaxed \ac{PERRK} scheme compared to the standard \ac{PERK} scheme.
	The relaxed scheme advances the simulation stable up to the final time $t_f = 3.2$ (see \cref{fig:DensityMesh_KelvinHelmholtzInstability} for a plot of the density and mesh at final time) while the standard scheme crashes at $t \approx 2.74$ due to violated positivity of density and pressure.
	We emphasize that application of the relaxation methodology does not exclusively reduce the timestep, i.e., the relaxation parameter $\gamma_n$ is not always less than $1$.
	In particular, for the simulation up to $t_f = 2.73$ we observe in fact reduced number of timesteps for the relaxed scheme ($1184$) compared to the standard \ac{PERK} scheme ($1267$).
	Although very encouraging, this should be taken with a grain of salt due to the highly volatile nature of the Kelvin-Helmholtz instability, especially when paired with \ac{AMR}.
	\section{Validation}
	\label{sec:Validation}
	\subsection{Entropy Conservation}
	\label{sec:EntropyConservation}
	To check for the exact entropy conservation property provided by the \ac{PERRK} methods we consider the compressible Euler equations in conjunction with an entropy-conservative spatial discretization.
	Here, we always employ the \ac{DGSEM} \cite{hesthaven2007nodal, black1999conservative, kopriva2009implementing} with entropy-conservative and kinetic energy preserving numerical flux by Ranocha \cite{Ranocha2020Entropy} to solve the Riemann problems at cell interfaces.
	Additionally, the volume integral is computed based on \ac{SBP} operators and flux-differencing \cite{gassner2016split, FISHER2013518, gassner2013skew, CHEN2017427, chan2018discretely} which is required to stabilize the entropy-conservative, i.e., diffusion-free surface flux.
	For the "volume flux" we use in the examples below also the flux by Ranocha, although we would like to remark that other choices such as the flux by Chandrashekar \cite{Chandrashekar_2013} are possible \cite{ranocha2018comparison}.
	All computations presented in this work are performed with Trixi.jl \cite{trixi1, trixi2, trixi3}, an open-source Julia \cite{bezanson2017julia} package for high-order, adaptive simulation of conservation laws and related systems.
	We provide a reproducibility repository publicly available on Github \cite{doehring2025PERRK_ReproRepo} for all results discussed in this paper.
	\subsubsection{Weak Blast Wave: Euler \& Ideal MHD Equations}
	\label{subsec:EC_WeakBlastWaveEulerMHD}
	We demonstrate the entropy conservation capabilities of the \ac{PERRK} methods for a discontinuous initial condition.
	To this end, we consider a 1D variant of the weak blast wave problem as proposed in \cite{hennemann2021provably} with discontinuous initial condition
	\begin{equation}
		\label{eq:WeakBlastWave_IC}
			\rho_0 = \begin{cases}
				1.1691 & \vert x \vert \leq 0.5 \\
				1 & \vert x \vert > 0.5 \\
			\end{cases},
			\:
			v_0 = \begin{cases}
				 0.1882 \cdot \text{sgn}(x) & \vert x \vert \leq 0.5 \\
				0 & \vert x \vert > 0.5 \\
			\end{cases},
			\:
			p_0 = \begin{cases}
				1.245 & \vert x \vert \leq 0.5 \\
				1 & \vert x \vert > 0.5 \\
			\end{cases}
			\: .
	\end{equation}
	The isentropic exponent is set to $\gamma = 1.4$ and we employ $k = 3$ \ac{DG} solution polynomials.
	The domain is $\Omega = [-2, 2]$ equipped with periodic boundary conditions.
	We discretize $\Omega$ with a non-uniform mesh consisting of three-levels which is refined towards the center.
	Precisely, cells in the $[-0.5, 0.5]$ interval are of size $\Delta x^{(1)} = \sfrac{1}{32}$, cells in the $[-1, -0.5] \cup [0.5, 1]$ interval are of size $\Delta x^{(2)} = \sfrac{1}{16}$, and cells in the $[-2, -1] \cup [1, 2]$ interval have size $\Delta x^{(3)} = \sfrac{1}{8}$.
	The mesh and initial condition are depicted in \cref{fig:WeakBlastWave_hydrodynamic_variables_init} which is provided in the appendix.
	We thus construct for every order of consistency $p = \{2, \allowbreak 3, \allowbreak 4\}$ a \ac{PERRK} scheme for which the maximum admissible timestep $\Delta t$ doubles for each family member.
	For the second-order method we use a $E_2 = \{3, \allowbreak 5, \allowbreak 9\}$ stage-evaluation method, for the third-order method a $E_3 = \{4, \allowbreak 7, \allowbreak 13\}$ method, and for the fourth-order method a $E_4 = \{6, \allowbreak 10, \allowbreak 18\}$ method.
	To solve the equation for the relaxation parameter $\gamma_n$ we use the standard Newton-Raphson method which is stopped after five steps or if the residual $r(\gamma_n)$ or steplength $\Delta \gamma_n$ are less or equal to machine epsilon $\epsilon \sim 2 \cdot 10^{-16}$.
	Here we enforce the first equality in \cref{eq:EntropyConservationDiscrete}, i.e., maintained entropy from timestep $n$ to $n+1$.
	Across orders of consistency, about two to three Newton iterations are performed on average per timestep.

\begin{figure}
		\centering
		\subfloat[{Standard \ac{PERK} methods.}]{
			\label{fig:WeakBlastWave_Standard}
			\centering
			\resizebox{.47\textwidth}{!}{\includegraphics{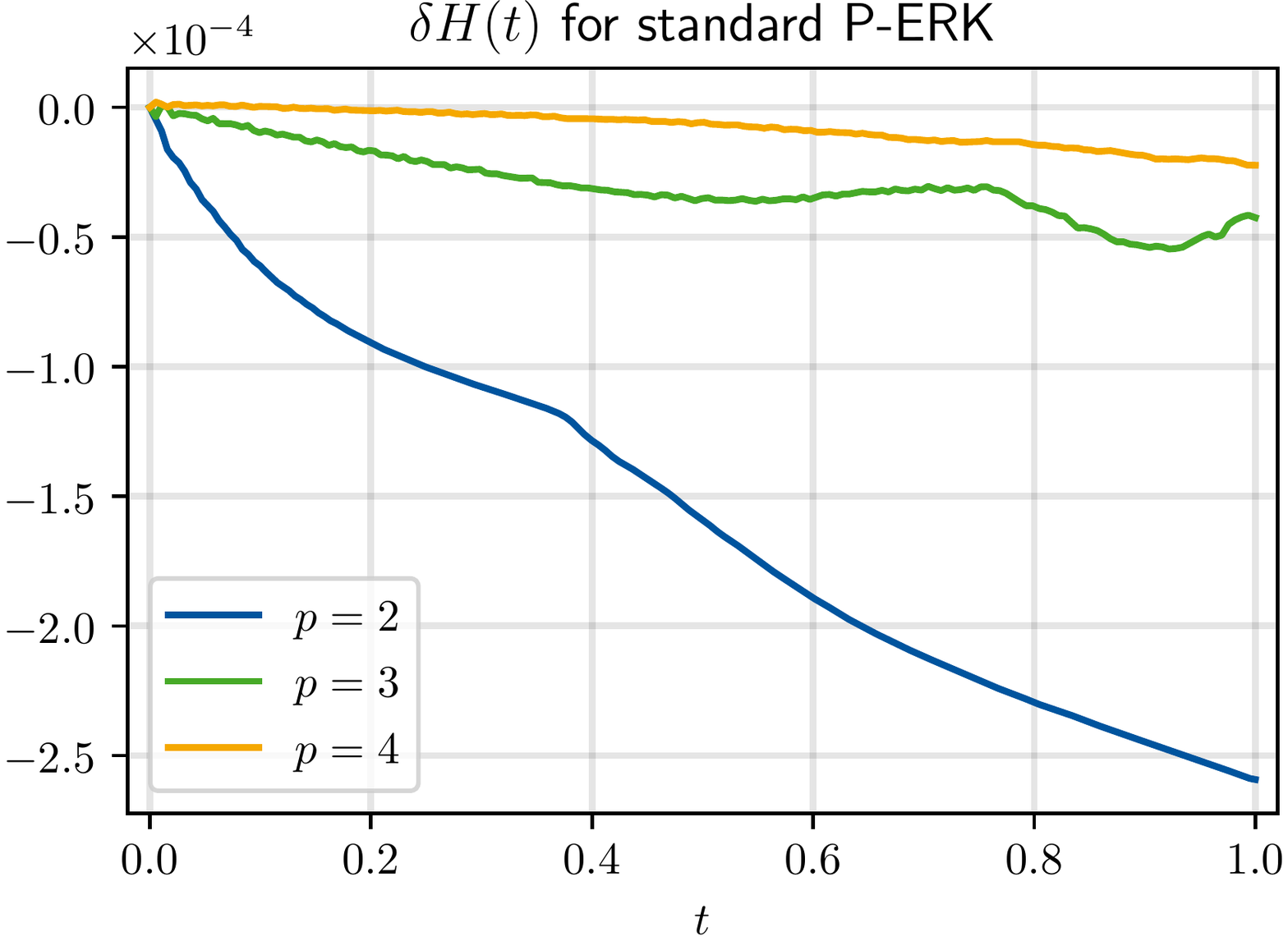}}
		}
		\hfill
		\subfloat[{Relaxation \ac{PERK}, i.e, \ac{PERRK} methods.}]{
			\label{fig:WeakBlastWave_Relaxation}
			\centering
			\resizebox{.47\textwidth}{!}{\includegraphics{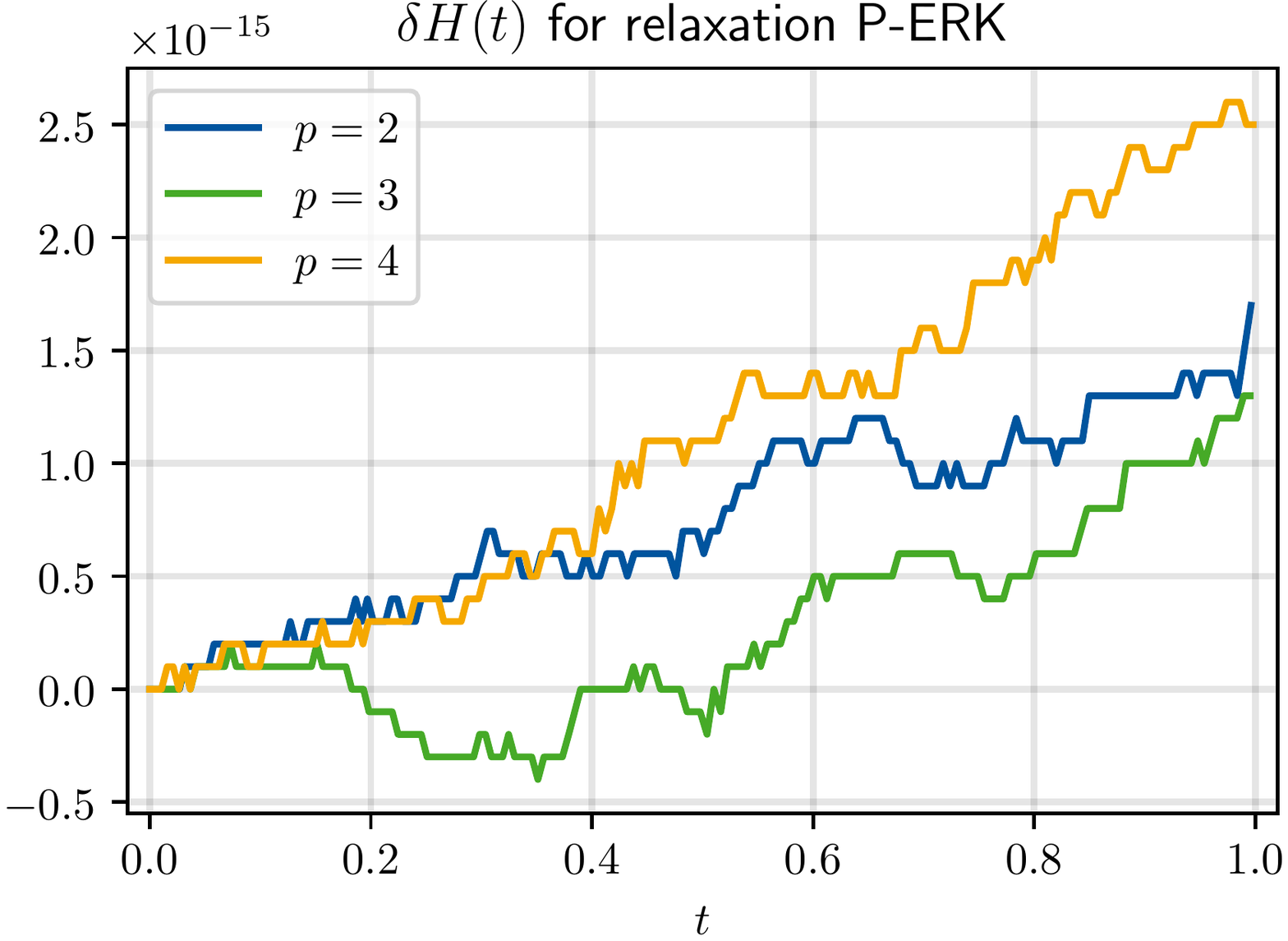}}
		}
		\caption[Temporal evolution of the total entropy defect $\delta H(t)$ for the weak blast wave problem (Euler equations) for standard \ac{PERK} methods and \ac{PERRK} methods with relaxation.]
		{Temporal evolution of the total entropy defect $\delta H(t)$ \eqref{eq:EntropyDefect} for the weak blast wave problem (Euler equations) \cref{eq:WeakBlastWave_IC} for standard \ac{PERK} methods (\cref{fig:WeakBlastWave_Standard}) and \ac{PERRK} methods with relaxation (\cref{fig:WeakBlastWave_Relaxation}).}
		\label{fig:WeakBlastWave_Standard_Relaxation}
	\end{figure}

	The evolution of the entropy change 
	\begin{equation}
		\label{eq:EntropyDefect}
		\delta H(t) \coloneqq H(t) - H \big(t_0 \big)
	\end{equation}
	is depicted in \cref{fig:WeakBlastWave_Standard_Relaxation}.
	For the standard methods we observe entropy dissipation, which is reduced for increasing order of consistency, cf. \cref{fig:WeakBlastWave_Standard}.
	For the \ac{PERRK} methods with relaxation we see that the total entropy is conserved up to orders within the rounding error regime, cf. \cref{fig:WeakBlastWave_Relaxation}.

\begin{figure}
		\centering
		\subfloat[{Standard \ac{PERK} methods.}]{
			\label{fig:WeakBlastWave_MHD_Standard}
			\centering
			\resizebox{.47\textwidth}{!}{\includegraphics{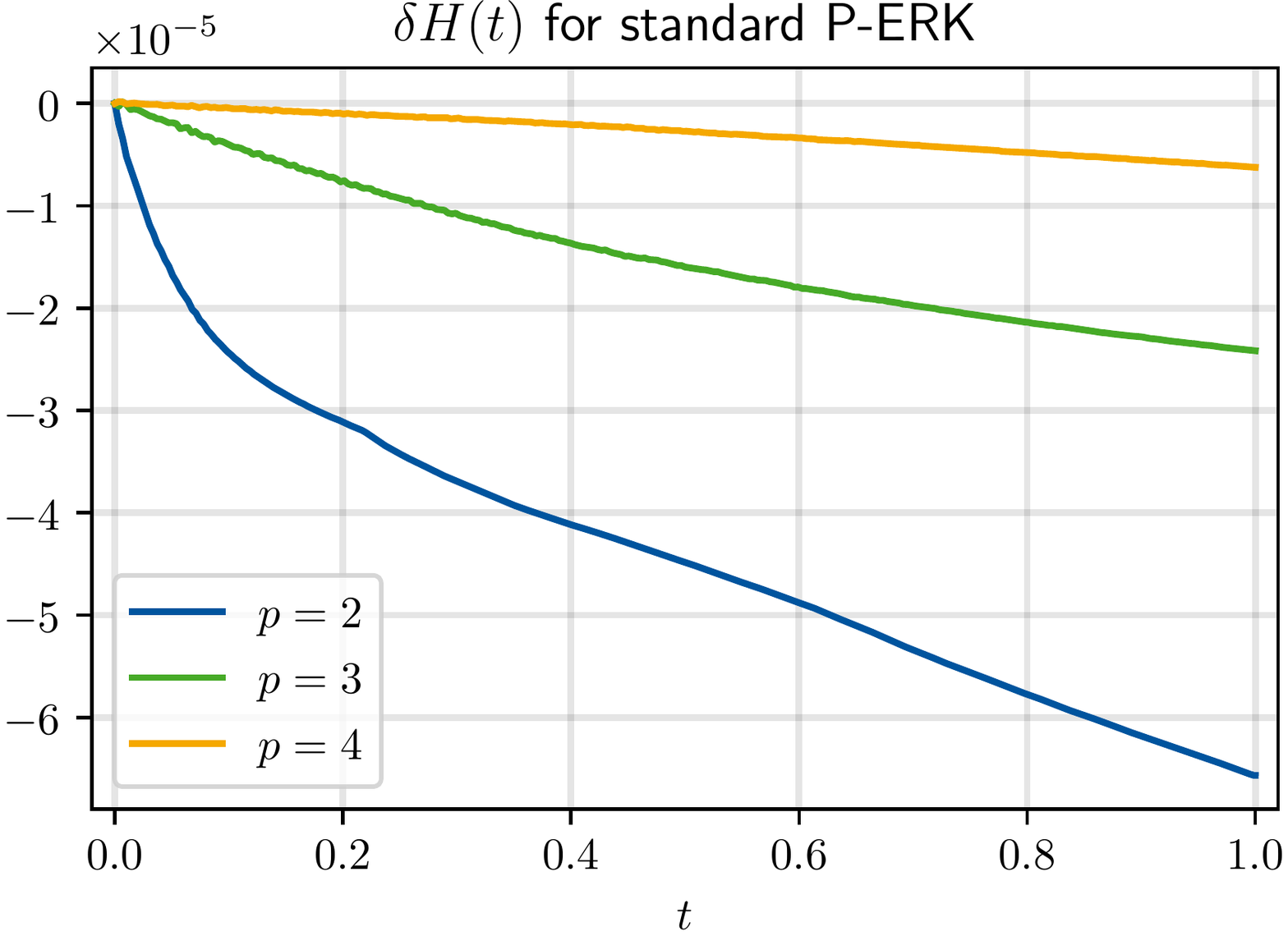}}
		}
		\hfill
		\subfloat[{Relaxation \ac{PERK}, i.e, \ac{PERRK} methods.}]{
			\label{fig:WeakBlastWave_MHD_Relaxation}
			\centering
			\resizebox{.47\textwidth}{!}{\includegraphics{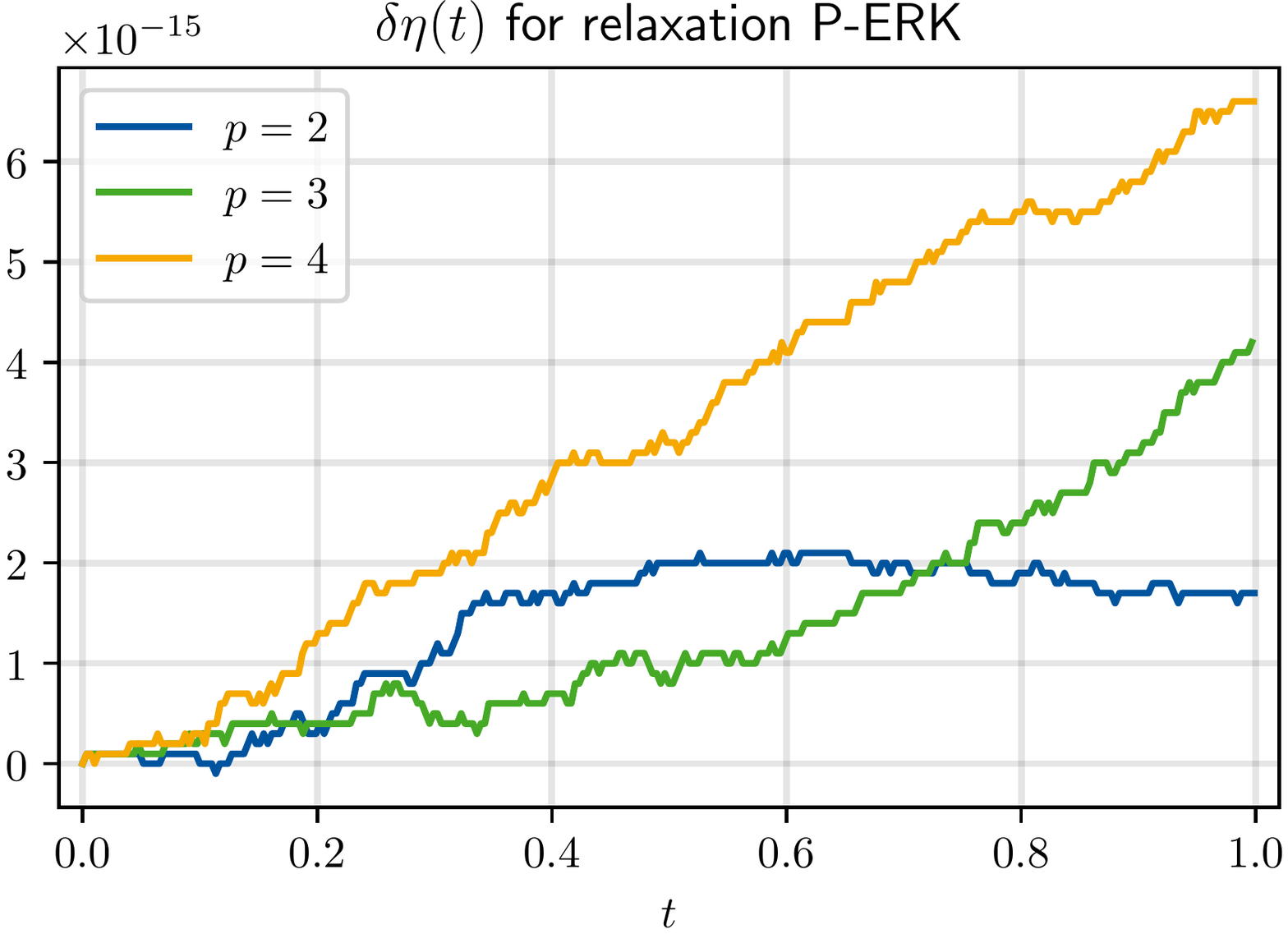}}
		}
		\caption[Temporal evolution of the total entropy defect $\delta H(t)$ for the weak blast wave problem (ideal MHD equations) for standard \ac{PERK} methods and \ac{PERRK} methods with relaxation.]
		{Temporal evolution of the total entropy defect $\delta H(t)$ \eqref{eq:EntropyDefect} for the weak blast wave problem (ideal MHD equations) for standard \ac{PERK} methods (\cref{fig:WeakBlastWave_MHD_Standard}) and \ac{PERRK} methods with relaxation (\cref{fig:WeakBlastWave_MHD_Relaxation}).}
		\label{fig:WeakBlastWave_MHD_Standard_Relaxation}
	\end{figure}

	Additionally, we consider the weak blast wave problem for the ideal \ac{MHD} equations.
	Here, we employ a "1.5"-dimensional variant of the \ac{MHD} equations by enhancing the state vector by velocity and magnetic fields in $y$ and $z$ direction, while allowing variations only in $x$ direction.
	The divergence-free constraint of the magnetic field $\boldsymbol B$
	\begin{equation}
		\label{eq:DivConstraint}
		\nabla \cdot \boldsymbol B = 0	
	\end{equation}
	is automatically satisfied for spatially univariate problems, as the $x$-component of the magnetic field stays constant for constant initializations.
	The hydrodynamic variables are initialized according to \cref{eq:WeakBlastWave_IC} while all components of the magnetic field $\boldsymbol B$ are set to unity.
	The $y$ and $z$ components of the velocity $\boldsymbol v$ are initialized to zero.
	An entropy-conservative semidiscretization is constructed analogous to the Euler case by means of an extension of the entropy-conservative flux by Ranocha \cite{Ranocha2020Entropy} to the \ac{MHD} equations \cite{hindenlang2019new}.
	The mesh, relaxation equation solver, and time integration schemes are kept identical to the previous example.
	The mathematical entropy $H$ of the \ac{MHD} equations is also the same as in the Euler case.
	On average, only one to two Newton iterations are performed per timestep for the \ac{PERRK} methods.
	We depict the evolution of the entropy defect $\delta H$ for the standard and relaxation \ac{PERRK} methods in
	\cref{fig:WeakBlastWave_MHD_Standard_Relaxation}.
	\subsubsection{Isentropic Vortex Advection}
	\label{subsec:IVA_EC}
	As a second testcase for the multirate, entropy conserving \ac{PERRK} methods we consider the classic isentropic vortex advection testcase \cite{shu1988efficient, wang2013high}.
	Here, we use the parametrization as proposed in \cite{vermeire2019paired}.
	The background flow is set to 
	\begin{equation}
		\label{eq:IsentropicVortexAdvectionBaseState}
		\rho_\infty = 1, \quad \boldsymbol v_\infty = \begin{pmatrix}
			1 \\ 1
		\end{pmatrix}, \quad p_\infty \coloneqq \frac{\rho_\infty^\gamma}{\gamma \text{Ma}_\infty^2}
	\end{equation}
	with $\text{Ma}_\infty = 0.4$ and $\gamma = 1.4$.
	To localize the effect of the vortex centered at $	\boldsymbol c(t, \boldsymbol x) \coloneqq \boldsymbol x - \boldsymbol v_\infty t$
	the perturbations are weighted with the Gaussian $g(t, x, y) \coloneqq \exp\left(\frac{1 - \Vert \boldsymbol c(t, \boldsymbol x) 
	\Vert^2}
	{2R^2}\right)$
	where $R= 1.5$.
	The size of the vortex is governed by $R$ and the intensity/strength is controlled by $I = 13.5$, following \cite{vermeire2019paired}.
	The density is then given by 
	\begin{equation}
		\label{eq:IsentropicVortexDensity}
		\rho(t, \boldsymbol x) = \rho_\infty \left( 1 - \frac{I^2 M^2 (\gamma-1) g^2(t, \boldsymbol x)}{8\pi^2} \right)^\frac{1}{\gamma -1}
	\end{equation}
	and the corresponding perturbed velocities are 
	\begin{equation}
		\boldsymbol v(t, \boldsymbol x) = \boldsymbol v_\infty + \frac{I g(t, \boldsymbol x)}{2 \pi R} \boldsymbol c(t, \boldsymbol x)
	\end{equation}
	while the pressure is computed analogous to the base pressure \cref{eq:IsentropicVortexAdvectionBaseState} as
	\begin{equation}
		\label{eq:IsentropicVortexPressure}
		p(t, \boldsymbol x) =  \frac{\rho^\gamma(t, \boldsymbol x)}{\gamma \text{Ma}_\infty^2} \: .
	\end{equation}
	The vortex is initially centered in the $\Omega = [-L, L]^2, L = 10$ domain equipped with periodic boundary conditions.

\begin{figure}
		\centering
		\includegraphics[width=0.55\textwidth]{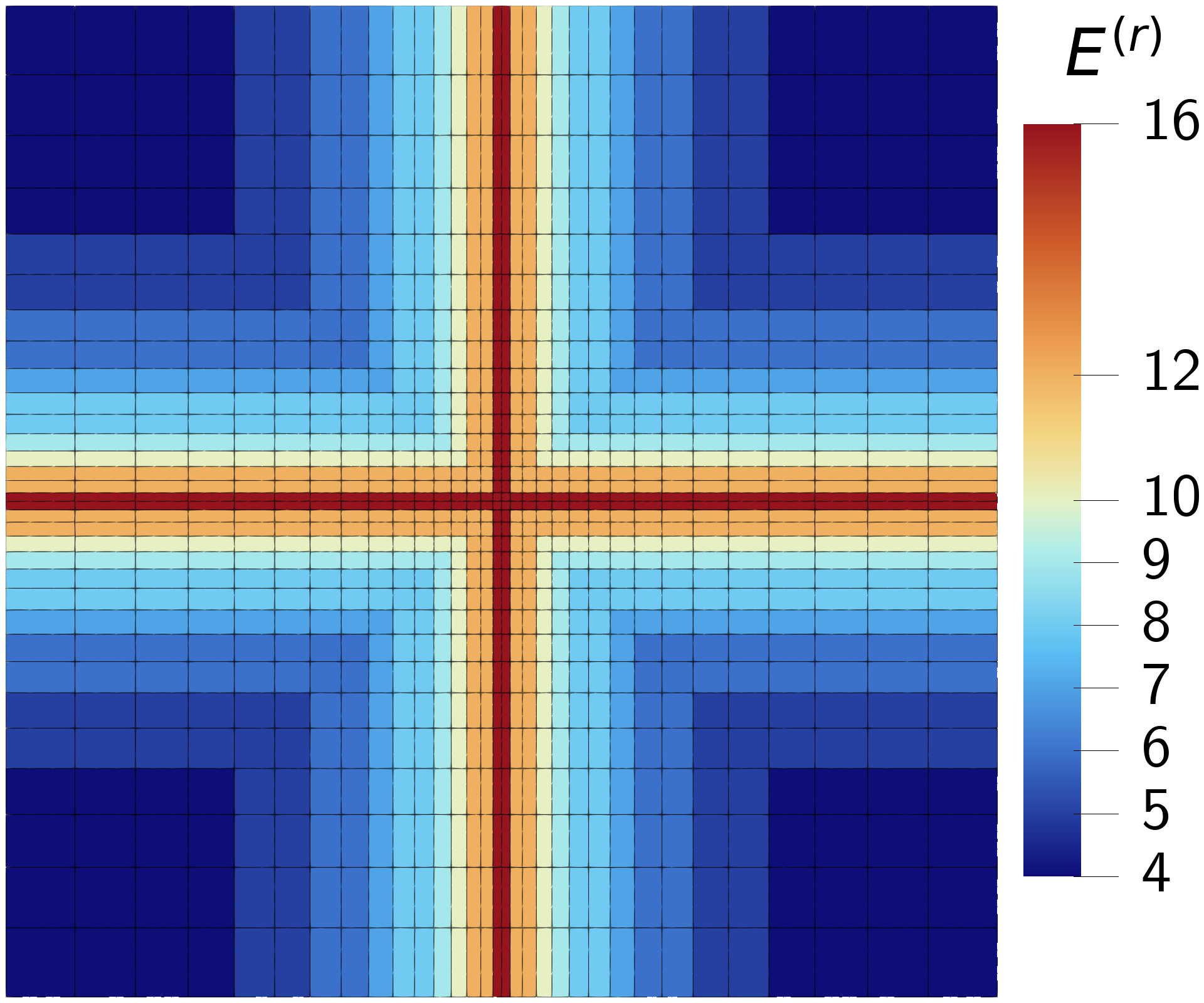}
		\caption[Non-uniform, conforming mesh for the isentropic vortex testcase to validate entropy conservation.]
		{Non-uniform, conforming mesh according to \cref{eq:IsentropicVortexMesh} for the isentropic vortex testcase to validate entropy conservation.
		The cells are colored according to the number of active stages $E_p^{(r)}$ for the third-order method $E_3$. The number of active stages $E_3^{(r)}$ is determined based on the shortest edge length $h$ of each cell.}
		\label{fig:IsentropicVortexMesh}
	\end{figure}

	As for the previous testcase, we employ $k=3$ local solution polynomials.
	To discretize $\Omega$, a non-uniform, conforming mesh consisting of $32 \times 32$ cells is used.
	This is obtained by mapping the uniformly discretized reference coordinates $\boldsymbol \xi \in [-1, 1]^2 $ to the physical domain via the transformation
	\begin{equation}
		\label{eq:IsentropicVortexMesh}
		\begin{pmatrix}
			x \\ y
		\end{pmatrix}
		= L \begin{pmatrix}
			\text{sgn}(\xi_1) \, \vert \xi_1 \vert^{1.4 + \vert \xi_1 \vert} \\
			\text{sgn}(\xi_2) \, \vert \xi_2 \vert^{1.4 + \vert \xi_2 \vert}
		\end{pmatrix}
	\end{equation}
	which results in "$+$"-like refined regions, cf. \cref{fig:IsentropicVortexMesh}.
	The ratio of the smallest to largest element edge length is in this case roughly a factor of eight.
	To integrate the \ac{PDE} on this uniform mesh efficiently, the following \ac{PERK}/\ac{PERRK} methods are constructed:
	$E_2 = \{4, \allowbreak 6, \allowbreak 8, \allowbreak 10, \allowbreak 12, \allowbreak 16 \}; \: \allowbreak E_3 = \{4, \allowbreak 5, \allowbreak 6, \allowbreak 7, \allowbreak 8, \allowbreak 9, \allowbreak 10, \allowbreak 12, \allowbreak 16 \}; \: \allowbreak E_4 = \{5, \allowbreak 6, \allowbreak 7, \allowbreak 9, \allowbreak 11, \allowbreak 16 \}$.
	The mesh according to \cref{eq:IsentropicVortexMesh} is shown in \cref{fig:IsentropicVortexMesh} with element assignment corresponding to $E_3$.
	%

\begin{figure}
		\centering
		\subfloat[{Standard \ac{PERK} methods.}]{
			\label{fig:IsentropicVortex_Standard}
			\centering
			\resizebox{.47\textwidth}{!}{\includegraphics{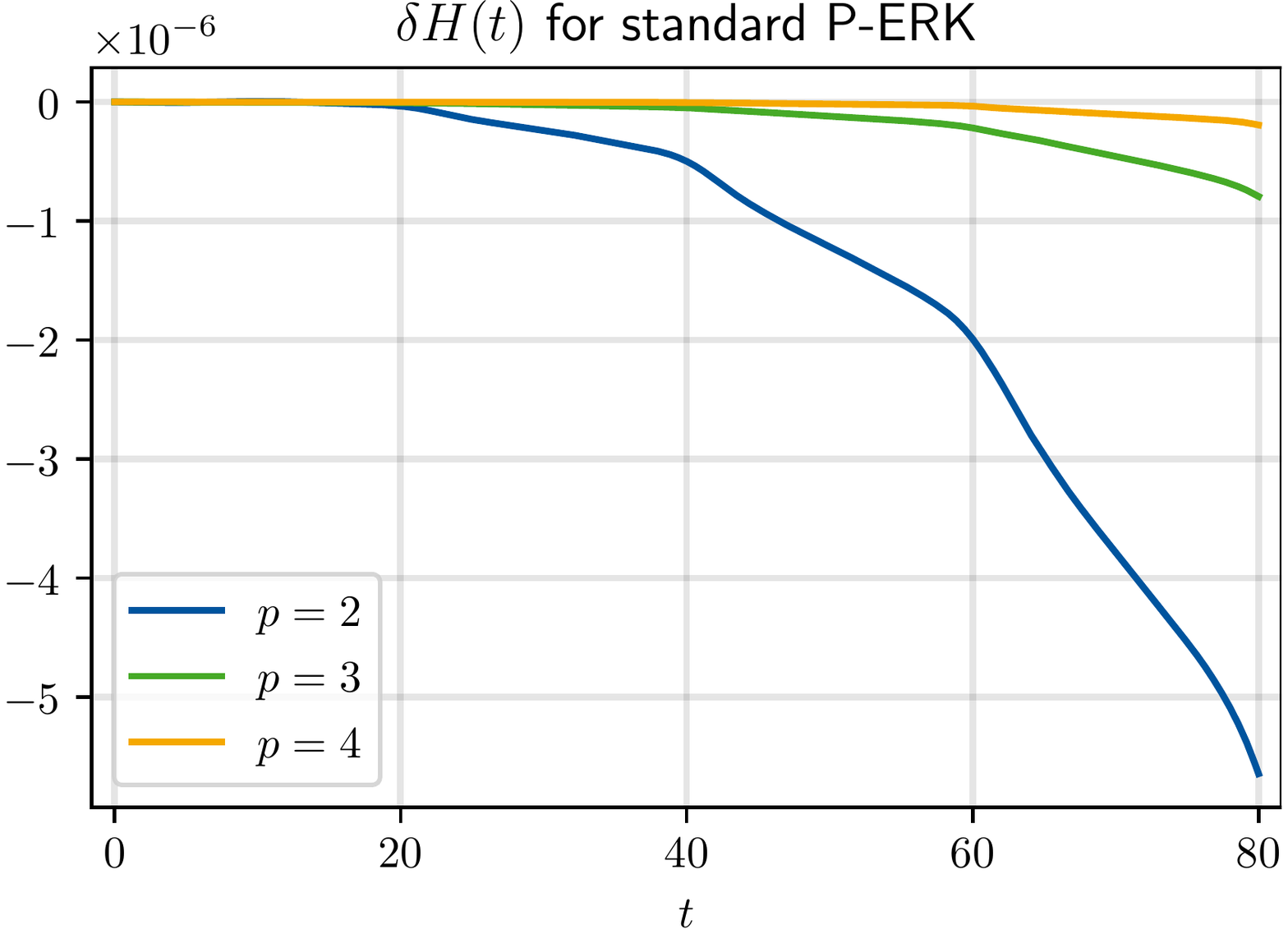}}
		}
		\hfill
		\subfloat[{Relaxation \ac{PERK}, i.e, \ac{PERRK} methods.}]{
			\label{fig:IsentropicVortex_Relaxation}
			\centering
			\resizebox{.47\textwidth}{!}{\includegraphics{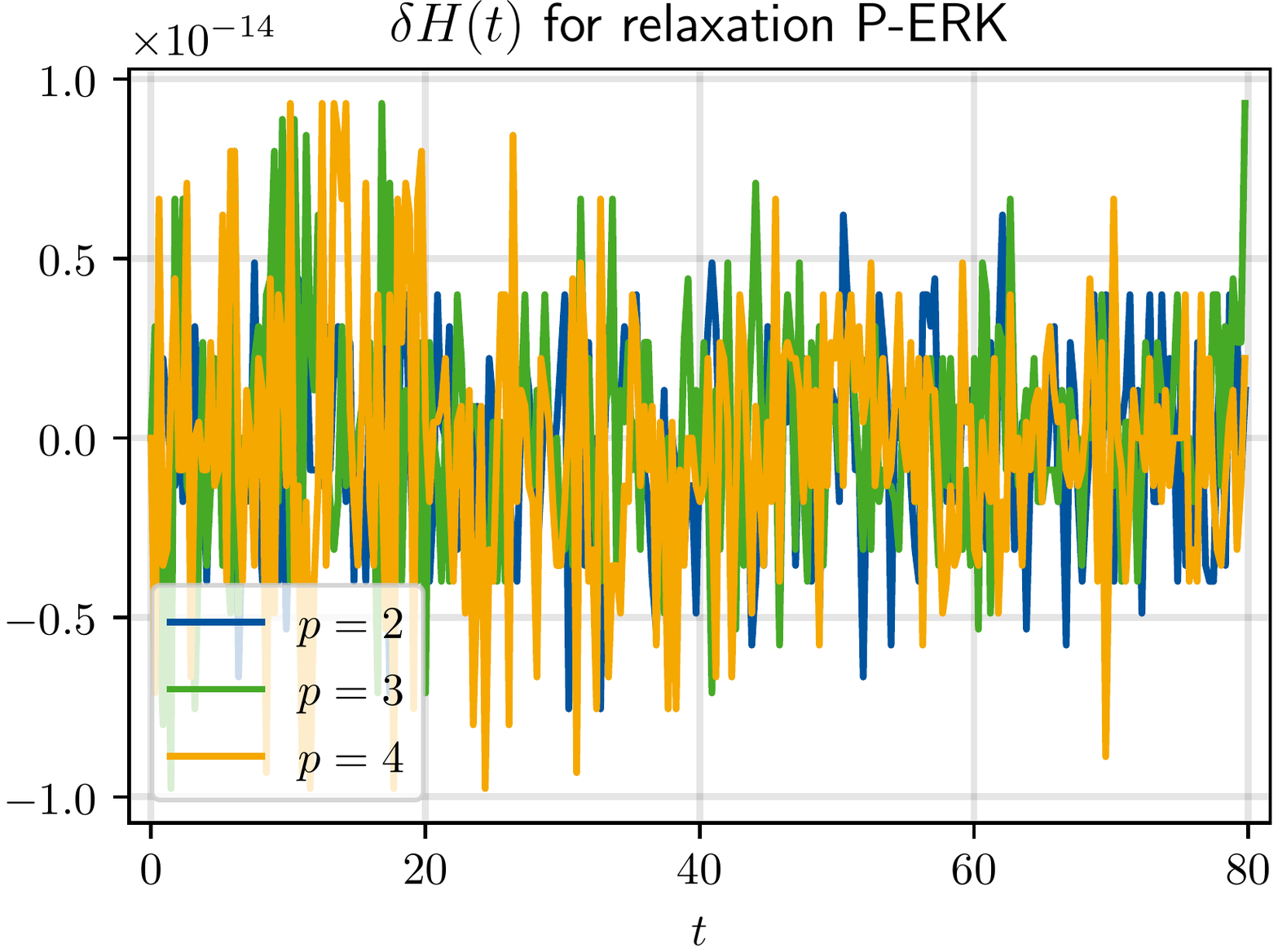}}
		}
		\caption[Temporal evolution of the total entropy defect $\delta H(t)$ for the isentropic vortex advection testcase for standard \ac{PERK} methods and \ac{PERRK} methods with relaxation.]
		{Temporal evolution of the total entropy defect $\delta H(t)$ \eqref{eq:EntropyDefect} for the isentropic vortex testcase for standard \ac{PERK} methods (\cref{fig:IsentropicVortex_Standard}) and \ac{PERRK} methods with relaxation (\cref{fig:IsentropicVortex_Relaxation}).}
		\label{fig:IsentropicVortex_Standard_Relaxation}
	\end{figure}

	We compute four passes through the domain ($t_f = 80$) with a uniform timestep for all methods.
	As for the previous testcase we observe that the standard \ac{PERK} methods exhibit entropy dissipation, cf. \cref{fig:IsentropicVortex_Standard}, which is more significant for lower order of consistency.
	In turn, the \ac{PERRK} methods with relaxation conserve the total entropy much better, with entropy fluctuations on the order of $10^{-15}$, cf. \cref{fig:IsentropicVortex_Relaxation}.
	The relaxation parameter is determined using the standard Newton-Raphson method which stops after three iterations, if the residual is less or equal $10^{-14}$ or the update step is less or equal $10^{-15}$.
	For this case, the second equality in \cref{eq:EntropyConservationDiscrete} is enforced, i.e., the entropy at any timestep $n$ is forced to be to equal to the entropy of the initial condition $H(\boldsymbol U_0)$.
	On average, these tolerances are met after one to two iterations.
	\subsection{Conservation of Linear Invariants}
	\label{subsec:ConservationOfLinearInvariants}
	We briefly demonstrate that the \ac{PERRK} schemes do not violate linear invariants.
	To this end, we resort again to the isentropic vortex testcase as described in the previous section.
	Here, the compressible Euler equations are discretized on a $32 \times 32$ uniform mesh with $k=3$ solution polynomials.
	The Riemann problems at cell interfaces are solved approximately with the \ac{HLLC} flux \cite{toro1994restoration} and we use a standard weak form volume integral.
	We compute ten passes through the domain, i.e., $t_f = 200$ with uniform timestep $\Delta t = 0.01$ resulting in $20000$ time steps.
	We limit ourselves to fourth-order, $E = \{5, \allowbreak 6, \allowbreak \dots, \allowbreak 14, \allowbreak 15 \}$ \ac{PERK} schemes, both with and without relaxation.
	The $1024$ grid cells are randomly assigned to one of the schemes to maximize the number of interfaces whose neighboring elements are integrated with different schemes.
	Note that interfaces between cells on different levels are always assigned to the higher-stage method.
	The conservation error
	\begin{equation}
		e_u^\text{Cons}(t_f) \coloneqq \left \vert \int_\Omega u(t_f, \boldsymbol x) \nid \boldsymbol x - \int_\Omega u \big (t_0, \boldsymbol x \big ) \nid \boldsymbol x \right \vert
	\end{equation}
	is computed for mass, momenta and energy at final time and tabulated in \cref{tab:ConservationErrors_IsentropicVortex}.
	As the median number of stage evaluations is $10$ we compare the conservation errors to the fourth-order, 10-stage, low-storage, \ac{SSP} method by Ketcheson \cite{Ketcheson2008highly} which is implemented in \texttt{OrdinaryDiffEq.jl} \cite{DifferentialEquations.jl-2017}.
	This results in approximately $2 \cdot 10^5$ \ac{ODE} \ac{RHS} calls across all methods, assuming even distribution of the grid cells among the \ac{PERRK} schemes.
	We see from \cref{tab:ConservationErrors_IsentropicVortex} that the conservation violation of both the \ac{PERK} and \ac{PERRK} schemes is on the order of round-off errors for double-precision floating point numbers.
	In particular, the conservation errors for the $10$-stage \ac{SSP} scheme are larger for every conserved quantity.
	\begin{table}
		\def\arraystretch{1.2}
		\centering
		\begin{tabular}{l?{2}c|c|c}
				& $\text{P-ERK}_{\{5, \dots, 15\}}$ & $\text{P-ERRK}_{\{5, \dots, 15\}}$ & $\text{SSP}_{10}$ \cite{Ketcheson2008highly} \\
					\Xhline{5\arrayrulewidth}
				$e_\rho^\text{Cons}(t_f)$       & $1.58 \cdot 10^{-13}$ & $1.57 \cdot 10^{-13} $ & $1.10 \cdot 10^{-12} $ \\
				$e_{\rho v_x}^\text{Cons}(t_f)$ & $2.56 \cdot 10^{-14}$ & $1.72 \cdot 10^{-14} $ & $9.70 \cdot 10^{-13} $ \\
				$e_{\rho v_y}^\text{Cons}(t_f)$ & $1.64 \cdot 10^{-14}$ & $2.23 \cdot 10^{-14} $ & $9.51 \cdot 10^{-13} $ \\
				$e_{\rho e}^\text{Cons}(t_f)$   & $4.53 \cdot 10^{-13}$ & $3.96 \cdot 10^{-13} $ & $1.36 \cdot 10^{-11} $ \\
		\end{tabular}
		\caption[Conservation-errors for isentropic vortex advection testcase.]
		{Conservation-errors for isentropic vortex advection testcase after ten passes through the domain at $t_f = 200$ for the fourth-order \ac{PERK} scheme with and without relaxation and the ten-stage, fourth-order \ac{SSP} method \cite{Ketcheson2008highly}.}
		\label{tab:ConservationErrors_IsentropicVortex}
	\end{table}
	\subsection{Order of Convergence}
	\label{subsec:OrderOfConvergence}
	To demonstrate that the relaxed \ac{PERK} schemes still meet their designed order of accuracy we consider two classical test cases for compressible flows.
	These were also used in \cite{RelaxationRanocha} to validate the relaxation methodology.
	In contrast to testcases based on manufactured solutions, these tests do not include source terms to match the constructed solution.
	First, to test inviscid, i.e., hyperbolic problems we employ the isentropic vortex advection problem \cite{shu1988efficient, wang2013high} as described in the previous section.
	Second, in order to examine viscous, i.e., hyperbolic-parabolic problems we consider the viscous shock propagation problem \cite{becker1922stosswelle, morduchow1949complete, margolin2017entropy}.
	For these testcases, the convergence studies are performed on non-uniform meshes, i.e., we employ multirate setups to test the convergence of the relaxed multirate \ac{PERK} methods.

	In addition, we test the relaxation methodology for systems involving non-conservative fluxes by means of the \ac{glm-mhd} equations \cite{munz2000divergence, dedner2002hyperbolic, derigs2018ideal}.
	This testcase is performed on a uniform mesh where the cells are randomly assigned to one out of five \ac{PERRK} family member schemes.
	\subsubsection{Isentropic Vortex Advection}
	To demonstrate that the \ac{PERRK} methods with relaxation maintain their designed order of accuracy we consider the isentropic vortex advection testcase with minimal spatial errors.
	To this end we employ $k=6$ \ac{DG} solution polynomials on a $64 \times 64$ base mesh, which is dynamically refined twice towards the center of the vortex.
	Precisely, the cells with distance $r \in [3, 2)$ to the center of the vortex are refined once while the cells with $r \in [2, 0]$ are refined twice.
	The vortex is parametrized exactly as in \cref{subsec:IVA_EC} where we demonstrated entropy conservation.
	Since currently only $L^2$ projection/interpolation mortars \cite{kopriva1996conservative} are implemented in Trixi.jl \cite{trixi1, trixi2, trixi3} this setup with \ac{AMR} is not exactly entropy-conservative.
	We remark that a possible remedy to be explored in future work is the flux-correction method proposed in \cite{chan2021mortar}.

	This setup allows for an initial domain-normalized $L^1$-error 
	\begin{equation}
		\label{eq:L1ErrorDomainNormalized_Density}
		e_\rho^{L^1}(t) \coloneqq \frac{1}{\vert \Omega \vert} \int_\Omega \left \vert \rho(t, \boldsymbol x) - \rho_{\text{exact}}(t, \boldsymbol x) \right \vert \nid \boldsymbol x
	\end{equation}
	of order $10^{-13}$.
	This error is computed from the quadrature defect of the solution representation using the $k=6$ polynomials to $\widetilde{k} = 2 \cdot k = 12$ degree polynomials.
	The mesh is adapted $100$ times over the course of the simulation, i.e., with time interval $\Delta t_{\text{AMR}} = 0.2$ which is kept constant over the convergence study.
	We again use the \ac{HLLC} \cite{toro1994restoration} flux for the Riemann problems at cell interfaces and a standard weak form volume integral.
	$L^1, L^\infty$-errors of the density $\rho$ are computed for $E_2 = \{3, \allowbreak 6, \allowbreak 12 \}; \, \allowbreak E_3 = \{4, \allowbreak 8, \allowbreak 16 \}; \, \allowbreak E_4 = \{5, \allowbreak 9, \allowbreak 15 \}$ \ac{PERRK} schemes after one pass of the vortex.
	The base timestep $\Delta t = 4 \cdot 10^{-3}$ is the same for all methods which is halved until the errors in density start to saturate due to the spatial discretization accuracy.
	To solve the relaxation equation we employ the Newton-Raphson method with a maximum of five iterations, a residual tolerance of $10^{-13}$ for the $p=2$ method and $10^{-14}$ for the $p=3, 4$ methods.
	For all orders of accuracy we supply a stepsize tolerance of double-precision machine epsilon.
	
\begin{figure}
		\centering
		\subfloat[{Domain-normalized $L^1$-error in density $\rho$, cf. \cref{eq:L1ErrorDomainNormalized_Density}.}]{
			\label{fig:Convergence_IVA_L1}
			\centering
			\resizebox{.47\textwidth}{!}{\includegraphics{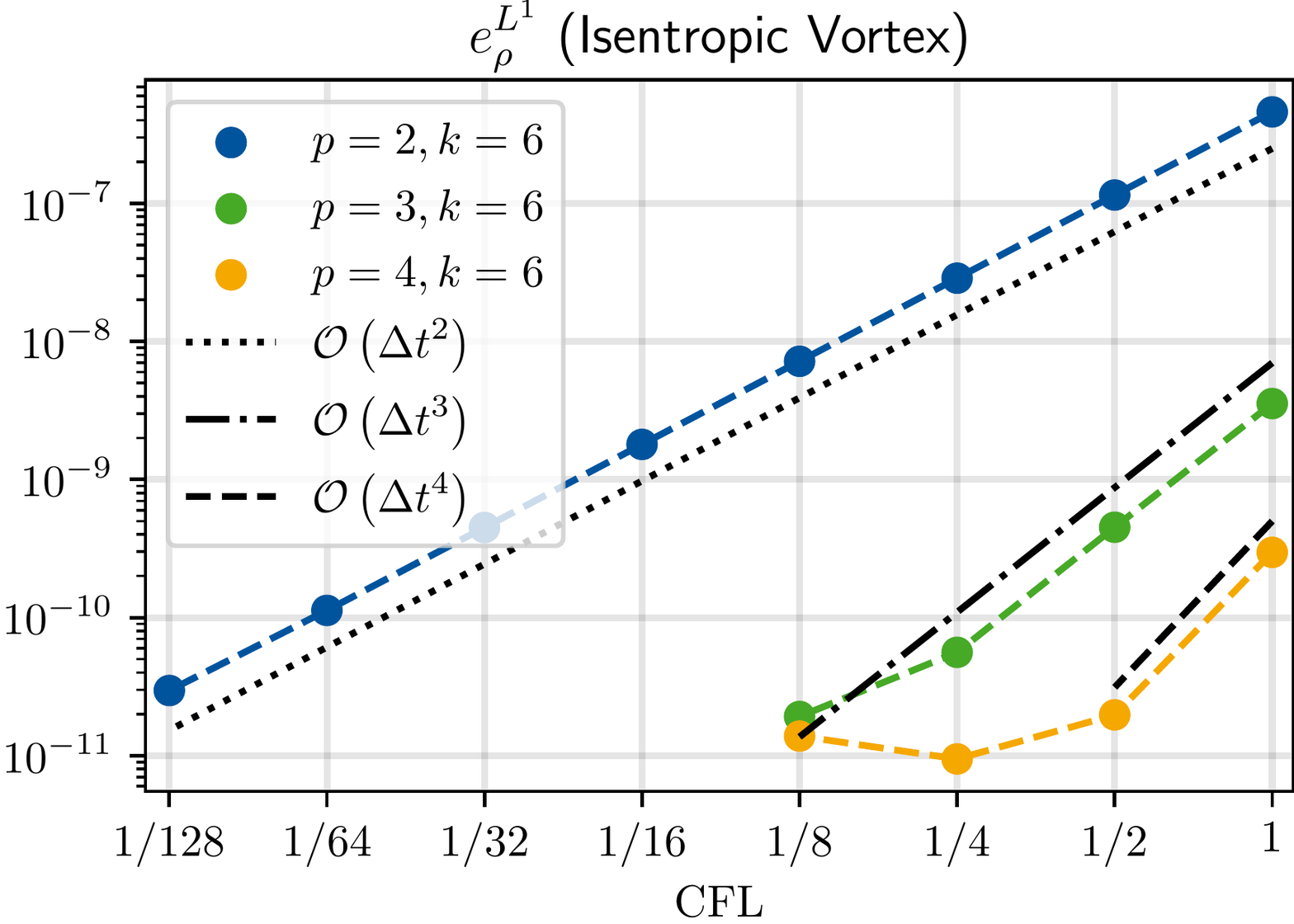}}
		}
		\hfill
		\subfloat[{$L^\infty$ error in density $\rho$.}]{
			\label{fig:Convergence_IVA_LInf}
			\centering
			\resizebox{.47\textwidth}{!}{\includegraphics{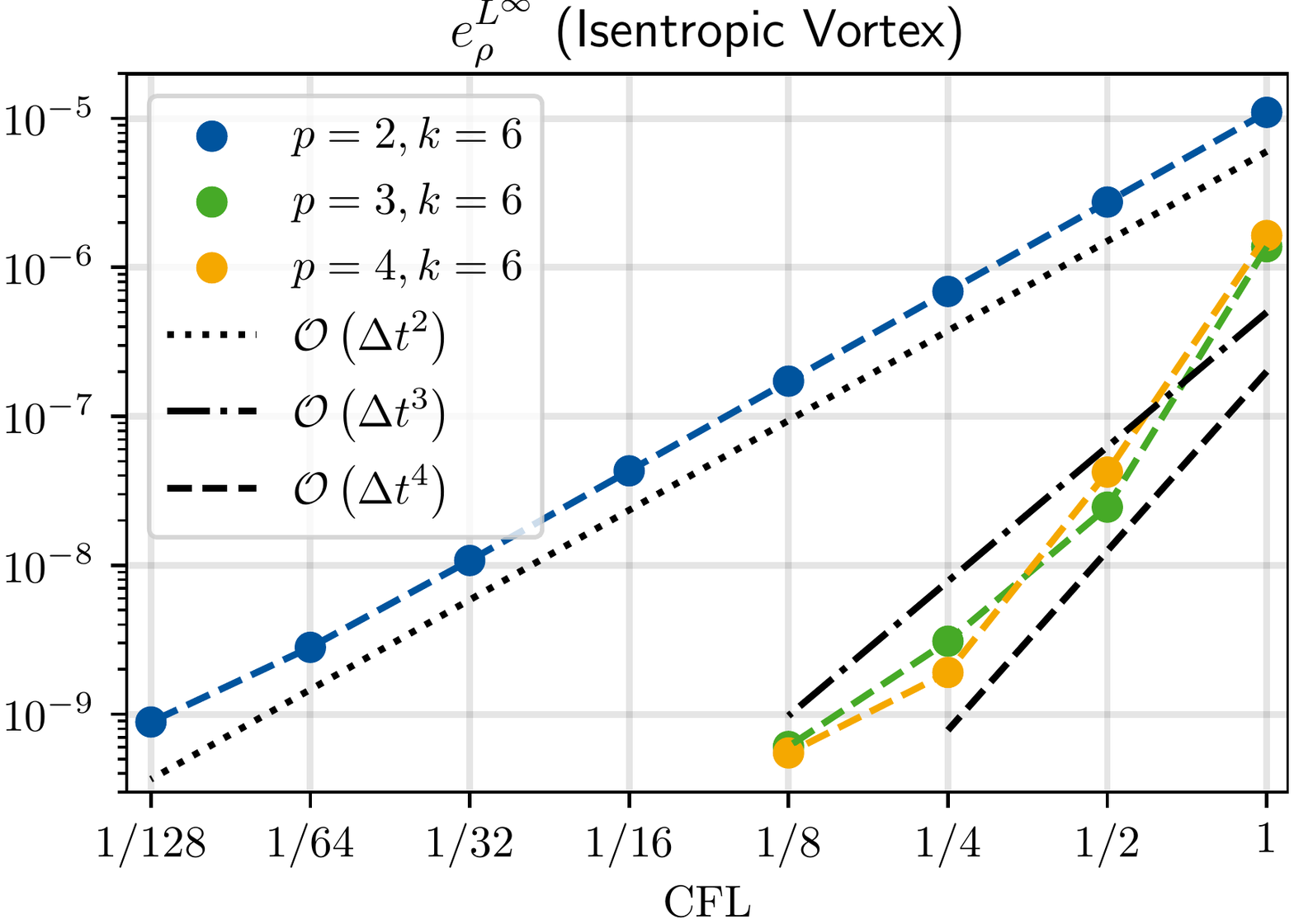}}
		}
		\caption[$L^1$ and $L^\infty$ errors of the density $\rho$ for the isentropic vortex testcase.]
		{$L^1$ (\cref{fig:Convergence_IVA_L1}) and $L^\infty$ (\cref{fig:Convergence_IVA_LInf}) errors of the density $\rho$ for the isentropic vortex advection testcase after one pass through the domain for $p=2, 3, 4$ \ac{PERRK} schemes.}
		\label{fig:Convergence_IVA_L1_LInf}
	\end{figure}

	For all schemes we observe the designed order of convergence in both $L^1$ and $L^\infty$ density errors, see \cref{fig:Convergence_IVA_L1_LInf}.
	While the $p=2$ requires a significant reduction of the reference timestep to reach the spatial discretization limit, the $p=3$ and especially the $p=4$ schemes reach this limit earlier.
	In particular, we observe for the $p=4$ scheme that only the step from $\text{CFL} = 1$ to $\text{CFL} = 0.5$ reduces the $L^1$ error by the expected factor of $16$, before flattening out, cf. \cref{fig:Convergence_IVA_L1}.
	For the $L^\infty$ error we observe superconvergence for the $p=3, \allowbreak p=4$ schemes for the first reduction from $\text{CFL} = 1$ to $\text{CFL} = 0.5$, then one step with designed order of convergence, followed by a saturation of the error due to the spatial accuracy.

\begin{figure}
		\centering
		\subfloat[{Relaxation parameter $\gamma_n$ over time for $p=2$.}]{
			\label{fig:gamma_IVA_p2}
			\centering
			\resizebox{.47\textwidth}{!}{\includegraphics{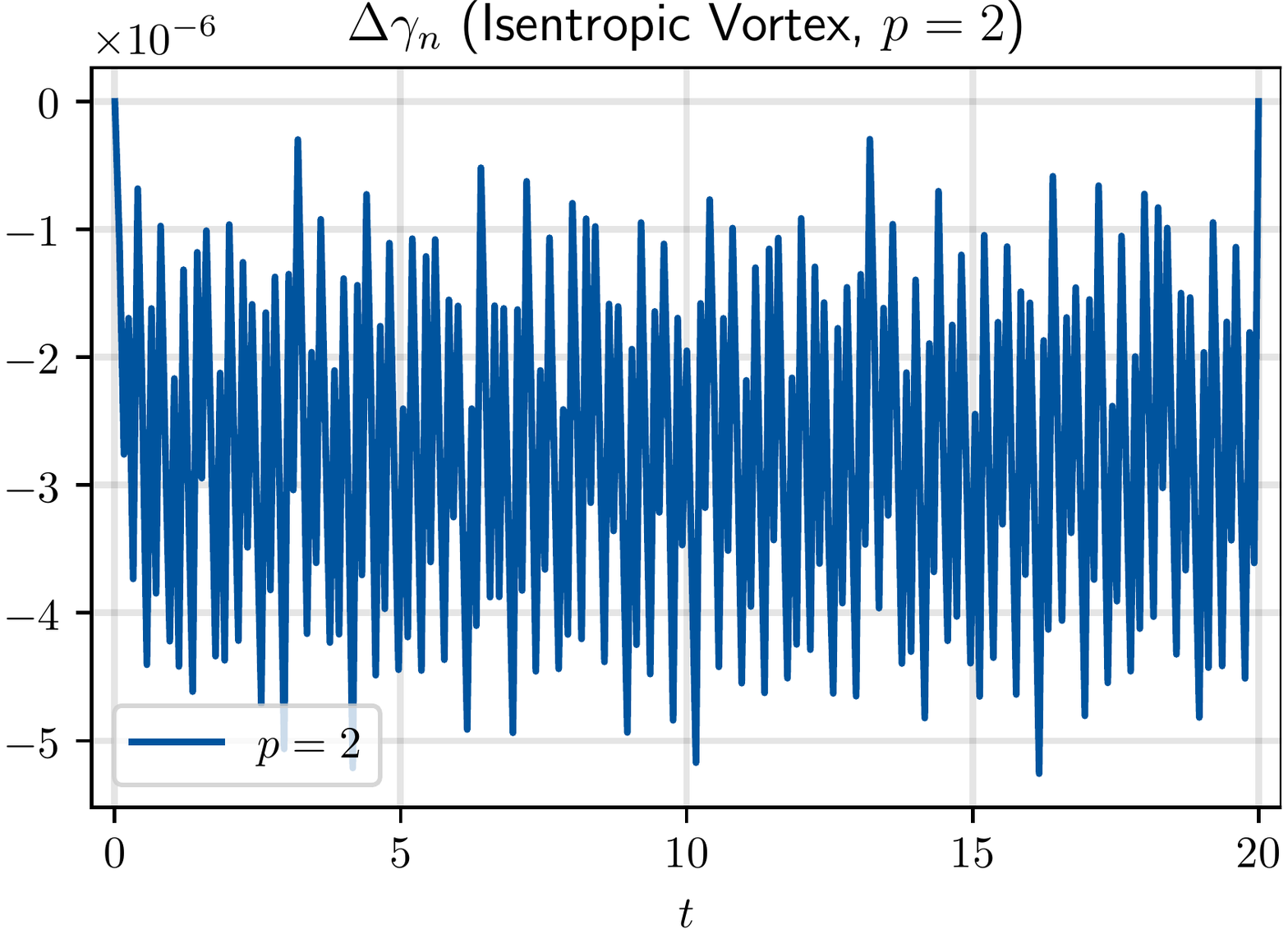}}
		}
		\hfill
		\subfloat[{Relaxation parameter $\gamma_n$ over time for $p=3, 4$.}]{
			\label{fig:gamma_IVA_p3p4}
			\centering
			\resizebox{.47\textwidth}{!}{\includegraphics{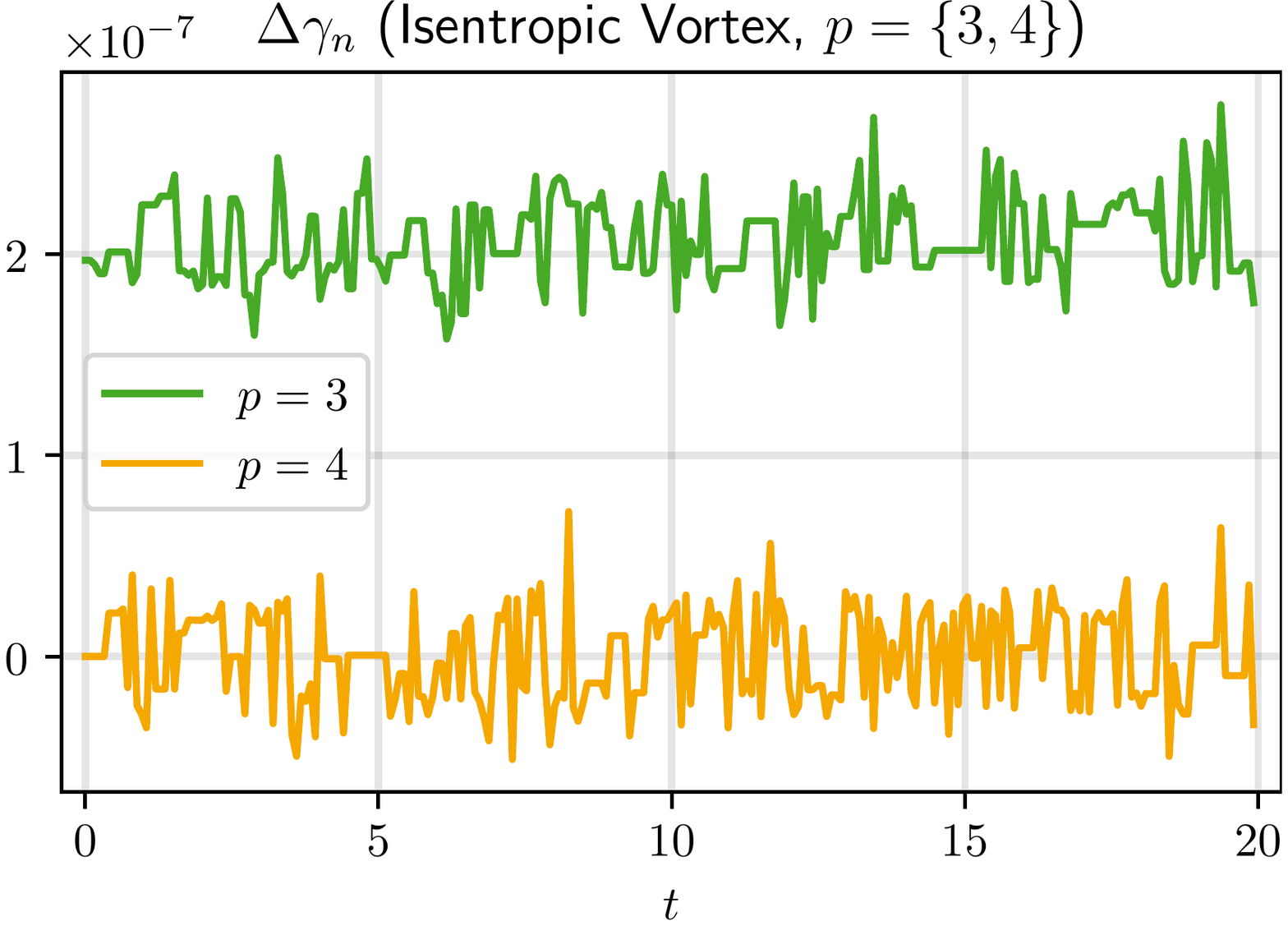}}
		}
		\caption[Evolution of the relaxation parameter $\gamma_n$ for the isentropic vortex convergence testcase.]
		{Evolution of the relaxation parameter $\gamma_n$ for the isentropic vortex convergence testcase with $\text{CFL} = 1$ for $p=2$ (\cref{fig:gamma_IVA_p2}) and $p=3, 4$ (\cref{fig:gamma_IVA_p3p4}) \ac{PERRK} schemes.
		In the plots above, the deviation from unity, i.e., $\Delta \gamma_n \coloneqq \gamma_n - 1$ is shown.}
		\label{fig:gamma_IVA_p2_p3p4}
	\end{figure}

	In \cref{subsec:RelaxationRungeKuttaMethods} we mentioned that the desired order of convergence $p$ should be reached for relaxation parameters that do not deviate too far from unity, i.e., $\gamma_n = 1 + \mathcal O \left(\Delta t^{p-1} \right)$ \cite{RelaxationKetcheson, RelaxationRanocha}.
	Thus, we exemplarily examine the evolution of the relaxation parameter $\gamma_n$ over time for the $\text{CFL} = 1$ case for all schemes.
	In \cref{fig:gamma_IVA_p2_p3p4} we present the relaxation parameter sampled at every 20th timestep.
	We observe that the deviations $\Delta \gamma_n$ from one scale as $\mathcal O \left( 10^{-(4 + p)} \right)$.
	In particular, we have for the $p=4$ scheme that $\Delta t^3 = 6.4 \cdot 10^{-8}$, thus we are in safe waters according to theory.
	For the other two cases, i.e., $p=\{2, 3\}$ the deviations are even smaller in relation to $\Delta t^{p-1}$ and thus convergence with the designed order of accuracy is expected and observed.
	\subsubsection{Viscous Shock Propagation}
	\label{subsec:ConvTest_ViscousShock}
	To demonstrate that the \ac{PERRK} methods with relaxation maintain their designed order of accuracy for entropy diffusive problems we consider the viscous shock propagation testcase.
	This testcase is constructed by balancing shock-formation due to the hyperbolic nonlinearity through diffusion, resulting in a time-invariant shock profile.
	Becker \cite{becker1922stosswelle} was the first one to derive the a solution for the shock profile of the 1D Navier-Stokes-Fourier equations, i.e., taking both viscosity and heat conduction into account.
	This requires setting the Prandtl number $\text{Pr} = \frac{3}{4} = c_p \frac{\mu}{\kappa} = \frac{\gamma}{\gamma - 1} \frac{\mu}{\kappa}$ which is not too far from physical values for air and other gases.
	Becker's work was revisited by Morduchow and Libby \cite{morduchow1949complete} which showed that some of Becker's assumptions could be relaxed.

	For a general parametrization, the solution for the 1D shock profile involves a nonlinear equation which needs to be solved numerically to obtain a solution to the Becker equation of motion \cite{morduchow1949complete, margolin2017entropy}.
	In \cite{margolin2017entropy}, however, the authors show that for the Mach number $\text{Ma} = \frac{2}{\sqrt{3 - \gamma}}, \allowbreak 1 \leq \gamma \leq \frac{5}{3}$ an explicit solution is available.
	Here, we set $\gamma = \frac{5}{3}$ (monatomic ideal gas), dynamic viscosity $\mu = 0.15$, reference density $\rho_0 = 1$, and shock speed $v = 1$.
	From this, we compute the Mach number, speed of sound ahead of the shock $c_0 = \frac{v}{\text{Ma}}$, and the reference pressure $p_0 = c_0^2 \frac{\rho_0}{\gamma}$ which is computed from the fact that the enthalpy is constant across the shock.
	Furthermore, we have the length scale $l = \frac{\mu}{\gamma - 1} \frac{2 \gamma}{\gamma - 1} \frac{1}{\rho_0 v}$ \cite{margolin2017entropy}.
	Since the shock profile is moving with velocity $v$, the solution at any time $t$ can be obtained by tracing back to the solution at initial position $x_0(x, t) \coloneqq x - v t$.
	With characteristic position 
	\begin{equation}
		\chi(x) \coloneqq 2 \exp\left( \frac{x}{2l} \right)	
	\end{equation}
	and "effective momentum"
	\begin{equation}
		w(\chi) \coloneqq 1 + \frac{1}{2 \chi^2} \left( 1 - \sqrt{1 + 2 \chi^2}\right)
	\end{equation}
	the primal variables are given by \cite{margolin2017entropy}
	\begin{subequations}
		\label{eq:ViscousShockSolution}
		\begin{align}
			\rho(t, x) &= \frac{\rho_0}{w} \\
			u(t, x) &= v (1 - w) \\
			p(t, x) &= \frac{p_0}{w} \left( 1 + \frac{\gamma - 1}{2} \text{Ma}^2 \left (1 - w^2 \right ) \right)
		\end{align}
	\end{subequations}
	where we abbreviate $w \equiv w(\chi(x_0(x, t)))$.
	The solution \cref{eq:ViscousShockSolution} at initial time is provided by \cref{fig:ViscousShock_Init} in the appendix.

	We solve this problem on $\Omega = [-2, 2]$ up to $t_f = 0.5$.
	For the hyperbolic boundary conditions, i.e., the values for the inviscid fluxes, we use inflow conditions obtained from the analytical solution at the left boundary and (free) outflow conditions at the right boundary.
	The parabolic boundary conditions required for the viscous gradient fluxes are realized as moving isothermal, no-slip walls where velocity and temperature are obtained from the analytical solution.
	The surface fluxes are computed using the \ac{HLLE} \cite{einfeldt1988godunov} flux and the viscous terms are discretized using the \ac{BR1} scheme \cite{bassi1997high, gassner2018br1}.
	Since this problem is diffusion-dominated we also tested the local Discontinuous Galerkin (LDG) method \cite{cockburn1998local} with zero penalty parameter \cite{cockburn2007analysis} which yields similar results.
	For the volume integral we use the standard weak form.

	In contrast to the previous example we now refine simultaneously in space and time.
	The $p = \{2, \allowbreak 3, \allowbreak 4 \}$ \ac{PERRK} schemes are combined with \ac{DG} representations with polynomial degrees $k = \{1, \allowbreak 2, \allowbreak 3\}$ to obtain discretizations of matching order of accuracy.
	The domain is refined in $[-1, 1]$ by a factor of two and thus two-member \ac{PERRK} families $E_2 = \{4, \allowbreak 8 \}$, $E_3 = \{5, \allowbreak 9 \}$, $E_4 = \{6, \allowbreak 10 \}$ are constructed for the semidiscretization at hand.
	Note that the considered flow is diffusion-dominated (with Reynolds number $\text{Re} = \frac{\rho_0 \, v \, l}{\mu} = 1.875$) and thus the maximum stable timestep scales quadratically with the smallest mesh size.
	We start the convergence study with a total of $N=6$ cells and refine the mesh five times up to $N=192$.
	The timestep is set for all $p$ to $\Delta t = 0.035 \frac{\Delta x_\text{min}^2}{v}$.
	To solve the relaxation equation \cref{eq:RelaxationEquation} we resort to the bisection procedure with parameters $\gamma_{\text{min}} = 0.8, \gamma_{\text{max}} = 1.2$.
	The root search is stopped at $25$ iterations or if either the length of the bracketing interval or the residual of the relaxation equation is less or equal to $10^{-15}$.

\begin{figure}
		\centering
		\subfloat[{Domain-normalized $L^1$-error in momentum $\rho u$, cf. \cref{eq:L1ErrorDomainNormalized_Density}.}]{
			\label{fig:Convergence_VSP_L1}
			\centering
			\resizebox{.47\textwidth}{!}{\includegraphics{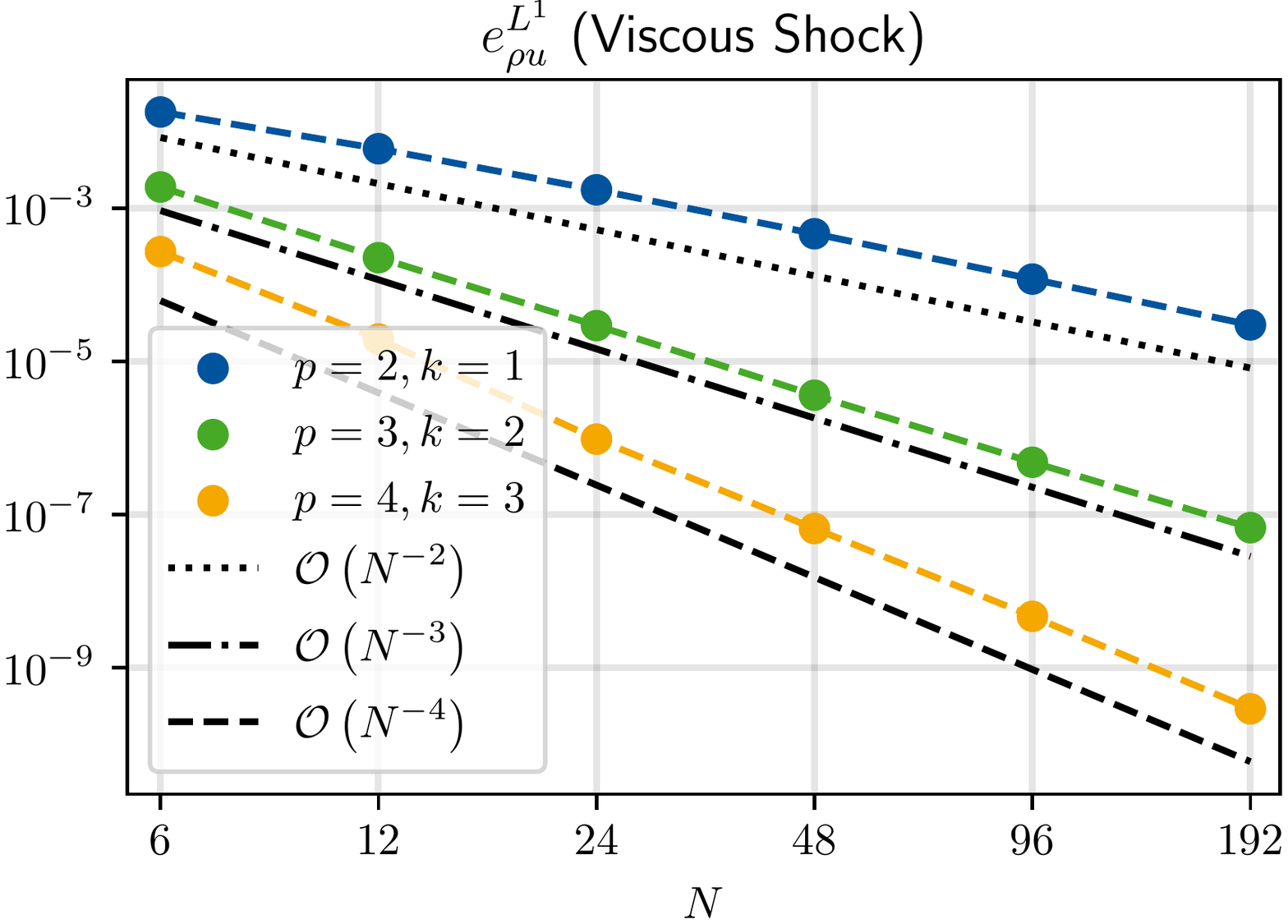}}
		}
		\hfill
		\subfloat[{$L^\infty$ error in momentum $\rho u$.}]{
			\label{fig:Convergence_VSP_LInf}
			\centering
			\resizebox{.47\textwidth}{!}{\includegraphics{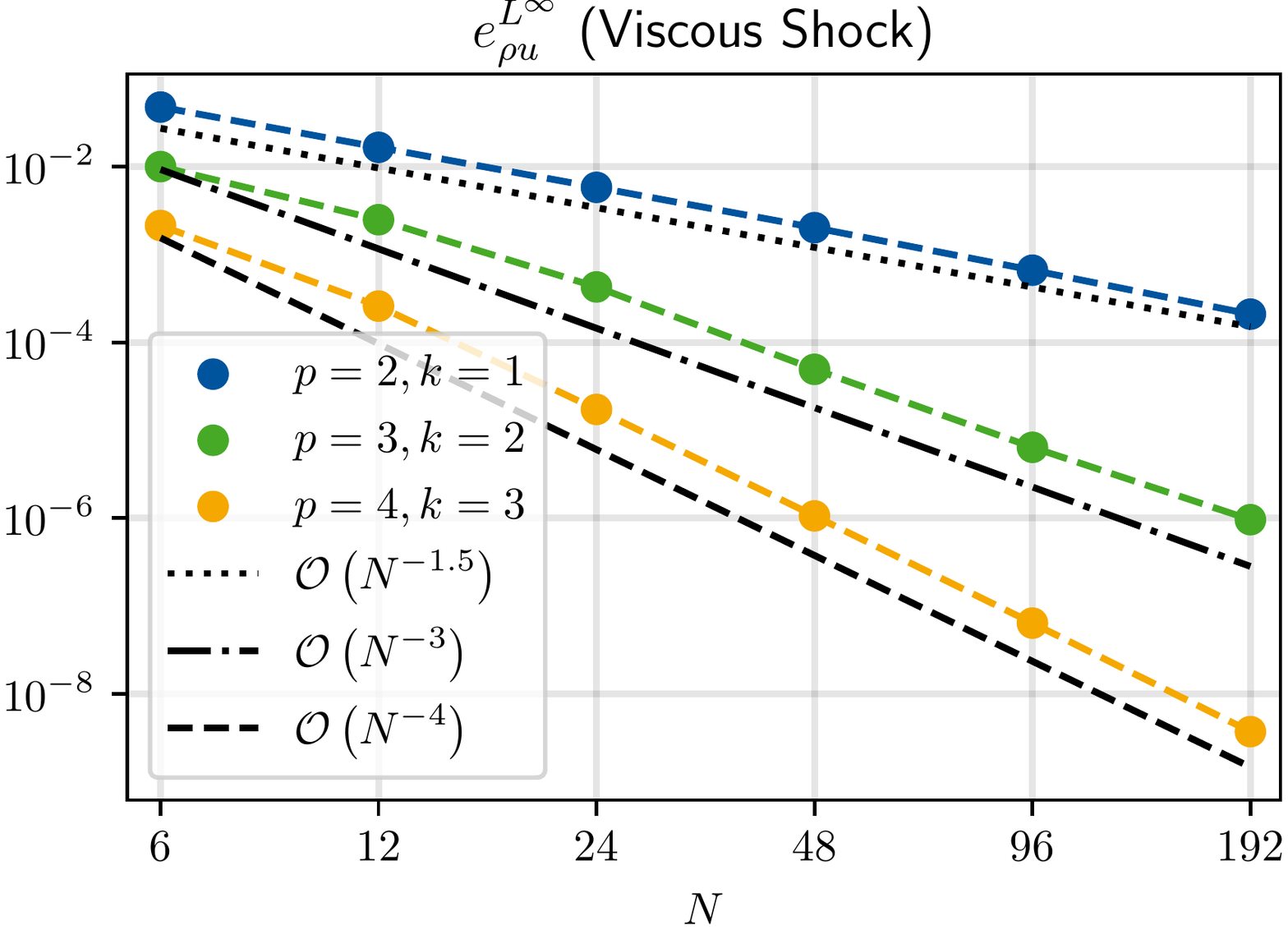}}
		}
		\caption[$L^1$ and $L^\infty$ errors of the momentum $\rho u$ for the viscous shock testcase.]
		{$L^1$ (\cref{fig:Convergence_IVA_L1}) and $L^\infty$ (\cref{fig:Convergence_IVA_LInf}) errors of the density $\rho$ for the viscous shock propagation testcase at $t_f = 0.5$ for $p=2, 3, 4$ \ac{PERRK} schemes.}
		\label{fig:Convergence_VSP_L1_LInf}
	\end{figure}

	$L^1$ and $L^\infty$ errors of the momentum $\rho u$ are presented in \cref{fig:Convergence_VSP_L1_LInf}.
	We observe the designed order of convergence for all schemes, except for the $L^\infty$ error of the $p=2$ which convergence with order $1.5$.
	This is, however, not due to the \ac{PERRK} scheme but instead a consequence of the linear solution polynomials.
	If $k=2$ degree solution polynomials are employed second order convergence in $e^{L^\infty}_{\rho u}$ is observed also for the second-order method.
	We remark that essentially the same convergence results are obtained with the Newton-Raphson method.
	We also examine the evolution of the relaxation parameter $\gamma_n$ for the viscous shock testcase for $N = 24$ grid cells.
	In this case, the deviations of $\gamma_n$ are of order $\mathcal O \left( 10^{-5} \right)$ for $p= 2$ and for both $p=3, p = 4$ of same magnitude $\mathcal O \left( 10^{-7} \right)$, see \cref{fig:gamma_VSP_p2p3p4} in the appendix.
	\subsubsection{Alfvén Wave}
	\label{subsec:ConvTest_AlfvenWave}
	As the final validation study we consider the \ac{MHD} equations in two spatial dimensions.
	The divergence constraint for the magnetic field 
	\begin{equation*}
		\tag*{\cref{eq:DivConstraint}}
		\nabla \cdot \boldsymbol B = 0	
	\end{equation*}
	is enforced using the generalized Lagrange multiplier method \cite{munz2000divergence, dedner2002hyperbolic, derigs2018ideal, bohm2020entropy}.
	The
	\ac{MHD} equations with thermodynamically consistent divergence cleaning using the \ac{GLM} methodology \cite{munz2000divergence, dedner2002hyperbolic, derigs2018ideal} read in $d$ spatial dimensions \cite{bohm2020entropy, warburton1999discontinuous}
	\begin{equation}
		\label{eq:MHDFluxForm}
			\partial_t \boldsymbol u 
		+ \left[ \sum_{i=1}^d  \partial_{i} \boldsymbol f_i(\boldsymbol u)  \right]+ \boldsymbol g(\boldsymbol u, \nabla \boldsymbol u) = \boldsymbol 0
	\end{equation}
	with unknown conserved variables $\boldsymbol u = (\rho, \rho \boldsymbol v, E, \boldsymbol B, \psi)$.
	The last variable $\psi$ is due to the \ac{GLM} extension of the system and corresponds to the propagation of the divergence errors.
	For a two-dimensional problem the inviscid fluxes in $x$ and $y$ direction read
	\begin{subequations}
		\label{eq:MHD_InviscidFluxes}
		\begin{align}
			\boldsymbol f_x(\boldsymbol u) &= 
			\begin{pmatrix}
				\rho v_x \\
				\rho v_x^2 + p + E_{\text{mag}} - B_x^2 \\
				\rho v_x v_y - B_x B_y \\
				v_x \left[ E_{\text{kin}} + \frac{\gamma}{\gamma - 1} p + 2 E_{\text{mag}} \right] - B_x \left[ \boldsymbol v \cdot \boldsymbol B \right] \\
				0 \\
				v_x B_y - v_y B_x \\ 
				c_h B_x
			\end{pmatrix}, \\
			\boldsymbol f_y(\boldsymbol u) &=
			\begin{pmatrix}
				\rho v_y \\
				\rho v_x v_y - B_x B_y \\
				\rho v_y^2 + p + E_{\text{mag}} - B_y^2 \\
				v_y \left[ E_{\text{kin}} + \frac{\gamma}{\gamma - 1} p + 2 E_{\text{mag}} \right] - B_y \left[ \boldsymbol v \cdot \boldsymbol B \right] \\
				v_y B_x - v_x B_y \\
				0 \\
				c_h B_y
			\end{pmatrix}
		\end{align}
	\end{subequations}
	with kinetic $E_{\text{kin}} \coloneqq \frac{1}{2} \rho \boldsymbol v^2$ and magnetic energy $E_{\text{mag}} = \frac{1}{2} \boldsymbol B^2$.
	The pressure $p$ is computed via $p = (\gamma - 1) \left[ E - E_{\text{kin}} - E_{\text{mag}} - \frac{1}{2} \psi^2\right]$ where $E = \rho e$ is the total energy.
	In the equations above, $c_h$ is the divergence cleaning speed which is computed from a \ac{CFL}-like condition.
	The nonconservative Powell term $\boldsymbol g(\boldsymbol u, \nabla \boldsymbol u)$ is given by \cite{derigs2018ideal, bohm2020entropy, powell1999solution}
	\begin{equation}
		\label{eq:MHD_PowellTerm}
		\boldsymbol g = \left(\nabla \cdot \boldsymbol B \right) 
		\begin{pmatrix} 
			0 \\ B_x \\ B_y \\ \boldsymbol v \cdot \boldsymbol B \\ v_x \\ v_y \\ 0
		\end{pmatrix}
		+ \partial_x \psi 
		\begin{pmatrix}
			0 \\ 0 \\ 0 \\ v_x \psi \\ 0 \\ 0 \\ v_x
		\end{pmatrix}
		+ \partial_y \psi
		\begin{pmatrix}
			0 \\ 0 \\ 0 \\ v_y \psi \\ 0 \\ 0 \\ v_y
		\end{pmatrix}
	\end{equation}
	which vanishes in the (continuous or one-dimensonal) divergence-free case.

	To test convergence we consider the propagation of an Alfvén wave in the two-dimensional domain $\Omega = [0, \sqrt{2}]^2$.
	The solution for all $t \in \mathbb N$ is given by \cite{TOTH2000605, derigs2016novel}
	\begin{equation}
		\boldsymbol u_\text{Prim} = 
		\begin{pmatrix}
			\rho \\
			v_x \\
			v_y \\
			p \\
			B_x \\
			B_y \\
			\psi
		\end{pmatrix}	
		= 
		\begin{pmatrix}
			1 \\
			\frac{-0.1}{\sqrt{2}} \sin \Big( \pi \sqrt{2} (x + y) \Big) \\
			\frac{0.1}{\sqrt{2}} \sin \Big( \pi \sqrt{2} (x + y) \Big) \\
			0.1 \\
			\frac{1}{\sqrt{2}} - \frac{0.1}{\sqrt{2}} \sin \Big( \pi \sqrt{2} (x + y) \Big) \\
			\frac{1}{\sqrt{2}} + \frac{0.1}{\sqrt{2}} \sin \Big( \pi \sqrt{2} (x + y) \Big) \\
			0
		\end{pmatrix}
		\: .
	\end{equation}
	Here, we run the simulations up to $t_f = 2$.
	The isentropic exponent is set to $\gamma =  \sfrac{5}{3}$.
	To discretize \cref{eq:MHDFluxForm} we resort to the flux-differencing with separate fluxes for surface fluxes and volume fluxes.
	For these, again different fluxes are used for the conserved and inviscid fluxes.
	To begin, the conserved interface flux is simply the central flux, while for the non-conservative Powell term we use the specialized flux described in \cite{derigs2018ideal, bohm2020entropy}.
	For the conservative volume fluxes the Rusanov/Local Lax-Friedrichs flux \cite{RUSANOV1962304} is employed while the non-conservative fluxes are discretized analogous to the surface fluxes.

	Similar to the previous example simultaneous refinement in space and time is performed.
	Here, however, the spatial discretization is always of increased accuracy, i.e., we choose $k=2$ polynomials for the $p=2$ \ac{PERRK} scheme, $k=3$ polynomials for the $p=3$ scheme, and $k=4$ polynomials for the $p=4$ scheme.
	This way we ensure that the temporal errors should dominate the spatial errors in the asymptotic limit.
	For every order of consistency, $E = \{10, \allowbreak 11, \allowbreak 12, \allowbreak 13, \allowbreak 14\}$ \ac{PERRK} schemes are constructed which are 
	randomly assigned to the uniform grid cells, similar to the testcase for preservation of linear invariants \cref{subsec:ConservationOfLinearInvariants}.
	The relaxation equation is solved with the Newton-Raphson method with a maximum of $10$ iterations, root-residual tolerance and minimum step tolerance of $10^{-15}$.

\begin{figure}
		\centering
		\subfloat[{Domain-normalized $L^1$-error in $B_x$, cf. \cref{eq:L1ErrorDomainNormalized_Density}.}]{
			\label{fig:Convergence_AW_L1}
			\centering
			\resizebox{.47\textwidth}{!}{\includegraphics{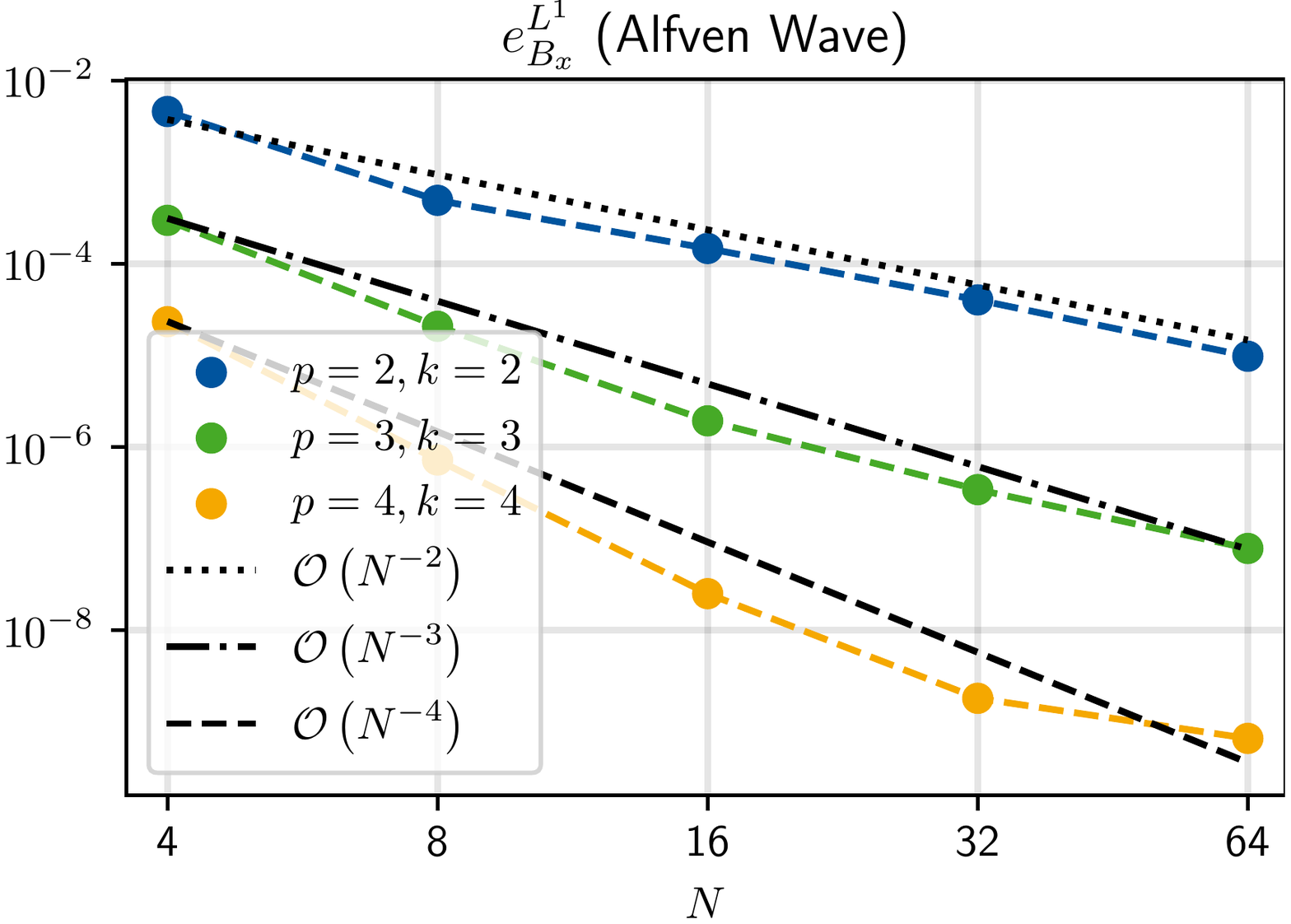}}
		}
		\hfill
		\subfloat[{$L^\infty$ error in $B_x$.}]{
			\label{fig:Convergence_AW_LInf}
			\centering
			\resizebox{.47\textwidth}{!}{\includegraphics{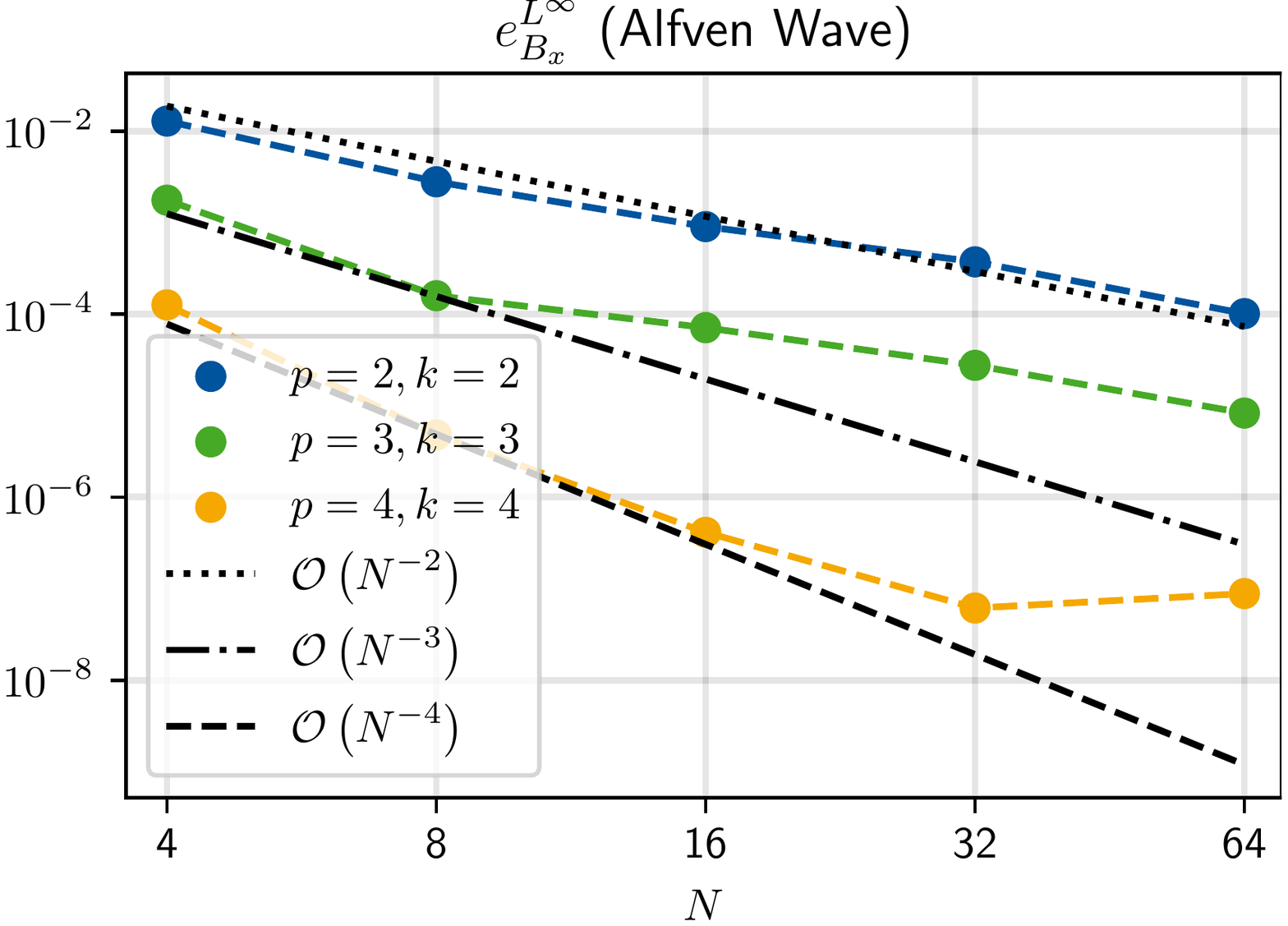}}
		}
		\caption[$L^1$ and $L^\infty$ errors of the $x$-component of the magnetic field for the Alfvén wave testcase.]
		{$L^1$ (\cref{fig:Convergence_IVA_L1}) and $L^\infty$ (\cref{fig:Convergence_IVA_LInf}) errors of the $x$-component of the magnetic field $B_x$ for the Alfvén wave testcase at $t_f = 2$ for $p=2, 3, 4$ \ac{PERRK} schemes with random assignment.}
		\label{fig:Convergence_AW_L1_LInf}
	\end{figure}

	As for the previous testcases, the expected order of convergence is observed for all schemes for the $L^1$ \cref{eq:L1ErrorDomainNormalized_Density} error in $x$-component of the magnetic field $B_x$, cf. \cref{fig:Convergence_AW_L1}.
	For the $L^\infty$ error we observe for $p=3$ and $p=4$ a reduced order of convergence which we attribute to the random assignment of the schemes to the grid cells.
	In particular, the constituting schemes of the \ac{PERRK} schemes possess different error constants.
	Since the individual schemes are no longer assigned to the same cells/locations across the convergence study, schemes with different error constants appear at points with different solution features during the study.
	Thus, convergence in $L^\infty$ is not guaranteed, as observed for the $p=3$ scheme and the last refinement of the $p=4$ scheme.
	We remark that this is not an issue inherent to the relaxation approach, as we observe the same behaviour for the standard/non-relaxed schemes.
	Plots of the relaxation parameter $\gamma_n$ for the $N=64^2$ discretization are provided by \cref{fig:gamma_AW_p2p3p4} in the appendix.
	Interestingly, $\gamma$ shows for $p=2$ and $p=3$ a periodic behaviour, in contrast to the chaotic behaviour observed for the other examples and the $p=4$ scheme.

	In addition, we also studied the convergence for a round-robin assignment of the $R=5$ schemes to the grid cells for $t_f = 4$.
	For $R = 5$ schemes this results for $N = 2^n, \allowbreak \, \allowbreak n = 3, \dots , 6$ gridcells in interfaces between schemes for every cell, both in $x$ and in $y$ direction.
	As for the random assignment, we observe the expected order of convergence in $L^1$ while the $L^\infty$ error deviates from the expected asymptotic convergence for the $p = 3$ and $p = 4$ scheme for the last two refinement steps, cf. \cref{fig:Convergence_AW_L1_LInf_RR}.
	As for the random assignment, this is due to the different error constants of the schemes and not an artifact of the relaxation approach, as we confirmed by performing the same convergence study with the standard \ac{PERRK} schemes.
	\section{Applications}
	\label{sec:Applications}
	In this section we apply the \ac{PERRK} methods to a range of different flow regimes of increasing complexity.
	We begin by considering three two-dimensional examples.
	First, we study the laminar, viscous, effectively incompressible flow over the SD7003 airfoil in \cref{subsec:SD7003}.
	This problem can be seen as an extended testcase as it allows validating the \ac{PERRK} methods by comparing aerodynamic coefficients to reference data and is not characterized by complicated flow features.
	Next, we consider the flow of a magnetized, viscous, compressible fluid around a cylinder in \cref{subsec:VRMHD_Cylinder}, which is described by the \ac{vrMHD} equations \cite{bohm2020entropy, warburton1999discontinuous}.
	Then, we turn our attention to the inviscid, transonic flow around the NACA0012 airfoil.
	The testcase discussed in \cref{subsec:NACA0012} corresponds to case C1.2 of the third international workshop on high-order computational fluid dynamics methods held in 2015 \cite{3HOCFD}.
	We apply \ac{AMR} to this problem to resolve the steady standing shock on the airfoil top surface.
	Coming to three-dimensional problems, we consider the classic ONERA M6 wing \cite{onera_m6_original_report, slater_study_1} in \cref{subsec:ONERA_M6} for an inviscid, compressible gas, i.e., also a steady setup.
	Additionally, we investigate the robustness of the relaxed multirate schemes by means of a viscous transonic flow across the NASA \ac{CRM} on a very coarse grid in \cref{subsec:NASA_CRM}.

	All computations are performed using the Trixi.jl \cite{trixi1, trixi2, trixi3} package and can be reproduced using the repository \cite{doehring2025PERRK_ReproRepo}.
	\subsection{Viscous Laminar Flow over SD7003 Airfoil}
	\label{subsec:SD7003}
	\subsubsection{Setup}
	As a first application case, we reconsider the laminar, viscous flow over the SD7003 airfoil presented in \cite{doehring2024fourth}.
	This testcase was also employed in the original \ac{PERK} papers \cite{vermeire2019paired,nasab2022third}.

	The flow parameters for this Navier-Stokes-Fourier simulation are $\text{Ma} = 0.2 = \frac{U_\infty}{c_\infty}$, $\text{Re} = 10^4 = \frac{\rho_\infty U_\infty l}{\mu_\infty}$ with airfoil chord length $l = 1.0$, angle of attack $\alpha = 4 \degree$, non-dimensionalized speed of sound $c_\infty = 1.0$, pressure $p_\infty = 1.0$, and density $\rho_\infty = 1.4$ corresponding to isentropic exponent $\gamma = 1.4$.
	The viscosity is kept constant at $\mu = 2.8 \cdot 10^{-5}$ and the Prandtl number is set to $\text{Pr} = 0.72$.
	At the domain boundaries the freestream flow is weakly enforced \cite{mengaldo2014guide} for both inviscid and viscous fluxes (gradients).
	The airfoil is modeled as a no-slip, adiabatic wall.
	The simulation is initialized with the freestream values.

\begin{figure}
		\centering
		\includegraphics[width=0.75\textwidth]{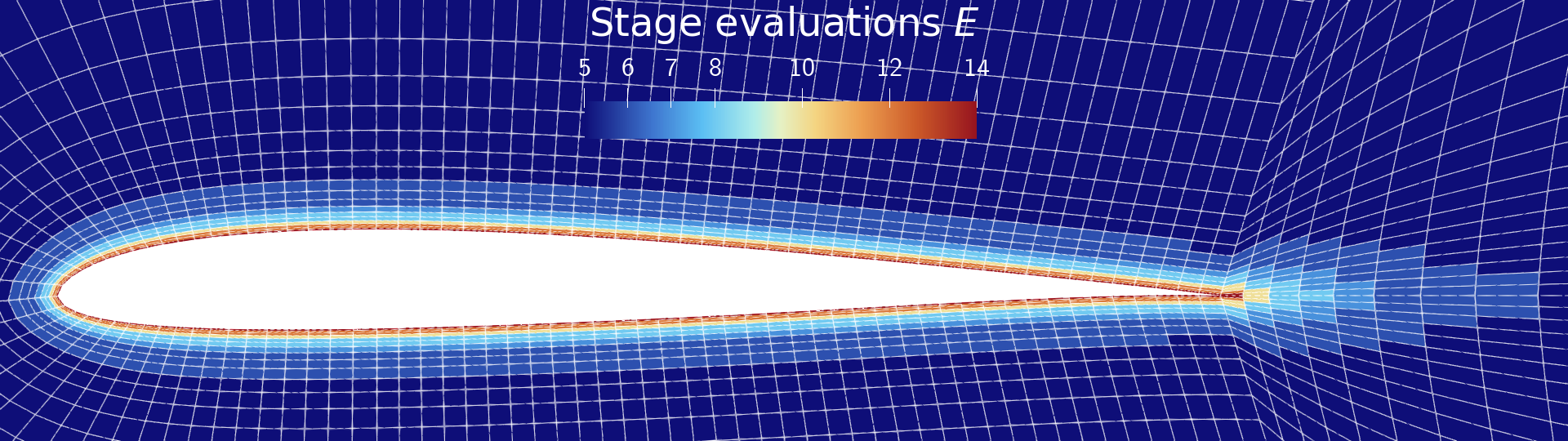}
		\caption[Method distribution for the SD7003 airfoil integration with the $E = \{5, 6, 7, 8, 10, 12, 14\}$ multirate \ac{PERK} family.]
		{Method distribution for the SD7003 airfoil integration with the $E = \{5, 6, 7, 8, 10, 12, 14\}$ multirate \ac{PERK} family.}
		\label{fig:SD7003_StageDistr}
	\end{figure}

	We employ a two-dimensional mesh consisting of $7605$ straight-sided quadrilateral elements.
	This mesh (with second-order elements) has been provided as electronic supplementary material to the second-order \ac{PERK} paper \cite{vermeire2019paired}.
	The spatial derivatives are discretized using the \ac{DGSEM} with $k=3$ polynomials and \ac{HLLC} numerical flux \cite{toro1994restoration} with volume flux differencing \cite{CHEN2017427} and entropy-conservative volume flux \cite{Chandrashekar_2013}.
	The gradients required for the viscous fluxes are computed with the BR1 scheme \cite{bassi1997high, gassner2018br1}.
	We optimize fourth-order, $E = \{5, \allowbreak 6, \allowbreak 7, \allowbreak 8, \allowbreak 10, \allowbreak 12, \allowbreak 14\}$ \ac{PERK}/\ac{PERRK} schemes for the free-stream state.
	The distribution of the composing schemes of the near-airfoil grid is depicted in \cref{fig:SD7003_StageDistr}.
	%

	Lift and drag coefficients $C_L, C_D$ are recorded over the $t \in [30 t_c, 35 t_c]$ interval where $t_c \coloneqq t \, \frac{U_\infty}{c}$ is a non-dimensionalized convective time.
	The provided lift coefficient is solely due to the pressure difference:
	\begin{equation}
		\label{eq:LiftCoeff}
		C_L = C_{L,p} = \oint_{\partial \Omega} \frac{ p \, \boldsymbol n \cdot \boldsymbol t^\perp}{0.5 \rho_{\infty} U_{\infty}^2 L_{\infty}} \, \text{d} S
	\end{equation}
	with $\boldsymbol t^\perp = \big( -\sin(\alpha), \cos(\alpha) \big)$ being the unit direction perpendicular to the free-stream flow.
	The vector $\boldsymbol n$ is the fluid-element outward-pointing normal, i.e., the normal pointing into the surface of the airfoil.
	The drag coefficient incorporates also the viscous stress contribution:
	\begin{subequations}
		\label{eq:DragCoeffs}
		\begin{align}
			C_D &= C_{D,p} + C_{D,\mu} \\
			\label{eq:DragCoeffPressure}
			C_{D,p} &= \oint_{\partial \Omega} \frac{p \, \boldsymbol n \cdot \boldsymbol t}
			{0.5 \rho_{\infty} U_{\infty}^2 L_{\infty}} \, \text{d} S \\
			C_{D,\mu} &= \oint_{\partial \Omega} \frac{\boldsymbol \tau_w \cdot \boldsymbol t}
			{0.5 \rho_{\infty} U_{\infty} L_{\infty}} \, \text{d} S
		\end{align}
	\end{subequations}
	where $\boldsymbol t = \big( \cos(\alpha), \sin(\alpha) \big)$.
	The wall shear stress vector $\boldsymbol \tau_w$ is obtained from contracting the viscous stress tensor $\underline{\underline{\boldsymbol \tau}}$ with the surface normal $\hat{\boldsymbol n} = - \boldsymbol n$:
	\begin{equation}
		\boldsymbol{\tau}_w = \underline{\underline{ \boldsymbol \tau}} \, \hat{\boldsymbol n} = \begin{pmatrix} \tau_{xx} & \tau_{xy} \\ \tau_{yx} & \tau_{yy} \end{pmatrix} \begin{pmatrix} -n_x \\ -n_y \end{pmatrix} \, .
	\end{equation}
	The unsteady lift and drag coefficients, cf. \cref{fig:LiftDragCoeffs},
	are sampled every $\Delta t_c = 5 \cdot 10^{-3}$ and averaged to enable comparison with values from the literature, see \cref{tab:SD7003_Lift_Drag}.
	We observe good agreement with values reported in the literature \cite{vermeire2019paired, nasab2022third, Uranga2011Implicit, LopezMorales2014Verification}.
	In particular, the lift and drag coefficients of the relaxed fourth-order \ac{PERRK} scheme deviate only by $\mathcal{O}\left(10^{-6} \right)$ from the values of the plain fourth-order \ac{PERK} scheme developed in \cite{doehring2024fourth}.
	\begin{table}
		\def\arraystretch{1.2}
		\centering
		\begin{tabular}{l?{2}c|c}
			Source & $\overline{C}_L 
			$ & $\overline{C}_D 
			$ \\
			\Xhline{5\arrayrulewidth}
			\ac{PERRK} $p=4, k=3$ 																   & $0.3827$ & $0.04995$ \\
			\hdashline
			\ac{PERK} $p=2, k=3$ \cite{vermeire2019paired}           & $0.3841$ & $0.04990$ \\
			\ac{PERK} $p=3, k=3$ \cite{nasab2022third}               & $0.3848$ & $0.04910$ \\
			Uranga et al. \cite{Uranga2011Implicit}                  & $0.3755$ & $0.04978$ \\
			López-Morales et al. \cite{LopezMorales2014Verification} & $0.3719$ & $0.04940$ \\
		\end{tabular}
		\caption[Average lift and drag coefficients $\overline{C}_L$ and $\overline{C}_D$ for the SD7003 airfoil.]
		{Time-averaged lift and drag coefficients $\overline{C}_L$ and $\overline{C}_D$ for the SD7003 airfoil over the $[30 t_c, \, 35 t_c]$ interval.}
		\label{tab:SD7003_Lift_Drag}
	\end{table}
	\subsubsection{Performance Comparison}
	To begin, we demonstrate performance improvements of the \ac{PERK} scheme with relaxation over the plain \ac{PERK} scheme.
	For the given configuration, the \ac{CFL} number can be increased by roughly \qty{20}{\percent}, which results in about $25 \cdot 10^3$ fewer total simulation steps, which manifests in a \qty{19}{\percent} decreased runtime, despite the additional overhead due to computation of the relaxation parameter $\gamma$.
	Similar values are obtained for the standalone $S=E=14$ \ac{PERRK} scheme when compared to the non-entropy-stable $S=14$ \ac{PERK} scheme.
	We attribute this observation to the increased robustness of the \ac{PERRK} scheme compared to the standard scheme, cf. \cref{subsec:EntropyRelaxationTimeLimiting}.
	This is particularly important given that the flow is not fully resolved, cf. \cref{fig:SD7003_vy} which results in artifacts at the cell boundaries.

	For the multirate \ac{PERRK} scheme, the relaxation solver makes up for a total of \qty{2.7}{\percent} of the total runtime, out of which \qty{1.9}{\percent} are spent on computing the entropy difference $\Delta H$ \cref{eq:EntropyChange_EntropyVariables_ReArranged}.
	This is a relatively expensive operation as it involves integration over the physical domain, which in turn demands the transform from reference to physical space for every element.
	Here, some code optimization could be possible e.g. by embodying this integration in the \ac{DG} solution process.

	Second, we compare the performance of the entropy-conservative, multirate scheme with entropy-stable, standalone schemes.
	Besides the natural comparison with the optimized $S=E=14$ \ac{PERRK} scheme, we consider three other fourth-order accurate Runge-Kutta schemes which permit a low-storage implementation in Butcher form.
	This guarantees a fair comparison of the \ac{PERRK} schemes to other entropy-stable methods.
	In particular, we consider two "purely sub-diagonal" schemes, i.e., methods with Butcher tableau
	\begin{equation}
		\label{eq:ButcherTableauSubDiagonal}
		\renewcommand\arraystretch{1.3}
		\begin{array}
				{c|c|c c c c}
				i & \boldsymbol c & \multicolumn{4}{c}{A} \\
				\hline
				1 & 0 & & & & \\
				2 & c_2 & c_2 & & &   \\
				3 & c_3 & 0 & c_3 & &  \\ 
				\vdots & \vdots & \vdots & \ddots & \ddots &  \\
				S & c_S & 0 & \dots & 0 & c_S \\
				\hline
				& & \multicolumn{4}{c}{ \boldsymbol b^T }
		\end{array}
		 \qquad .
	\end{equation}
	The classic four-stage, fourth-order Runge-Kutta schemes can be classified as a purely sub-diagonal scheme of form \eqref{eq:ButcherTableauSubDiagonal}.
	Furthermore, the six-stage, fourth-order low-dispersion, low-dissipation scheme proposed in \cite{tselios2005runge} is also of this form.

	In addition, we consider the five-stage, fourth-order Runge-Kutta scheme optimized for compressible Navier-Stokes equations from \cite{kennedy2000low}.
	The method considered here is implemented in van-der-Houwen form \cite{van1972explicit}, i.e., with Butcher tableau archetype
	\begin{equation}
		\label{eq:ButcherTableau_vdH2N}
		\renewcommand\arraystretch{1.3}
		\begin{array}
				{c|c|c c c c c c}
				i & \boldsymbol c & \multicolumn{5}{c}{A} & \\
				\hline
				1 & 0 & & & & & & \\
				2 & c_2 & a_{21} & & & & & \\
				3 & c_3 & b_1 & a_{32} & & & & \\ 
				4 & c_4 & b_1 & b_2 & a_{43} & & & \\ 
				\vdots & \vdots & \vdots & \vdots & \ddots & \ddots & & \\
				S & c_S & b_1 & b_2 & \dots & b_{S-2} & a_{S, S-1} & \\
				\hline
				& & b_1 & b_2 & \dots & b_{S-2} & b_{S-1} & b_S
		\end{array}
		\qquad .
	\end{equation}
	For all considered methods, the relaxation equation \cref{eq:RelaxationEquation} is solved with the Newton-Raphson method which stops after five iterations, for objective residual less or equal $10^{-13}$ and stepsizes less or equal to $10^{-13}$.
	This set of parameters is motivated by our findings in section \cref{sec:EntropyConservation}, in particular the two-dimensional examples.
	For this configuration, very few Newton iterations are required to reach the desired accuracy which mainly conenctrate on the first few time steps after the initialization.

	We compare the observed ratios of measured runtime to optimal (i.e., \ac{PERRK}-Multi) runtime $\sfrac{\tau}{\tau^\star}$ and scalar (i.e., per \ac{DoF}) \ac{RHS} evaluations
	\begin{equation}
		\label{eq:RHSEvals}
		N_\text{RHS} = N_t \cdot N_K
	\end{equation}
	for the aforementioned schemes in \cref{tab:Runtimes_SD7003}.
	In \eqref{eq:RHSEvals}, $N_t$ the number of timesteps taken and $N_K$ is the number of scalar \ac{RHS} evaluations involving all stages.
	For the multirate schemes the latter is computed as 
	\begin{equation}
		N_K = \sum_{r=1}^R E^{(r)} \cdot M^{(r)}
	\end{equation}
	where $M^{(r)}$ is the number of \acp{DoF} integrated with the $r$-th partition.
	For standard schemes, $N_K = S \cdot M$ where $M$ the dimension of the semidiscretization \eqref{eq:Semidiscretization} and $S$ is the number of stages.
	The values in \cref{tab:SD7003_Lift_Drag} correspond to the run from $t = 0$ to $t = 30 t_c$ to eliminate skewed results due to the frequent computation of the lift and drag coefficients.

	The multirate \ac{PERRK} scheme is the most efficient method both in terms of runtime and number of scalar \ac{RHS} evaluations.
	In particular, the \ac{PERRK}-Multi scheme is twice as fast as the optimized $p=4, S = E =14$ scheme, which is the most efficient standalone scheme.
	The number of \ac{RHS} calls is reduced by a factor of $2.42$, i.e., the \ac{PERRK}-Multi scheme requires roughly \qty{42}{\percent} of the scalar \ac{RHS} evaluations of the optimized standalone scheme.
	The discrepancy between the saved computations and runtime is a consequence of the bookkeeping overhead for the multirate scheme and the loss of data locality which leads to suboptimal performance due to, e.g., cache misses.
	The other considered schemes require at least $2.5$ times scalar \ac{RHS} evaluations and are at least $2.3$ times slower than the optimized \ac{PERRK}-Multi scheme.
	\begin{table}
		\def\arraystretch{1.2}
		\centering
		\begin{tabular}{l?{2}c|c}
			Method & $\sfrac{\tau}{\tau^\star}$ & $\sfrac{N_\text{RHS}}{N_\text{RHS}^\star}$ \\
			\Xhline{5\arrayrulewidth}
			$\text{P-ERRK}_{4;\{5,6, \dots, 14\}}$ 			& $1.0$  & $1.0$ \\
			$\text{P-ERRK}_{4;14}$ 											& $2.13$ & $2.42$ \\
			\hdashline
			$\text{P-ERK}_{4;\{5,6, \dots, 14\}}$ 			& $1.19$  & $1.20$ \\
			\hdashline
			$\text{R-RK}_{4;4}$ 												& $2.55$ & $2.77$ \\
			$\text{R-TS}_{4;6}$ \cite{tselios2005runge} & $2.98$ & $3.16$ \\
			$\text{R-CKL}_{4;5}$ \cite{kennedy2000low}  & $2.38$ & $2.53$
		\end{tabular}
		\caption[Runtimes and number of scalar RHS evaluations of different fourth-order integrators compared to the optimized $p = 4, E = \{5, \allowbreak 6, \allowbreak \dots, \allowbreak 14\} $ integrator for the laminar flow around the SD7003 airfoil.]
		{Runtimes and number of scalar \ac{RHS} evaluations \cref{eq:RHSEvals} of different entropy conservation, relaxed fourth-order integrators compared to the optimized $p = 4, E = \{5, \allowbreak 6, \allowbreak 7, \allowbreak 8, \allowbreak 10, \allowbreak 12, \allowbreak 14\} $ integrator for the laminar flow around the SD7003 airfoil.
		We also included the non-relaxed $\text{P-ERK}_{4;\{5,6, \dots, 14\}}$ scheme.}
		\label{tab:Runtimes_SD7003}
	\end{table}
	\subsection{Visco-Resistive MHD Flow past a Cylinder}
	\label{subsec:VRMHD_Cylinder}
	In this example we consider the \ac{vrMHD} equations in two spatial dimensions.
	The
	\ac{vrMHD} equations with 
	divergence cleaning using the \ac{GLM} methodology \cite{munz2000divergence, dedner2002hyperbolic, derigs2018ideal} read in $d$ spatial dimensions \cite{bohm2020entropy, warburton1999discontinuous} (cf. \cref{eq:MHDFluxForm})
	\begin{equation}
		\label{eq:VRMHDFluxForm}
			\partial_t \boldsymbol u 
		+ \left[ \sum_{i=1}^d  \partial_{i} \big( \boldsymbol f_i(\boldsymbol u) - \boldsymbol f_i^{\text{vr}} (\boldsymbol u, \nabla \boldsymbol u) \big) \right]+ \boldsymbol g(\boldsymbol u, \nabla \boldsymbol u) = \boldsymbol 0
	\end{equation}
	with unknown conserved variables $\boldsymbol u = (\rho, \rho \boldsymbol v, E, \boldsymbol B, \psi)$.
	The last variable $\psi$ is due to the \ac{GLM} extension of the system and corresponds to the propagation of the divergence errors.
	The visco-resistive fluxes $\boldsymbol f_i^{\text{vr}}$ are given by \cite{bohm2020entropy, warburton1999discontinuous}
	\begin{subequations}
		\begin{align}
			\boldsymbol f_x^{\text{vr}}(\boldsymbol u, \nabla \boldsymbol u) &= 
			\begin{pmatrix}
				0 \\
				\mu \tau_{xx} \\
				\mu \tau_{xy} \\
				\mu \left[ \boldsymbol v \cdot \boldsymbol \tau_{x} - q_x \right] + 
				\eta B_y \left[\partial_x B_y - \partial_y B_x \right]  \\
				0 \\
				\eta \left[ \partial_x B_y - \partial_y B_x \right] \\
				0
			\end{pmatrix},
			\\
			\boldsymbol f_y^{\text{vr}}(\boldsymbol u, \nabla \boldsymbol u) &=
			\begin{pmatrix}
				0 \\
				\mu \tau_{yx} \\
				\mu \tau_{yy} \\
				\mu \left[ \boldsymbol v \cdot \boldsymbol \tau_{y} - q_y \right] + 
				\eta B_x \left[ \partial_y B_x - \partial_x B_y \right] \\
				\eta \left[\partial_y B_x - \partial_x B_y \right] \\
				0 \\
				0
			\end{pmatrix}
			\: .
		\end{align}
	\end{subequations}
	In the definitions above, $\mu$ is the dynamic viscosity and $\eta$ the magnetic resistivity.
	The symmetric stress tensor $\underline{\underline{\boldsymbol \tau}}$ is given by
	\begin{equation}
		\def\arraystretch{1.4}
		\underline{\underline{\boldsymbol \tau}} \coloneqq \begin{pmatrix}
			\boldsymbol \tau_x \\ \boldsymbol \tau_y
		\end{pmatrix}
		= \begin{pmatrix}
			\frac{4}{3} \partial_x v_x - \frac{2}{3} \partial_y v_y & \partial_x v_y + \partial_y v_x \\
			\partial_y v_x + \partial_x v_y & \frac{4}{3} \partial_y v_y - \frac{2}{3} \partial_x v_x
		\end{pmatrix} \: .
	\end{equation}
	The heat flux $\boldsymbol q$ is given by Fourier's law $\boldsymbol q = - \kappa \nabla T$ with thermal conductivity $\kappa = \frac{\gamma}{\gamma -1} \frac{1}{\text{Pr}}$ and temperature $T$.
	The Prandtl number $\text{Pr}$ is a parameter of the considered gas and treated here as a constant with value $\text{Pr} = 0.72$.
	The conserved fluxes $\boldsymbol f_i$ are given by \cref{eq:MHD_InviscidFluxes} and the nonconservative Powell term $\boldsymbol g$ is given by \cref{eq:MHD_PowellTerm}.
	\subsubsection{Setup}
	The flow of a magnetized fluid past a cylinder has been studied both numerically and experimentally \cite{crawford1995control, weier1998experiments, kanaris2013three}.
	The configuration employed here is loosely based on the setup described in \cite{warburton1999discontinuous}.
	Most importantly, we also set the Mach number $\text{Ma} = 0.5$, rendering the flow truly compressible.
	The ratio of specific heats is set to $\gamma = \frac{5}{3}$ and the Reynolds number based on the diameter of the cylinder $D=1$ is set to $\text{Re} = \frac{\rho_\infty v_\infty D}{\mu} = 200$.
	Following \cite{warburton1999discontinuous}, the Alfvénic Mach number $A = \frac{v_A}{v_\infty}$ is set to $A = 0.1$.
	This defines the magnitude of the magnetic field $B_\infty = v_A \sqrt{\rho_\infty}$.
	With this Alfvénic Mach number the viscous Lundquist number is $S_\mu = A \cdot \text{Re} = \frac{\rho_\infty v_A D}{\mu} = 50$.
	We choose $v_\infty = p_\infty = 1$, thus $\rho_\infty$ follows from the speed of sound $c = \frac{v_\infty}{\text{Ma}}$ of an ideal gas as $\rho_\infty = p_\infty \frac{\gamma}{c^2}$.
	Again, as in \cite{warburton1999discontinuous} we use the same values for the viscous and resistive Lundquist numbers $S_\mu = S_\eta$.
	Then, the magnetic resistivity is given by $\eta = \frac{v_A D}{S_\eta}$.

	The cylinder is placed at $(0,0)$ in the rectangular domain $\Omega = [-5D, 60D] \times [-10D, 10D]$.
	The computational mesh is generated using \texttt{HOHQMesh} \cite{kopriva2024hohqmesh} which enables the generation of conforming high-order, i.e., curved-boundary quadrilateral meshes.
	In this case, the mesh uses fourth-order boundaries, consists of $2258$ elements and is symmetric around the $y=0$ axis, cf. \cref{fig:Cylinder_Mesh}.

	The hyperbolic boundary conditions imposed at $\partial \Omega$ are weakly enforced Dirichlet conditions \cite{mengaldo2014guide} based on the initial/freestream condition.
	At the cylinder surface a slip-wall boundary condition \cite{Vegt2002SlipFB} with zero magnetic field is imposed.
	The parabolic boundary conditions of the cylinder are no-slip, adiabatic, and isomagnetic $\boldsymbol B = \boldsymbol 0$ walls.
	At the inlet, i.e., left boundary, the parabolic boundary condition imposes a moving, adiabatic, isomagnetic fluid with velocity $v_x = v_\infty$ and purely horizontal magnetic field $B_x = B_\infty$.
	On the other boundaries, i.e., top, bottom and outlet/right boundary the gradients are just copied which corresponds to domain-extending boundary conditions.
	The simulation is initialized with the uniform inflow state.

	The solution fields are approximated with $k=3$ \ac{DG} solution polynomials.
	The Riemann problems due to the inviscid fluxes \eqref{eq:MHD_InviscidFluxes} are approximately solved with the \ac{HLL} two-wave solver \cite{harten1983upstream} with Davis-type wave speed estimates \cite{doi:10.1137/0909030}.
	As for the previous example, the viscous fluxes are computed with the BR1 scheme \cite{bassi1997high, gassner2018br1}.
	The nonconservative Powell term \eqref{eq:MHD_PowellTerm} is solved with the specialized numerical flux described in \cite{derigs2018ideal, bohm2020entropy}.
	Opposed to the standard Navier-Stokes-Fourier equations the simulation of a magnetized fluid requires stabilization.
	To this end, we employ the entropy-stable schock-capturing technique by means of the first-order finite volume method on the subcells of the \ac{DG} scheme \cite{hennemann2021provably, rueda2021entropy}.
	The required volume fluxes are the entropy conservation and kinetic energy preserving flux \cite{hindenlang2019new} which extends the flux by Ranocha \cite{Ranocha2020Entropy} to the MHD equations.
	The nonconservative volume fluxes are treated with the same numerical flux as the nonconservative surface fluxes.

	We simulate the flow until $t_f = 120$ such that the vortex street is fully developed and reaches the outlet boundary.
	Upon inspection of the mesh (cf. \cref{fig:Cylinder_Mesh}) it becomes apparent that roughly four base cell sizes (minimum edge length) are present: Cells with $h \approx 0.25$ at the cylinder, $h \approx 0.5$ in the near wake, $h \approx 1$ in the far wake, and $h \approx 2$ at the top and bottom boundary.
	Thus, we seek to find a four-member \ac{PERRK} family for which the admissible timestep roughly doubles for each family member.
	For the fourth-order methods the $E= \{5, \allowbreak 6, \allowbreak 8, \allowbreak 13 \}$ scheme fulfills this requirement and is thus used.

	To solve the relaxation equation \cref{eq:RelaxationEquation} we employ the Newton-Raphson procedure which stops at $3$ iterations or if the residual tolerance is less or equal $10^{-14}$.
	For this simulation we employ a linear ramp-up of the $\text{CFL}$ number, cf. \eqref{eq:TimestepConstraintGlobal}, from a small value $\text{CFL}_0$ which is required to overcome the naive initial condition to a maximum value for the developed flow $\text{CFL}_\text{max}$.
	Thus, the $\text{CFL}$ number as a function of time $t$ is given by
	\begin{equation}
		\label{eq:CFLRampUp}
		\text{CFL}(t) = \min \left \{ \text{CFL}_\text{max}, \, \text{CFL}_0 + \left( \text{CFL}_\text{max} - \text{CFL}_0 \right)\frac{t}{t_\text{CFL}} \right \}	
	\end{equation}
	where the ramp-up time $t_\text{CFL}$ determines how quickly the maximum $\text{CFL}$ number is reached.
	The maximum CFL numbers $\text{CFL}_\text{max}$ have been obtained by starting the simulations from a physical initial condition at $t=10$ and increasing the $\text{CFL}$ number until the simulation produces unphysical results.

\begin{figure}
		\centering
		\subfloat[{Navier-Stokes-Fourier simulation.}]{
			\label{fig:Density_Cylinder_NSF}
			\centering
			\includegraphics[width=0.75\textwidth]{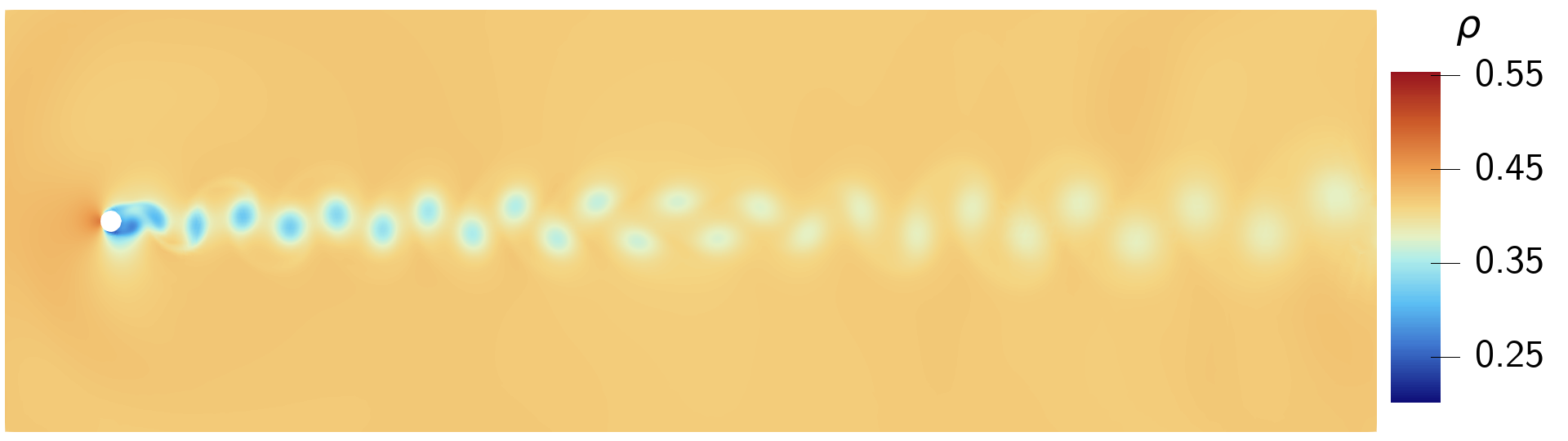}
		}
		\hfill
		\subfloat[{Visco-Resistive MHD (VRMHD) simulation.}]{
			\label{fig:Density_Cylinder_VRMHD}
			\centering
			\includegraphics[width=0.75\textwidth]{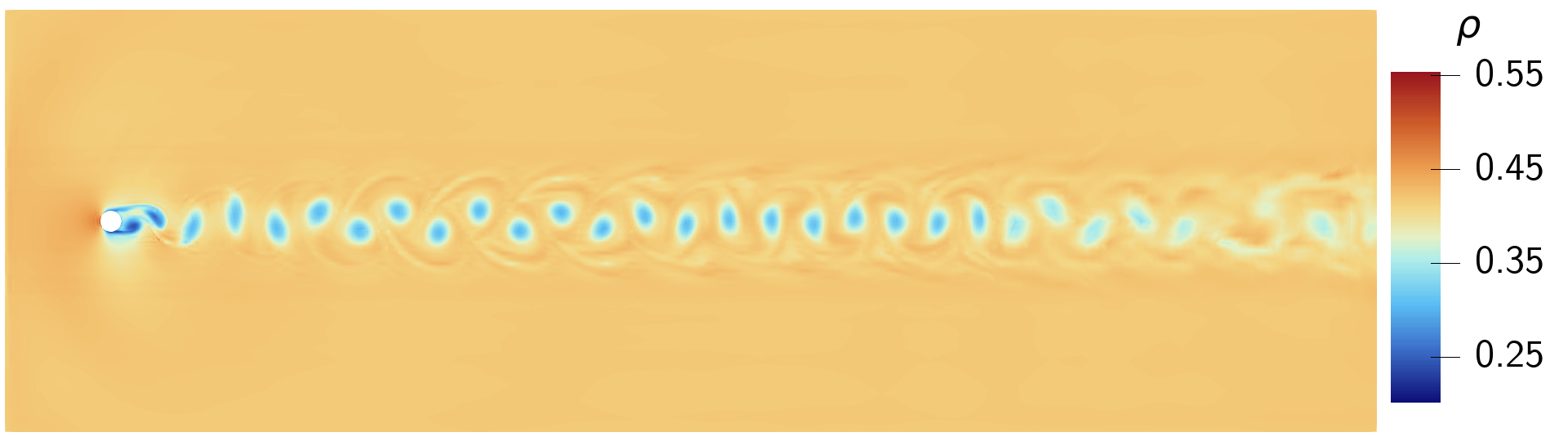}
		}
		\caption[Density $\rho$ at final time $t_f = 120$ for the flow past a cylinder.]
		{Density $\rho$ at final time $t_f = 120$ for the flow past a cylinder.
		In \cref{fig:Density_Cylinder_NSF} the Navier-Stokes-Fourier simulation is shown, while \cref{fig:Density_Cylinder_VRMHD} depicts the visco-resistive MHD simulation.
		Both results are obtained using fourth-order \ac{PERRK} schemes optimized for the respective semidiscretization.}
		\label{fig:Density_Cylinder_NSF_VRMHD}
	\end{figure}

	We plot the density $\rho$ at final time $t_f = 120$ for a standard Navier-Stokes-Fourier simulation and the visco-resistive MHD simulation in \cref{fig:Density_Cylinder_NSF_VRMHD}.
	As observed for instance in \cite{warburton1999discontinuous} the magnetic field alters the regular vortex structures in the wake of the cylinder.
	This is in contrast to the plain diffusion of the vortical structures observed in the Navier-Stokes-Fourier simulation.
	\subsubsection{Performance Comparison}
	For this example, the \ac{PERRK} does not come with improved \ac{CFL} numbers compared to the standard \ac{PERK}.
	The additional work due to the relaxation process accounts for \qty{3.6}{\percent} of the total runtime.
	Additionally, the multirate \ac{PERRK} scheme requires one more timestep ($4537$ vs. $4536$ steps corresponding to the base scheme.)
	Across schemes, we observe that on average a single Newton iteration is required to solve the relaxation equation \cref{eq:RelaxationEquation} to the desired accuracy.

	Turning to the comparison to the standalone relaxed entropy-stable schemes we observe that these require at least \qty{76}{\percent} more scalar \ac{RHS} evaluations, manifesting in \qty{81}{\percent} to \qty{207}{\percent} longer runtimes.
	As for the previous example, the five-stage scheme $\text{CKL}_{4;5}$ \cite{kennedy2000low} is the most efficient out of the non-optimized, off-the-shelf schemes.

	Upon closer inspection of the results it is evident that the off-the-shelf schemes are overproportionally slow in relation to the number of \ac{RHS} evaluations.
	This is partly due to the dense weight vector $\boldsymbol b$ of these schemes, which requires in turn more integrations to compute the entropy change $\Delta H$, cf. \cref{eq:EntropyChange_EntropyVariables_ReArranged}.
	Note that the $p=4$ \ac{PERRK} schemes have only non-zero $b_1$ and $b_S$, independent of the number of stages $S$ \cite{doehring2024fourth}.
	In the present implementation, each stage is integrated on-the-fly, resulting for the standard schemes in $S$ integrations per total timestep.
	This could be reduced, however, by accumulating the integrand of \cref{eq:EntropyChange_EntropyVariables_ReArranged} in an additional vector.
	This would then store the scalar product of the stage values $\boldsymbol W_{n,i}$ and stages $\boldsymbol K_i$ evaluated at every \ac{DG} node $j$.
	In particular, one can not immediately accumulate each elements' contribution into a single scalar sum as in this case information on the geometry is lost which is required to compute the integral.
	Nevertheless, computing the integrand makes up for about only \qty{5}{\percent} of the total runtime which could be reduced to $1-2 \%$ with the aforementioned optimization.
	\begin{table}
		\def\arraystretch{1.2}
		\centering
		\begin{tabular}{l?{2}c|c}
			Method & $\sfrac{\tau}{\tau^\star}$ & $\sfrac{N_\text{RHS}}{N_\text{RHS}^\star}$ \\
			\Xhline{5\arrayrulewidth}
			$\text{P-ERRK}_{4;\{5,6, 8, 13\}}$ 			    & $1.0$  & $1.0$  \\
			$\text{P-ERRK}_{4;13}$ 											& $1.81$ & $1.76$ \\
			\hdashline
			$\text{R-RK}_{4;4}$ 												& $2.69$ & $2.53$ \\
			$\text{R-TS}_{4;6}$ \cite{tselios2005runge} & $3.07$ & $2.73$ \\
			$\text{R-CKL}_{4;5}$ \cite{kennedy2000low}  & $2.01$ & $1.76$
		\end{tabular}
		\caption[Runtimes and number of scalar RHS evaluations of different fourth-order integrators compared to the optimized $p = 4, E = \{5, \allowbreak 6, \allowbreak 8, \allowbreak 13\} $ integrator for the \ac{vrMHD} flow around a cylinder.]
		{Runtimes and number of scalar \ac{RHS} evaluations \cref{eq:RHSEvals} of different entropy conservation, relaxed fourth-order integrators compared to the optimized $p = 4, E = \{5, \allowbreak 6, \allowbreak 8, \allowbreak 13\} $ integrator for the \ac{vrMHD} flow around a cylinder.}
		\label{tab:Runtimes_VRMHD_Cylinder}
	\end{table}
	\subsection{Inviscid Transonic Flow over NACA0012 Airfoil with Adaptive Mesh Refinement}
	\label{subsec:NACA0012}
	As a first inviscid example we consider the transonic flow over the NACA0012 airfoil.
	We follow the C1.2 testcase of the third installment of the International Workshop on High-Order CFD Methods which took place in 2015 \cite{3HOCFD}.
	Furthermore, we demonstrate that the \ac{PERRK} schemes can also be used efficiently in conjunction with \ac{AMR}, following upon a previous study \cite{doehring2024multirate}.
	\subsubsection{Setup}
	For this testcase, the Mach number is set to $\text{Ma} = 0.8$, the ratio of specific heats is $\gamma = 1.4$, and the angle of attack is $\alpha = 1.25\degree$.
	A range of uniformly refined fifth-order quadrilateral large meshes (farfield boundary placed 1000 cords from the airfoil center) are provided online \cite{3HOCFD} in \texttt{gmsh} \cite{geuzaine2009gmsh} format.
	We employ the medium resolution mesh with $2240$ cells truncated to second-order and $k=3$ solution polynomials for our computations.
	The surface fluxes are discretized using the Rusanov/Local Lax-Friedrichs flux \cite{RUSANOV1962304} and we use the flux-differencing shock-capturing procedure \cite{hennemann2021provably} with the entropy-stable flux by Chandrashekar \cite{Chandrashekar_2013}.
	For time integration $p = 4$, $E = \{5, \allowbreak 6, \allowbreak 7, \allowbreak 8, \allowbreak 10, \allowbreak 12, \allowbreak 14\}$ \ac{PERRK} schemes are constructed.
	The airfoil surface is modeled as a slip wall with no penetration \cite{Vegt2002SlipFB} involving the analytic solution of the pressure Riemann problem \cite{toro2013riemann}.
	At the outer boundaries the freestream flow is weakly enforced \cite{mengaldo2014guide}.

	We split the overall simulation into two steps:
	First, we use the provided mesh to advance the simulation until the shock position is approximately constant.
	Here, we run the simulation until $t_f = 100$ which corresponds to a convective time $t_c = t \, \frac{U_\infty}{c} \allowbreak = t \cdot \text{Ma} = 80$.
	As in the previous example, a ramp up of the $\text{CFL}$ number \cref{eq:CFLRampUp} is used to avoid immediate breakdown after initialization the flow field with the freestream state and unnecessary small timesteps at later times.
	For all considered time integration schemes the ramp-up time is set to $t_\text{CFL} = 1.5$.
	This is then followed by a restarted run where we allow for cells to being refined up to three times, i.e., $8$ times smaller than the base cell size which allows to resolve the shock much sharper, cf. \cref{fig:NACA0012StaticAMR}.
	As indicator we reuse the shock sensor that also governs the stabilization \cite{hennemann2021provably} which employs the product of density and pressure $\rho \cdot p$ as indicator variable.
	The simulation with adaptive mesh refinement is run for $80$ additional convective time units, i.e., until $t_f = 200$.

	For both simulation runs we use the Newton method to solve the relaxation equation.
	The iterations are forced to stop after $5$ steps, relaxation equation \cref{eq:RelaxationEquation} residual less or equal to $10^{-12}$, or stepsize $\Delta \gamma$ less or equal to $10^{-13}$.
	Interestingly, the entropy defect due to the considered Runge-Kutta methods is for this setup negligible, i.e., the unrelaxed ($\gamma = 1$) Runge-Kutta scheme satisfies the relaxation equation to desired accuracy.
	This is in contrast to the previous viscous examples where usually some Newton-Raphson updates to the relaxation parameter were required.
	\subsubsection{Performance Comparison}
	Comparisons of recorded runtime and computational effort are provided in \cref{tab:Runtimes_Inviscid_NACA0012}.
	For the first run on the static mesh we observe that the optimized standalone \ac{PERRK} scheme is performing remarkably well, which is due to a comparatively high stable \ac{CFL} number.
	For the second run on the adaptively refined mesh, however, $\text{R-CKL}_{4;5}$ \cite{kennedy2000low} is the most efficient standalone scheme.
	Thus, linear stability is no longer the timestep-limiting factor in this case, but nonlinear stability, i.e., positivity preservation.
	We would like to remark that the re-initialization of the additional partitioning \ac{PERRK} datastructures is of negligible cost, as it contributes only about \qty{0.1}{\percent} to the overall runtime.
	The entire \ac{AMR} procedure is also very efficient, as it amounts to merely \qty{3.3}{\percent} of the recorded runtime.
	\begin{table}
		\def\arraystretch{1.2}
		\centering
		\begin{tabular}{l?{2}c|c?{2}c|c}
			& \multicolumn{2}{c?{2}}{$1^\text{st}$ run (static)} & \multicolumn{2}{c}{$2^\text{nd}$ run (AMR)} \\
			Method & $\sfrac{\tau}{\tau^\star}$ & $\sfrac{N_\text{RHS}}{N_\text{RHS}^\star}$ & $\sfrac{\tau}{\tau^\star}$ & $\sfrac{N_\text{RHS}}{N_\text{RHS}^\star}$ \\
			\Xhline{5\arrayrulewidth}
			$\text{P-ERRK}_{4;\{5, \dots, 14\}}$ 			  & $1.0$  & $1.0$  & $1.0$  & $1.0$ \\
			$\text{P-ERRK}_{4;14}$ 											& $1.49$ & $1.55$ & $2.27$ & $2.46$ \\
			\hdashline
			$\text{R-RK}_{4;4}$ 												& $3.59$ & $3.22$ & $2.54$ & $2.18$ \\
			$\text{R-TS}_{4;6}$ \cite{tselios2005runge} & $4.08$ & $3.55$ & $2.95$ & $2.68$ \\
			$\text{R-CKL}_{4;5}$ \cite{kennedy2000low}  & $2.54$ & $2.33$ & $1.64$ & $1.88$
		\end{tabular}
		\caption[Runtimes and number of scalar RHS evaluations of different fourth-order integrators compared to the optimized $p = 4, E = \{5, \allowbreak 6, \allowbreak 7, \allowbreak 8, \allowbreak 10, \allowbreak 12, \allowbreak 14\} $ integrator for the transonic inviscid flow around the NACA0012 airfoil.]
		{Runtimes and number of scalar \ac{RHS} evaluations \cref{eq:RHSEvals} of different entropy conservation, relaxed fourth-order integrators compared to the optimized $p = 4, E = \{5, \allowbreak 6, \allowbreak 7, \allowbreak 8, \allowbreak 10, \allowbreak 12, \allowbreak 14\} $  integrator for the transonic inviscid flow around the NACA0012 airfoil.}
		\label{tab:Runtimes_Inviscid_NACA0012}
	\end{table}

\begin{figure}
		\centering
		\resizebox{.55\textwidth}{!}{\includegraphics{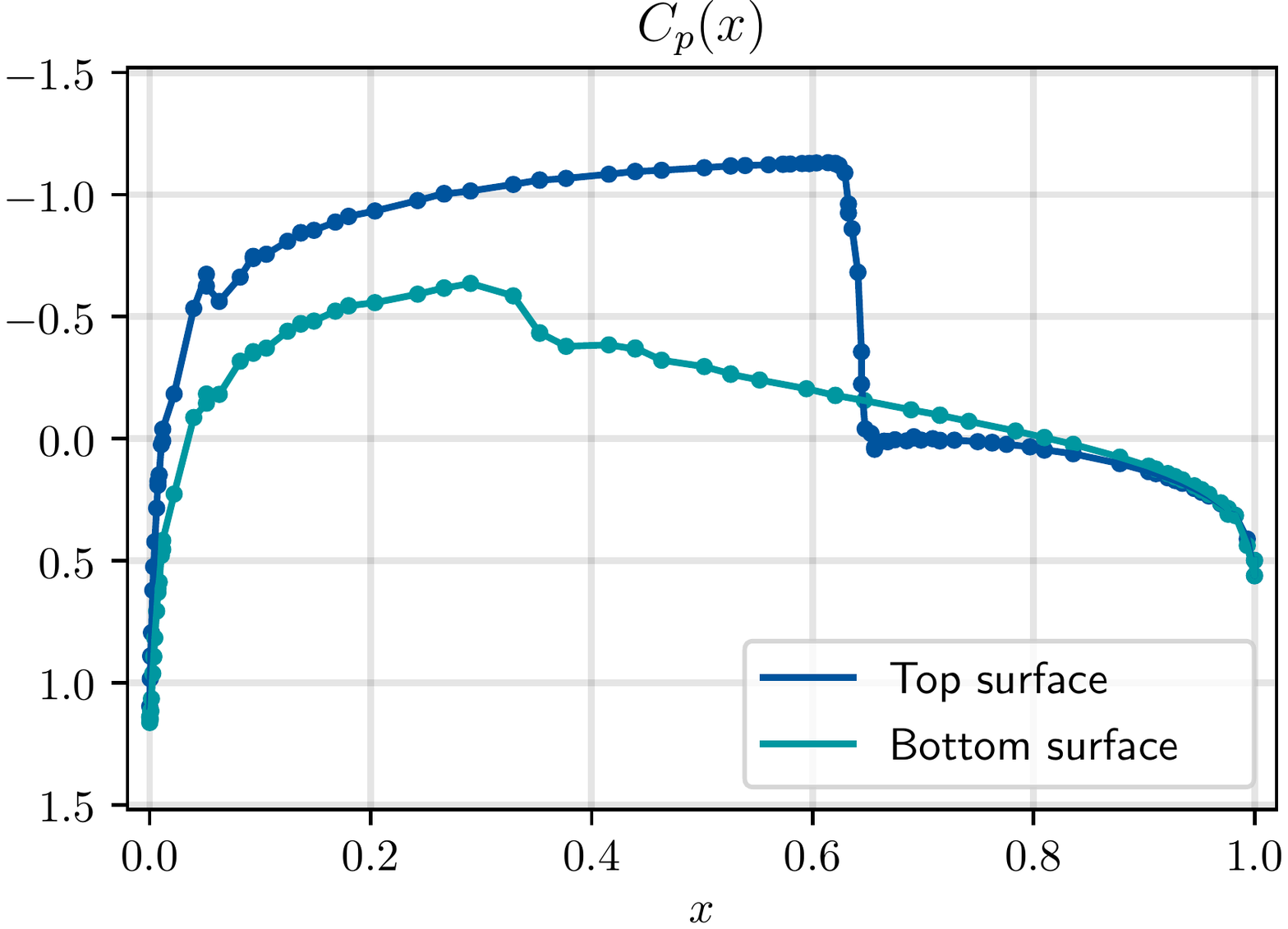}}
		\caption[Computational mesh for the VRMHD flow past a cylinder.]
		{Pressure coefficient \cref{eq:PressureCoefficient} recorded at top (blue) and bottom (petrol) NACA0012 airfoil surface.
		The dots correspond to the values recorded at the surface nodes.}
		\label{fig:PressureCoefficientNACA0012}
	\end{figure}

	We compute pressure-based lift \cref{eq:LiftCoeff} and drag \cref{eq:DragCoeffPressure} coefficients $C_{L,p}, \allowbreak C_{D,p}$ at final time $t_f = 200$, i.e., on the adapted mesh.
	These are compared to the data of Giangaspero et al. \cite{CaseC12_Twente} who employ fifth-order \ac{SBP} operators on a shock-fitted mesh with $577 \times 513$ vertices.
	Thus, this mesh has about $2.96 \cdot 10^5$ \ac{DoF} per solution field which is significantly more than the $\sim 4.2 \cdot 10^4$ \ac{DoF} of the adaptively refined mesh used in our simulation.
	Reference values for the lift and drag coefficients given in \cite{CaseC12_Twente} are $C_{L, \text{ref}} = 3.51607 \cdot 10^{-1}$ and $C_{D, \text{ref}} = 2.26216 \cdot 10^{-2}$.
	Based on these, we compute relative errors of the lift and drag coefficients
	\begin{equation}
		\label{eq:RelativeErrorsLiftDrag}
		e_{C_{\{L; D\}}} \coloneqq \frac{\left \vert C_{\{L; D\}} - C_{\{L; D\}, \text{ref} } \right \vert }{C_{\{L; D\}, \text{ref} } }
	\end{equation}
	at final time $t_f = 200$ and averaged over the last $20$ time units.
	For all time integrators, the final lift and drag coefficients match the given reference data quite well, with relative errors for $e_{C_L} \sim 0.15 \% - 0.26 \%$ and $e_{C_D} \sim 1.4 \% - 2.3 \%$.
	For the averaged coefficients $\bar C_{\{L;D\}}$ we obtained similar values for the relative errors of the mean values.
	Additionally, the pressure coefficient
	\begin{align}
		\label{eq:PressureCoefficient}
		C_p(x) \coloneqq \frac{p(x) - p_{\infty}}
		{0.5 \, \rho_{\infty} U_{\infty}^2 L_{\infty}}
	\end{align}
	obtained from the multirate \ac{PERRK} simulation is provided in \cref{fig:PressureCoefficientNACA0012}. This is again in good agreement with the submitted data of the participants of the workshop \cite{3HOCFD, CaseC12_Twente}, except for a small defect at $x \approx 0.08$ on the top surface where a higher resolution seems to be required.
	\subsection{Inviscid Transonic Flow over ONERA M6 Wing}
	\label{subsec:ONERA_M6}
\begin{figure}
		\centering
		\subfloat[{Rescaled ONERA M6 wing geometry following \cite{onera_m6_original_report, slater_study_1}.
		}]{
			\label{fig:ONERA_M6_geometry}
			\centering
			\includegraphics[width=0.47\textwidth]{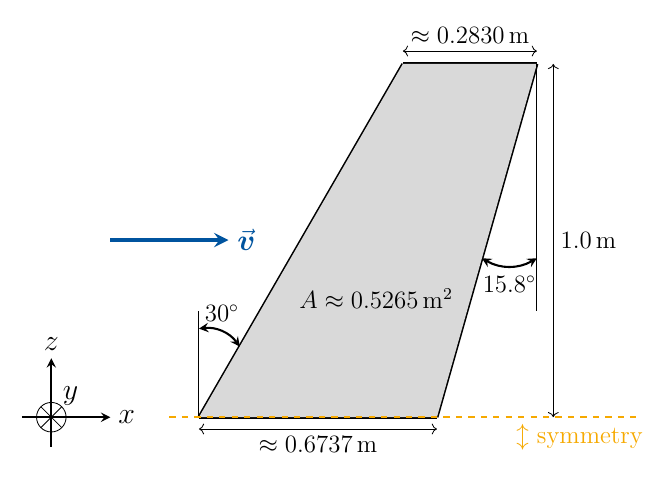}
		}
		\hfill
		\subfloat[{Computational mesh (2D/surface element edges).
		}]{
			\label{fig:ONERA_M6_Mesh}
			\centering
			\includegraphics[width=0.45\textwidth]{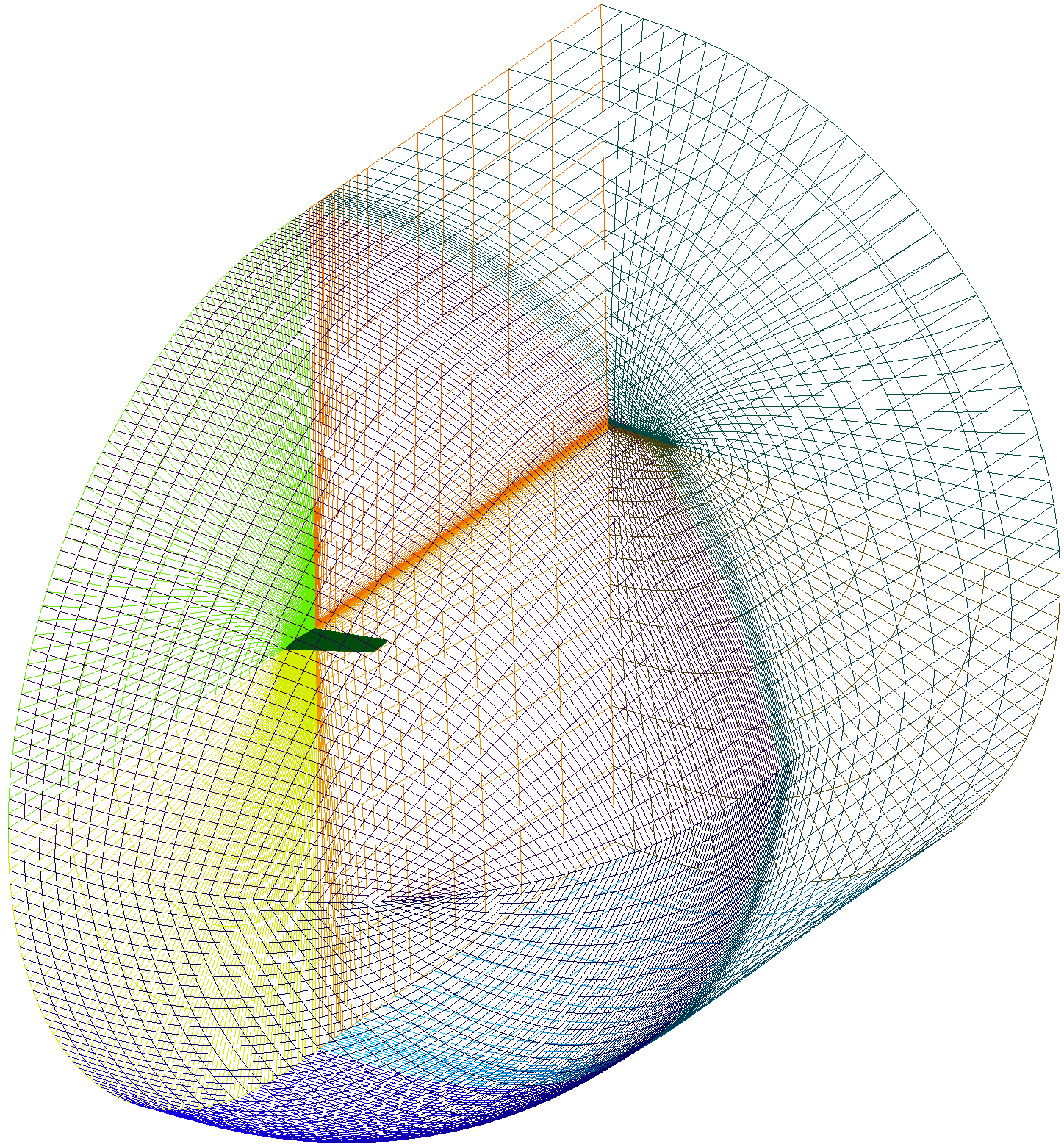}
		}
		\caption[ONERA M6 wing geometry and used computational mesh.]
		{ONERA M6 wing geometry \cref{fig:ONERA_M6_geometry} and used computational mesh \cref{fig:ONERA_M6_Mesh}.
		The displayed mesh is the sanitized version of the mesh from \cite{slater_study_1}.
		Note that the $y$-axis points into the plane in \cref{fig:ONERA_M6_geometry} and correspondingly upwards in \cref{fig:ONERA_M6_Mesh}.
		}
		\label{fig:ONERA_M6_geometry_mesh}
	\end{figure}

	We consider the transonic flow over the ONERA M6 wing, a well-known benchmark problem for compressible flow solvers.
	In our study we follow the geometry as presented in \cite{onera_m6_original_report}, although we employ a rescaled variant with wingspan set to $1$.
	A sketch of the wing geometry in the $x-z$ plane is given by \cref{fig:ONERA_M6_geometry}.
	Details on the chamber of the wing are also provided in the cited reference.
	%
	\subsubsection{Setup}
	The classic configuration is given by $\text{Ma} = 0.84$, angle of attack of $\alpha = 3.06\degree$ and Reynolds numbers $\text{Re} \sim \mathcal{O}\left(12 \cdot 10^6 \right)$.
	Here, we consider the purely hyperbolic scenario ($\text{Re} = \infty$) to avoid any need for turbulence modeling.

	For the geometry proposed in \cite{onera_m6_original_report} a purely hexahedral mesh is available for download from publicly available NASA webpages as part of the NPARC Alliance Validation Archive \cite{slater_study_1}.
	The mesh as published is split into four different blocks (cf. edge coloring in \cref{fig:ONERA_M6_Mesh}) where a united version in \texttt{gmsh} format is available from the HiSA repository \cite{heyns2014modelling}.
	Subsequently, we deleted duplicated nodes originating from different blocks and replaced node labels accordingly.
	Additionally, some mesh sanitizing was necessary to make this grid compatible with the \ac{DG} methodology as implemented in Trixi.jl.
	In particular, $72$ hexahedra are degenerated into prisms, i.e., neighboring hexahedra that share six instead of four nodes and edges.
	Those were easily sanitized by merging the two neighboring effective prisms into a single hexahedron.
	Unfortunately, two additionally elements on the wing surface were corrupted in a way that they share $5$ nodes/edges, which prevents such an easy fix.
	To resolve this, we deleted the corrupted elements and remeshed locally using \texttt{gmsh} \cite{geuzaine2009gmsh} which finally enabled us to use the mesh which consists now of $294838$ straight-sided elements.
	The final sanitized mesh is publicly available as part of the reproducibility repository which we provide \cite{doehring2025PERRK_ReproRepo}.
	The mesh extends in $x$-direction from $-6.373$ to $ 7.411$, in $y$-direction from $-6.375$ to $6.375$ and in $z$-direction from $0$ to $7.411$.
	The 2D element edges are depicted for illustration in \cref{fig:ONERA_M6_Mesh}.

	At $z=0$ plane symmetry is enforced by rotating the symmetry-plane crossing velocity $v_z$ to zero and copying density and pressure from the domain.
	At all other outer boundaries we weakly enforce \cite{mengaldo2014guide} the freestream flow.
	The wing surface is modeled as a slip-wall \cite{Vegt2002SlipFB} with an analytic solution of the pressure Riemann problem \cite{toro2013riemann}.

	We represent the solution using $k=2$ degree element polynomials which results in about $40$ million total \ac{DoF}, i.e., roughly $8$ million unknowns per solution field.
	The surface fluxes are computed with the Rusanov/Local Lax-Friedrichs flux \cite{RUSANOV1962304}.
	To stabilize the simulation it suffices in this case to use the flux-differencing approach \cite{gassner2016split, FISHER2013518, gassner2013skew, CHEN2017427, chan2018discretely}.
	The arising volume fluxes are discretized using the entropy-stable and kinetic-energy preserving flux by Ranocha \cite{Ranocha2020Entropy}.

	The simulation is initialized with the freestream values and we run the simulation until $t_f = 6.05$ for which the timings are obtained from the last $0.01$ time units.

\begin{figure}
		\centering
		\includegraphics[width=0.55\textwidth]{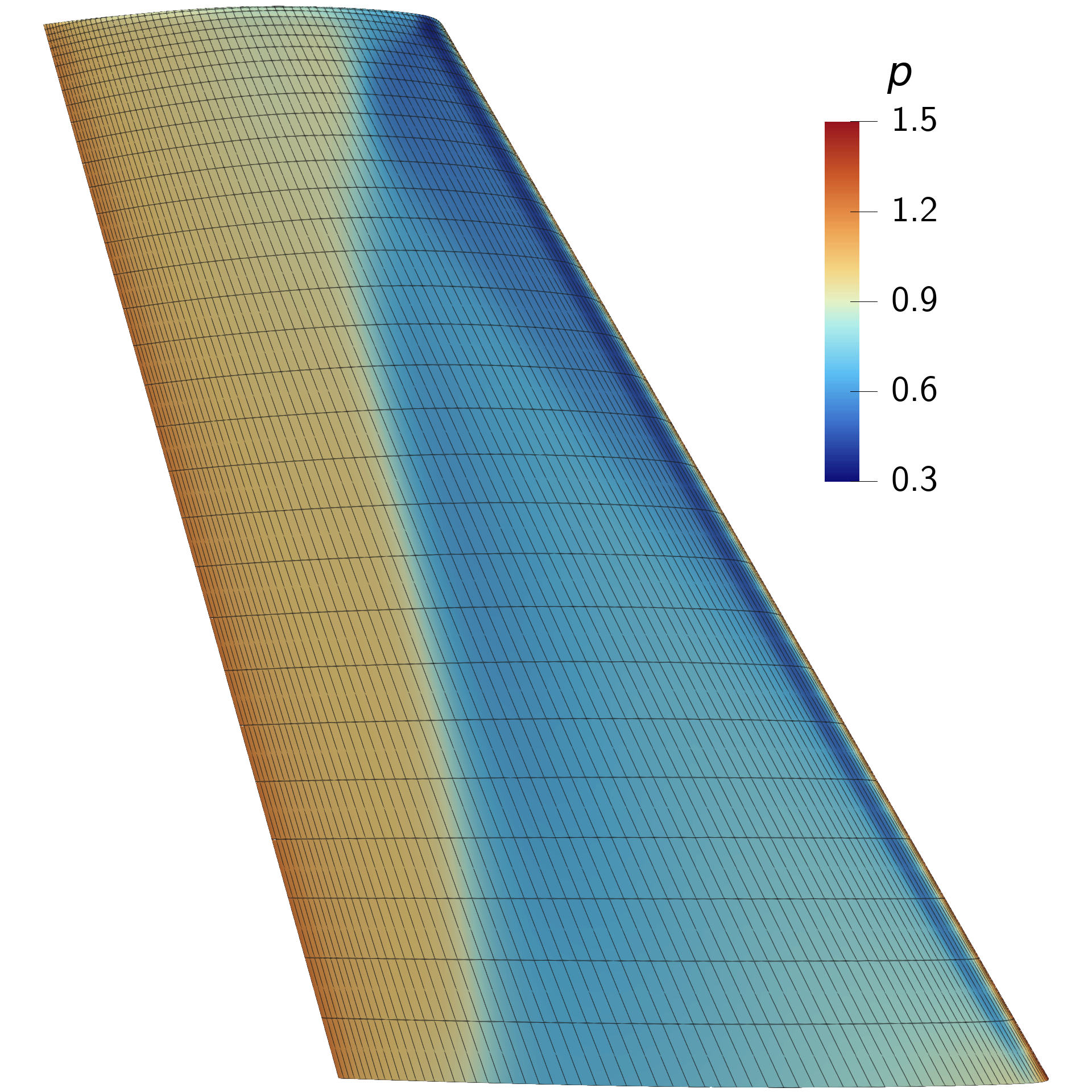}
		\caption[Initial density for the isentropic vortex testcase.]
		{Pressure $p$ on the ONERA M6 wing's upper surface at $t_f = 6.05$.
		The lambda shock pattern refers to the shape formed by the shock at the leading edge, where the pressure drops signficantly below the ambient value of $p=1$, and the standing shock across the wing across which the pressure increases again.
		Note that in constrast to \cref{fig:ONERA_M6_geometry} the $y$-axis points in this figure out of the view plane.
		}
		\label{fig:ONERA_M6_Pressure_UpperWing}
	\end{figure}

	We compute the lift coefficient \cref{eq:LiftCoeff} for $\alpha = 3.06\degree$ based on the reference area $A_\infty \approx 0.5265$.
	We obtain a lift coefficient of $C_L \approx 0.2951$ which is about \qty{3}{\percent} larger of what is reported in one of the tutorials of the \texttt{SU2} code \cite{economon2016su2}.
	This might be due to the fact that in the latter case the actual area of the wing is used 
	which is necessarily larger than the projected referenced area $A_\infty$ (which we employ) due to the curvature of the wing.
	In that tutorial a lift coefficient of $C_L \approx 0.2865$ is reported, obtained from a second-order in space, first-order in time simulation on a $\sim 5.8 \cdot 10^5$ element tetrahedral mesh.
	The classic lambda shock profile on the upper wing is displayed in \cref{fig:ONERA_M6_Pressure_UpperWing}.
	%
	\subsubsection{Performance Comparison}
	For multirate time integration third-order \ac{PERRK} schemes with stage evaluations 
	$ E = \{3, \allowbreak 4, \dots \allowbreak 14, \allowbreak 15 \}$
	are optimized for the discretized freestream state.
	The Runge-Kutta tableaus are constructed using the variant described in \cite{doehring2024multirate} which avoids negative entries, i.e., stage-downwinding in contrast to the original formulation \cite{nasab2022third}, see also \cref{eq:PERK_ButcherTableauOwn_p3}.

	As for the previous cases, we compare the \ac{PERRK} schemes to low-storage, standalone relaxed entropy-stable schemes of matching order in Butcher form.
	A natural choice for the "subdiagonal" methods cf. \cref{eq:ButcherTableauSubDiagonal} is the three-stage scheme $\text{RK}_{3;3}$ by Ralston \cite{ralston1962runge} since it has the minimum error bound among the $p=3, S = 3$ methods \cite{ralston1962runge}.
	In \cite{kennedy2000low} a four-stage, third-order scheme with $2N$-storage van-der-Houwen \cref{eq:ButcherTableau_vdH2N} implementation is presented which we employ here for comparison.
	The multirate scheme is significantly more efficient than the considered standalone schemes, requiring only about a quarter of the number of \ac{RHS} evaluations.
	This improvement compared to the fourth-order schemes is mainly due to the fact that here the minimum number of stage evaulations is $\min_r E^{(r)} = 3$, compared to $5$ for the fourth-order schemes.
	In terms of runtime, however, we observe that the savings are unfortunately not as pronounced.
	This is mostly due to non-ideal parallel scaling of the partitioned Runge-Kutta solver, which is due to distribution of too few data across too many computing units (threads) on the finest levels.
	A possible remedy to this would be to use a constant ratio of \acp{DoF} per computing unit.
	Unfortunately, reducing the number of parallel computing units, i.e., threads in certain parts of the program is currently not (natively) supported by the Julia programming language.
	\begin{table}
		\def\arraystretch{1.2}
		\centering
		\begin{tabular}{l?{2}c|c}
			Method & $\sfrac{\tau}{\tau^\star}$ & $\sfrac{N_\text{RHS}}{N_\text{RHS}^\star}$ \\
			\Xhline{5\arrayrulewidth}
			$\text{P-ERRK}_{3;\{3, \dots, 15\}}$		    & $1.0$  & $1.0$  \\
			$\text{P-ERRK}_{3;15}$ 											& $2.98$ & $3.92$ \\
			\hdashline
			$\text{R-RK}_{3;3}$ 												& $3.79$ & $4.26$ \\
			$\text{R-CKL}_{3;4}$ \cite{kennedy2000low}  & $3.54$ & $4.14$
		\end{tabular}
		\caption[Runtimes and number of scalar RHS evaluations of different third-order integrators compared to the optimized $p = 3, E = \{3, \allowbreak 4, \allowbreak \dots , \allowbreak 14, , \allowbreak 15\} $ integrator for the flow over the ONERA M6 wing.]
		{Runtimes and number of scalar \ac{RHS} evaluations \cref{eq:RHSEvals} of different entropy conservation, relaxed third-order integrators compared to the optimized $p = 3, E = \{3, \allowbreak 4, \allowbreak \dots , \allowbreak 14, , \allowbreak 15\} $ integrator for the flow over the ONERA M6 wing.}
		\label{tab:Runtimes_ONERAM6}
	\end{table}
	\subsection{Viscous Transonic Flow over NASA Common Research Model (CRM)}
	\label{subsec:NASA_CRM}
	As the final example we consider the viscous transonic flow over the NASA \ac{CRM} \cite{vassberg2008development} which was also considered in the 2015 International Workshop on High-Order CFD Methods \cite{3HOCFD}.
	Here, we also use the transonic Mach number $\text{Ma} = 0.85$ but increase the Reynolds number to $\text{Re} = 200 \cdot 10^6$ which is based on the wing reference chord of \qty{7.005}{\metre}.
	This configuration allows us to demonstrate increased robustness of the \ac{PERRK} schemes compared to the standard ones.
	As for the NACA0012 airfoil, a mesh consisting of $79505$ 
	hexahedra is provided online \cite{3HOCFD} which we use for our computations.
	This mesh is way too coarse to resolve the flow accurately, but yet allows us to demonstrate the benefit of the relaxation approach for the multirate \ac{PERK} schemes.
	\subsubsection{Setup}
	The solution fields are reprensented using $k=2$ polynomials, which results in about $2.1$ million \ac{DoF} per solution field.
	The surface fluxes are discretized using the \ac{HLL} two-wave solver \cite{harten1983upstream} with Davis-type wave speed estimates \cite{doi:10.1137/0909030}.
	We stabilize the simulation using the subcell limiting approach \cite{hennemann2021provably} with the entropy-stable and kinetic energy preserving flux by Ranocha \cite{Ranocha2020Entropy} for the volume fluxes.
	The viscous gradients are discretized with the \ac{BR1} scheme \cite{bassi1997high, gassner2018br1}.
	The ratio of specific heats is set to $\gamma = 1.4$, the Prandtl number is set to $0.72$, and the viscosity is computed from the Reynolds number and is constant across the domain.
	The aircraft is modeled as a no-slip adiabatic wall, the freestream values are weakly enforced at the farfield boundaries \cite{mengaldo2014guide}.
	Symmetry across the $z$-plane is enforced by rotating the velocity component $v_z$ to zero and eliminating heat transfer across the plane (adiabatic).
	No turbulence modeling is employed.
	For this semidiscretization a third-order \ac{PERRK} scheme with $E = \{ 3, \allowbreak 4, \allowbreak 5,\allowbreak 7, \allowbreak 8, \allowbreak 9, \allowbreak 10, \allowbreak 11, \allowbreak 12, \allowbreak 15\}$ is constructed and we compare against the $S=15$ optimized standalone, $\text{CKL}_{3;4}$ \cite{kennedy2000low} and $\text{RK}_{3;3}$ \cite{ralston1962runge} schemes.
	The relaxation equation is solved using Newton iterations, limited to five steps, objective residual less or equal to $10^{-13}$, and stepsize $\Delta \gamma$ less or equal to $10^{-13}$.
	\subsubsection{Performance Comparison}
\begin{figure}
		\centering
		\includegraphics[width=0.95\textwidth]{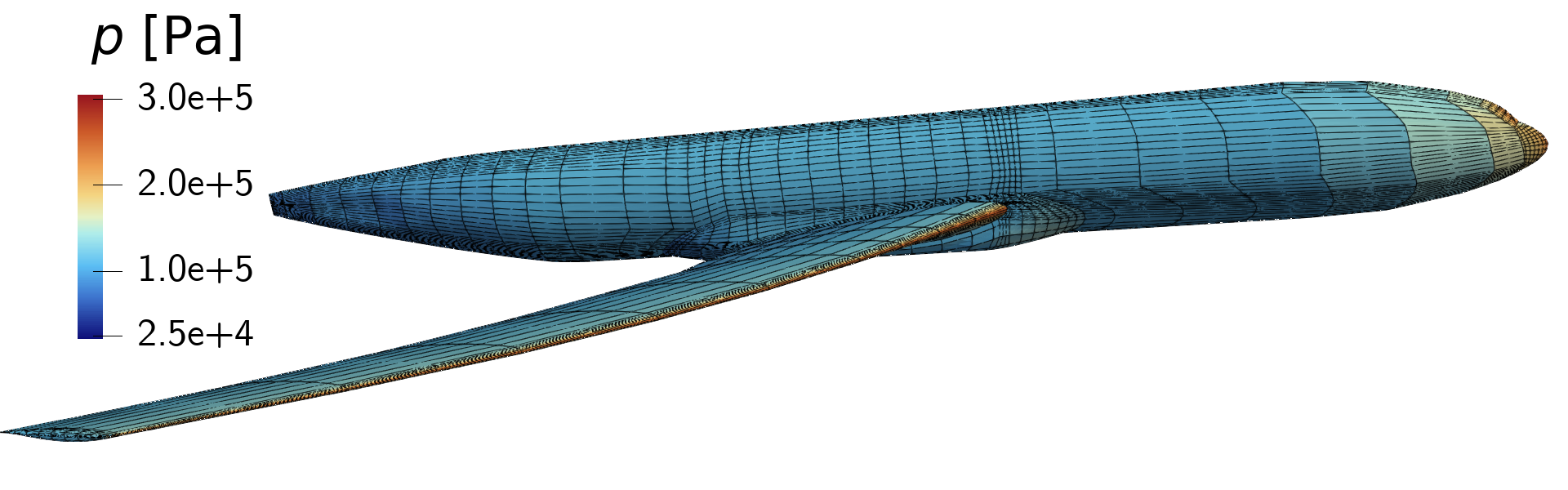}
		\caption[Pressure and surface mesh on the NASA Common Reserach Model aircraft at $t_f = 1.5 \cdot 10^{-5} \text{ s}$.]
		{Pressure $p \text{ [Pa]}$ and surface mesh for the NASA \ac{CRM} at $t_f = 1.5 \cdot 10^{-5} \text{ s}$.
		}
		\label{fig:CRM_Pressure}
	\end{figure}

	To begin, we highlight that the relaxed scheme maintains positivity 
	\qty{22}{\percent} longer than the standard scheme ($\sim 1.97 \cdot 10^{-5} \text{ s}$ vs. $\sim 1.61 \cdot 10^{-5}\text{ s}$).
	Restarting the simulation from $1.5 \cdot 10^{-5} \text{ s}$ with a smaller \ac{CFL} number shows the same behaviour, where the relaxed scheme runs again about \qty{20}{\percent} longer (after restart) than the plain \ac{PERK} method before breaking down ($3.26 \cdot 10^{-5} \text{ s}$ vs. $2.97 \cdot 10^{-5} \text{ s}$).
	Thus, we observe the same stabilizing effect of the relaxation methodology as encountered for the Kelvin-Helmholtz instability, recall \cref{subsec:KelvinHelmholtzInstability}.
	For illustration, we provide a snapshot of the pressure field and the surface mesh at $t_f = 1.5 \cdot 10^{-5} \text{ s}$ in \cref{fig:CRM_Pressure}.

	In terms of improved performance compared to the standalone schemes, we observe significant svaings in runtime and number of \ac{RHS} evaluations for the multirate \ac{PERRK} scheme.
	In particular, \ac{RHS} evaluations are reduced by a factor of about three, while the runtime is decreased even further, with speedups ranging from $3.51$ to $3.87$.
	\begin{table}
		\def\arraystretch{1.2}
		\centering
		\begin{tabular}{l?{2}c|c}
			Method & $\sfrac{\tau}{\tau^\star}$ & $\sfrac{N_\text{RHS}}{N_\text{RHS}^\star}$ \\
			\Xhline{5\arrayrulewidth}
			$\text{P-ERRK}_{3;\{3, \dots, 15\}}$		    & $1.0$  & $1.0$  \\
			$\text{P-ERRK}_{3;15}$ 											& $3.67$ & $3.02$ \\
			\hdashline
			$\text{R-RK}_{3;3}$ 												& $3.51$ & $2.85$ \\
			$\text{R-CKL}_{3;4}$ \cite{kennedy2000low}  & $3.87$ & $3.04$
		\end{tabular}
		\caption[Runtimes and number of scalar RHS evaluations of different third-order integrators compared to the optimized $p = 3, E = \{ 3, \allowbreak 4, \allowbreak 5,\allowbreak 7, \allowbreak 8, \allowbreak 9, \allowbreak 10, \allowbreak 11, \allowbreak 12, \allowbreak 15\} $ integrator for the flow over the NASA CRM aircraft.]
		{Runtimes and number of scalar \ac{RHS} evaluations \cref{eq:RHSEvals} of different entropy conservation, relaxed third-order integrators compared to the optimized $p = 3, E = \{ 3, \allowbreak 4, \allowbreak 5,\allowbreak 7, \allowbreak 8, \allowbreak 9, \allowbreak 10, \allowbreak 11, \allowbreak 12, \allowbreak 15\} $ integrator for the flow over the NASA CRM aircraft.}
		\label{tab:Runtimes_CRM}
	\end{table}
	\section{Conclusions}
	\label{sec:Conclusions}
	In this work we presented a novel class of entropy conservation/stable high-order, optimized, multirate time integration schemes suitable for the integration of convection-dominated compressible flow problems.
	These schemes are based on the \ac{PERK} methods and allow for the same non-intrusive implementation based on a coefficient-based partitioned semidiscretization.
	The \ac{PERRK} methods maintain the same order of accuracy as the corresponding \ac{PERK} methods, preserve linear invariants such as mass, momentum, and energy, and introduce no defects in the global entropy.
	We have shown that the \ac{PERK} schemes benefit significantly from the relaxation approach, as it enhances robustness of the multirate schemes.
	In particular, this is demonstrated for underresolved viscous flows, where the schemes with relaxation maintain positivity for longer time intervals or are stable for larger \ac{CFL} numbers.
	This is due to the relaxation parameter $\gamma$ acting as an adaptive timestep controller, which in turn has a limiting effect.

	We have shown that the \ac{PERRK} schemes outperform standard Relaxation-Runge-Kutta methods for a collection of problems, ranging from low Mach and low Reynolds number flow regimes to transonic inviscid and viscous flows.
	We highlight here the flow of a magnetized fluid around a cylinder governed by the \ac{vrMHD} equations, the transonic inviscid flow over a NACA0012 airfoil with \ac{AMR} and the simulation of the ONERA M6 wing.
	Depending on the problem, \ac{RHS} evaluations are reduced by factors up to four and speedups in excess of three are observed.

	In future work we will explore the application of the \ac{PERRK} schemes to flows involving explicit turbulence modeling.
	In addition, the construction of implicit-explicit (IMEX) \ac{PERK}/\ac{PERRK} schemes is of interest, for instance for simulations involving chemical reactions.
	\section*{Data Availability}
	All data generated or analyzed during this study are included in this published article and its supplementary information files.
	\section*{Code Availability \& Reproducibility}
	We provide a reproducibility repository publicly available on GitHub \cite{doehring2025PERRK_ReproRepo}.
	\section*{Acknowledgments}
	Funding by German Research Foundation (DFG) under Research Unit FOR5409: 
	\\ "Structure-Preserving Numerical Methods for Bulk- and Interface-Coupling of Heterogeneous Models ~(SNuBIC)~" (grant \#463312734). \\
	The authors are grateful for Andrew Winters' assistance in generation of the symmetric mesh for the visco-resistive MHD flow past a cylinder.
	%
	%
	\section*{Declaration of competing interest}
	The authors declare the following financial interests/personal relationships which may be considered as potential competing interests:
	Daniel Doehring's financial support was provided by German Research Foundation.

	\section*{CRediT authorship contribution statement}
	\noindent
	\textbf{Daniel Doehring}: Conceptualization, Methodology, Software, Validation, Investigation, Writing - original draft. \\
	\textbf{Hendrik Ranocha}: Methodology, Software, Writing - review \& editing. \\
	\textbf{Manuel Torrilhon}: Conceptualization, Funding acquisition, Supervision, Writing - review \& editing.
	
	\newpage 
	%
	\appendix
	\section{P-ERK Butcher Tableau Archetypes}
	\label{sec:P-EKR_ButcherTableauArchetypes}
	In \cref{subsec:PERKMs} we discussed the central properties of the \ac{PERK} schemes.
	Here, we present the structural archetypes for the \ac{PERK}/\ac{PERRK} schemes for orders of consistency $p = 2, 3, 4$.
	\subsection{Second-Order Scheme}
	The archetype for the second-order \ac{PERK} is given by \cref{eq:PERK_ButcherTableauClassic_p2_E36}.
	\subsection{Third-Order Scheme}
	The archetype of the third-order \ac{PERK} schemes as presented in \cite{nasab2022third} is given by
	\begin{equation}
		\label{eq:PERK_ButcherTableauClassic_p3}
		\renewcommand\arraystretch{1.3}
		\begin{array}
				{c|c|c c c c c c}
				i & \boldsymbol c & & & A^{(r)} & & &   \\
				\hline
				1 & 0 & & & & & &             					\\
				2 & c_2 & c_2 & & & &                   \\
				3 & c_3 & c_3 - a_{3,2} & a_{3,2} & & & \\ 
				\vdots & \vdots & \vdots & & & & &      \\
				S - 2 & c_{S-2} & c_{S-2} - a_{S-2, S-3} & \dots & a_{S-2, S-3} & & &         \\
				S-1 & \sfrac{1}{3} & \sfrac{1}{3} - a_{S-1,S-2} & 0 & 0 & a_{S-1,S-2} & &     \\
				S & 1 & 1 - a_{S,S-1} & 0 & 0 & 0 & a_{S,S-1} &                               \\
				\hline
				& \boldsymbol b^T & 0 & 0 & 0 & 0 & \sfrac{3}{4} & \sfrac{1}{4}
			\end{array}
	\end{equation}
	where the superscript $(\cdot)^{(r)}$ has been truncated in the coefficients $a_{i,j}$ of the $A$-matrix for the sake of compactness.
	While this scheme has the sparsest weight vector $\boldsymbol b^T$, it comes at the cost that the $E=S=3$ scheme has negative entry $a_{3, 1} = -1$ \cite{doehring2024multirate}.
	This corresponds to undesired downwinding \cite{gottlieb1998total} of the $\boldsymbol K_1$ stage in the construction of the $i=3$ intermediate approximation $\boldsymbol U_3$, cf. \cref{eq:PartitionedRKSecondEq}.
	Thus, we derived in \cite{doehring2024multirate} a new scheme which avoids this issue by basing the third-order scheme on the $E=S=3$ Shu-Osher scheme \cite{shu1988efficient}.
	This leads to the following Butcher tableau
	\begin{equation}
		\label{eq:PERK_ButcherTableauOwn_p3}
		\renewcommand\arraystretch{1.3}
		\begin{array}
				{c|c|c c c c c c}
				i & \boldsymbol c & & & A^{(r)} & & &   \\
				\hline
				1 & 0 & & & & & &             					\\
				2 & c_2 & c_2 & & & &                   \\
				3 & c_3 & c_3 - a_{3,2} & a_{3,2} & & & \\ 
				\vdots & \vdots & \vdots & & & & &      \\
				S - 2 & c_{S-2} & c_{S-2} - a_{S-2, S-3} & \dots & a_{S-2, S-3} & & &         \\
				S-1 & 1 & 1 - a_{S-1,S-2} & 0 & 0 & a_{S-1,S-2} & &     \\
				S & \sfrac{1}{2} & 1 - a_{S,S-1} & 0 & 0 & 0 & a_{S,S-1} &                               \\
				\hline
				& \boldsymbol b^T & \sfrac{1}{6} & 0 & 0 & 0 & \sfrac{1}{6} & \sfrac{2}{3}
			\end{array}
	\end{equation}
	We emphasize that the usage of the first Runge-Kutta stage $\boldsymbol K_1$ in the final update step does not increase storage requirements of the scheme as $\boldsymbol K_1$ is always kept in memory to compute the intermediate approximations $\boldsymbol U_i$.
	Also, the computational costs per stage stay at two scalar-vector multiplications and one vector-vector addition.
	
	We have not yet mentioned how to choose the free abscissae $c_i, i = 2, \dots, S-2$. 
	In this work, we set
	\begin{equation}
		\label{eq:Free_Timesteps}
		c_i = \frac{i-1}{S-3}, \quad i = 1, \dots, S-2
	\end{equation}
	which corresponds to a linear distribution of timesteps between $c_1=0$ and $c_{S-2} = 1$, similar to the second-order case \cite{vermeire2019paired}.
	We remark that more sophisticated choices of the abscissae $\boldsymbol c$ are possible, for instance to improve the internal stability properties of the scheme \cite{ketcheson2014internal}.
	\subsection{Fourth-Order Scheme}
	In \cite{doehring2024fourth} we derived a fourth-order \ac{PERK} scheme which is based on a $S=5$ method.
	In contrast to the second- and third-order schemes, additional coupling conditions in the order conditions arise which render the construction of the Butcher tableau more involved.
	We were able to find the one-parameter family 
	\begin{equation}
		\label{eq:PERK_ButcherTableau_p4}
		\renewcommand\arraystretch{1.3}
		\begin{array}{c|c|c c c c c c}
			i & \boldsymbol c & & & A^{(r)} & & & \\
			\hline
			1 & 0 & & & & & &             										 \\
			2 & c_2 & c_2 & & & & &        \\
			3 & c_3 & c_3 - a_{3,2} & a_{3,2} & & & &      \\ 
			\vdots & \vdots & \vdots & \ddots & &  & \\
			S - 2 & c_{S-2} & c_{S-2} - a_{S-2, S-3} & \dots & a_{S-2, S-3} & & &         \\
			S - 1 & 0.5 + \frac{\sqrt{3}}{6} & 0.5 + \frac{\sqrt{3}}{6} - a_{S-1,S-1}  & 0 & \dots & a_{S-1,S-2}  & & \\
			S & 0.5 - \frac{\sqrt{3}}{6} & 0.5 - \frac{\sqrt{3}}{6} - a_{S,S-1}  & 0 & \dots & 0 & a_{S,S-1} & \\
			\hline
			& & 0 & 0 & 0 & 0 & 0.5 & 0.5
		\end{array}
	\end{equation}
	with 
	\begin{subequations}
		\label{eq:4thOrderPERK_BaseCoeffs}
		\begin{align}
			c_{S-2} &= 0.479274057836310 \\
			c_{S-1} &= 0.5 + \frac{\sqrt{3}}{6} = 0.788675134594813 \\
			c_{S} &= 0.5 - \frac{\sqrt{3}}{6} = 0.211324865405187 \\
			\label{eq:4thOrderPERK_BaseCoeffs_aS2}
			a_{S-2, S-3} &= \frac{0.114851811257441}{c_{S-3}} \\
			a_{S-1, S-2} &= 0.648906880894214 \\
			a_{S, S-1} &= 0.0283121635129678
		\end{align}	
	\end{subequations}
	where $c_{S-3}$ is a free parameter, as the remaining abscissae $c_i, i = 2, \dots, S-3$.
	In \cite{doehring2024fourth} we demonstrated that setting all of these to unity results in the scheme with the best internal stability properties, given the restriction that $c_i \leq 1$.
	\newpage 
	\section{Supplementary Figures}
	\begin{figure}[!ht]
		\centering
		\subfloat[{Spectrum of the update matrix $D$ corresponding to the standard \ac{PERK} scheme.}]{
			\label{fig:LinearStabilityStandard}
			\centering
			\resizebox{.47\textwidth}{!}{\includegraphics{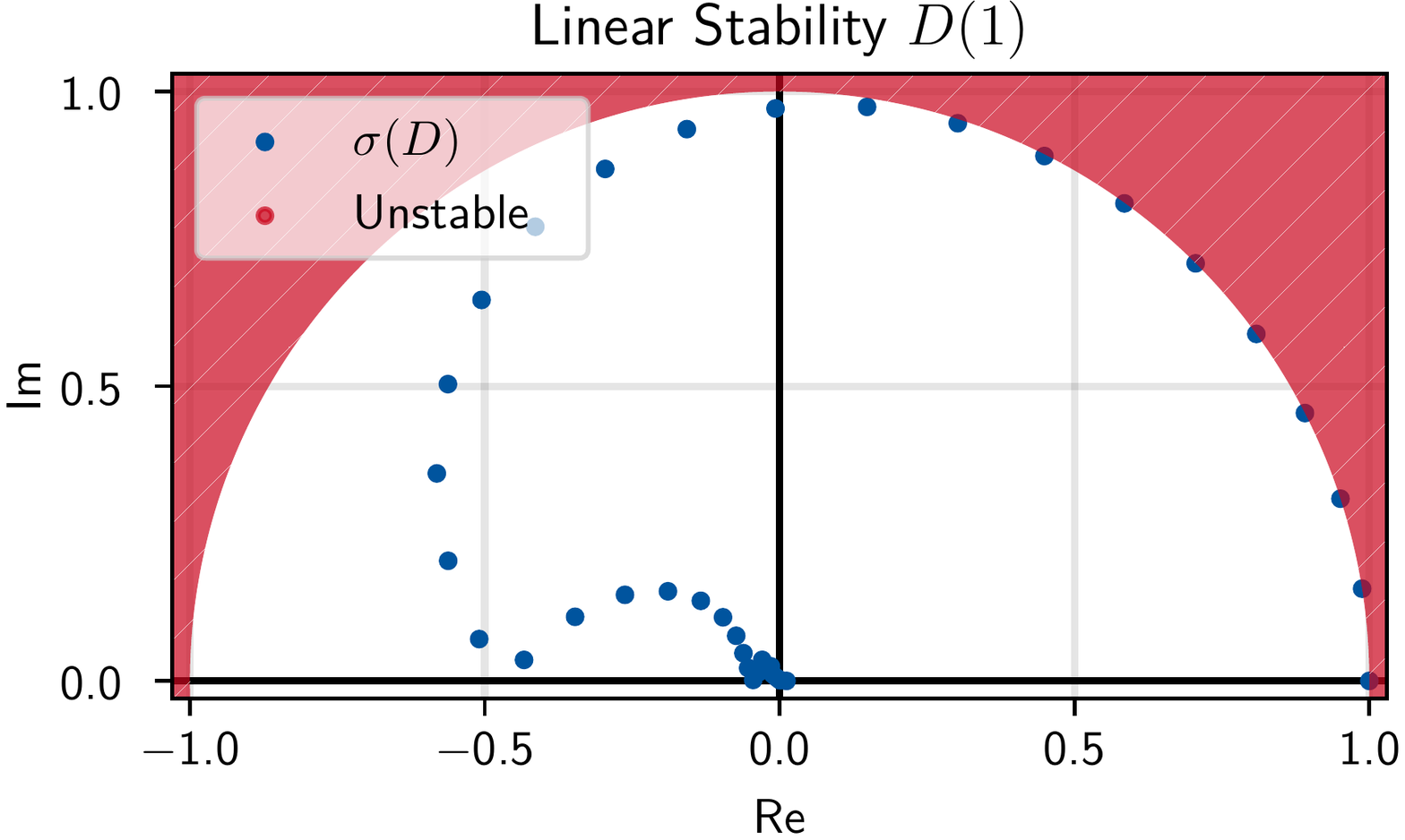}}
		}
		\hfill
		\subfloat[{Spectrum of the relaxed update matrix $D(\gamma_1)$ corresponding to the relaxed \ac{PERK} scheme.}]{
			\label{fig:LinearStabilityRelaxation}
			\centering
			\resizebox{.47\textwidth}{!}{\includegraphics{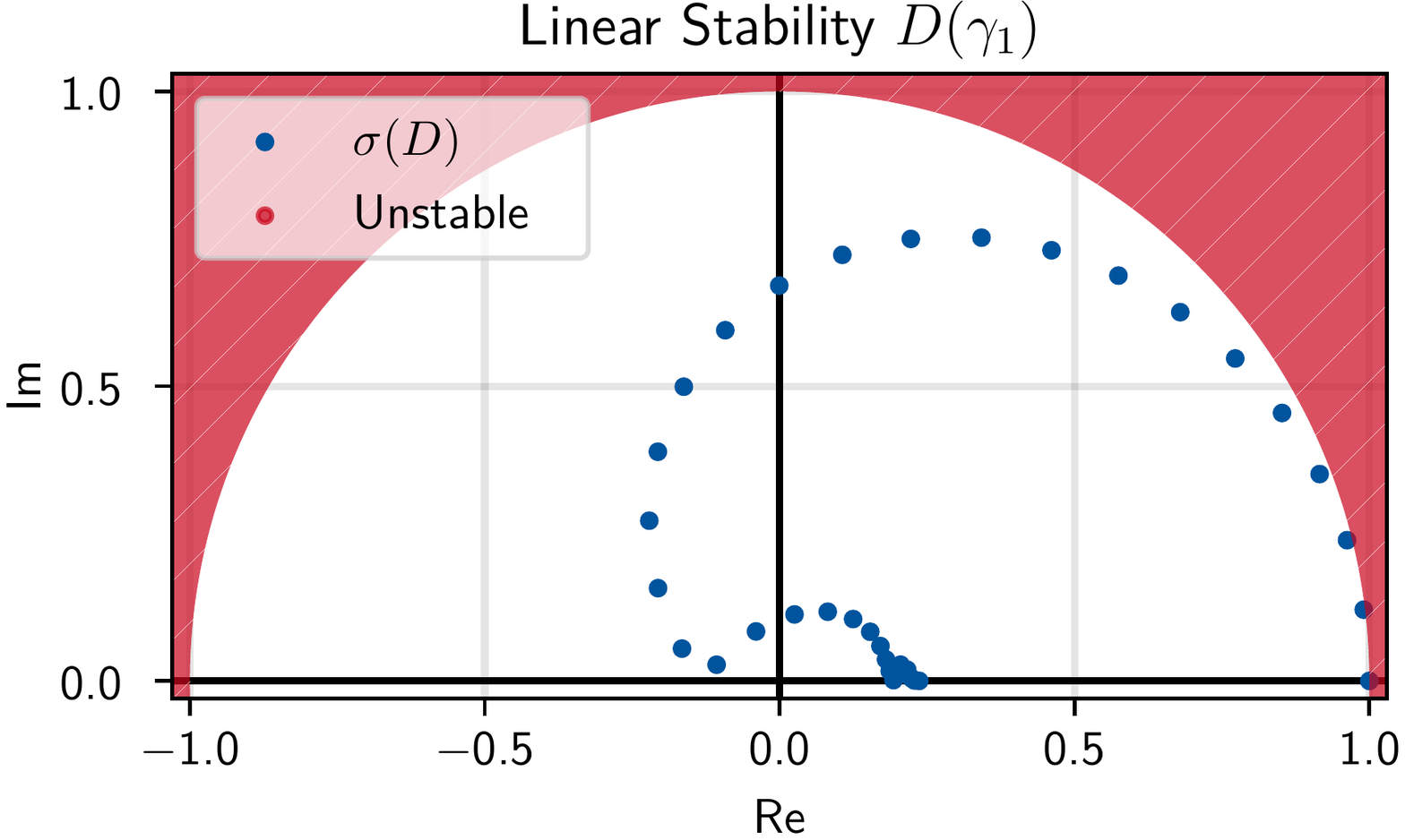}}
		}
		\caption[Spectra $\boldsymbol \sigma(D)$ of the fully discrete matrix $D(\gamma)$ of the $\text{P-ERK}_{2; \{8, 16\}}$ scheme applied to \eqref{eq:LinearTwoLevelPartitionedODESys}.]
		{\cref{subsec:EntropyRelaxationTimeLimiting_1DAdvection}: Spectra $\boldsymbol \sigma(D)$ of the fully discrete matrices $D(\gamma)$ of the $\text{P-ERK}_{2; \{8, 16\}}$ scheme applied to \eqref{eq:LinearTwoLevelPartitionedODESys}.
		The unstable region with $\vert \lambda \vert > 1$ is shaded in red.
		The spectrum of the update matrix $D$ corresponding to the standard \ac{PERK} scheme is shown in \cref{fig:LinearStabilityStandard}, while the spectrum of the relaxed update matrix $D(\gamma_1)$ corresponding to the relaxed \ac{PERK} scheme is shown in \cref{fig:LinearStabilityRelaxation}.
		}
		\label{fig:LinearStability}
	\end{figure}
	\begin{figure}[!ht]
		\centering
		\subfloat[{Density $\rho$ at final time $t_f = 3.2$ for the Kelvin-Helmholtz instability testcase.}]{
			\label{fig:Density_KelvinHelmholtzInstability}
			\centering
			\includegraphics[width=0.45\textwidth]{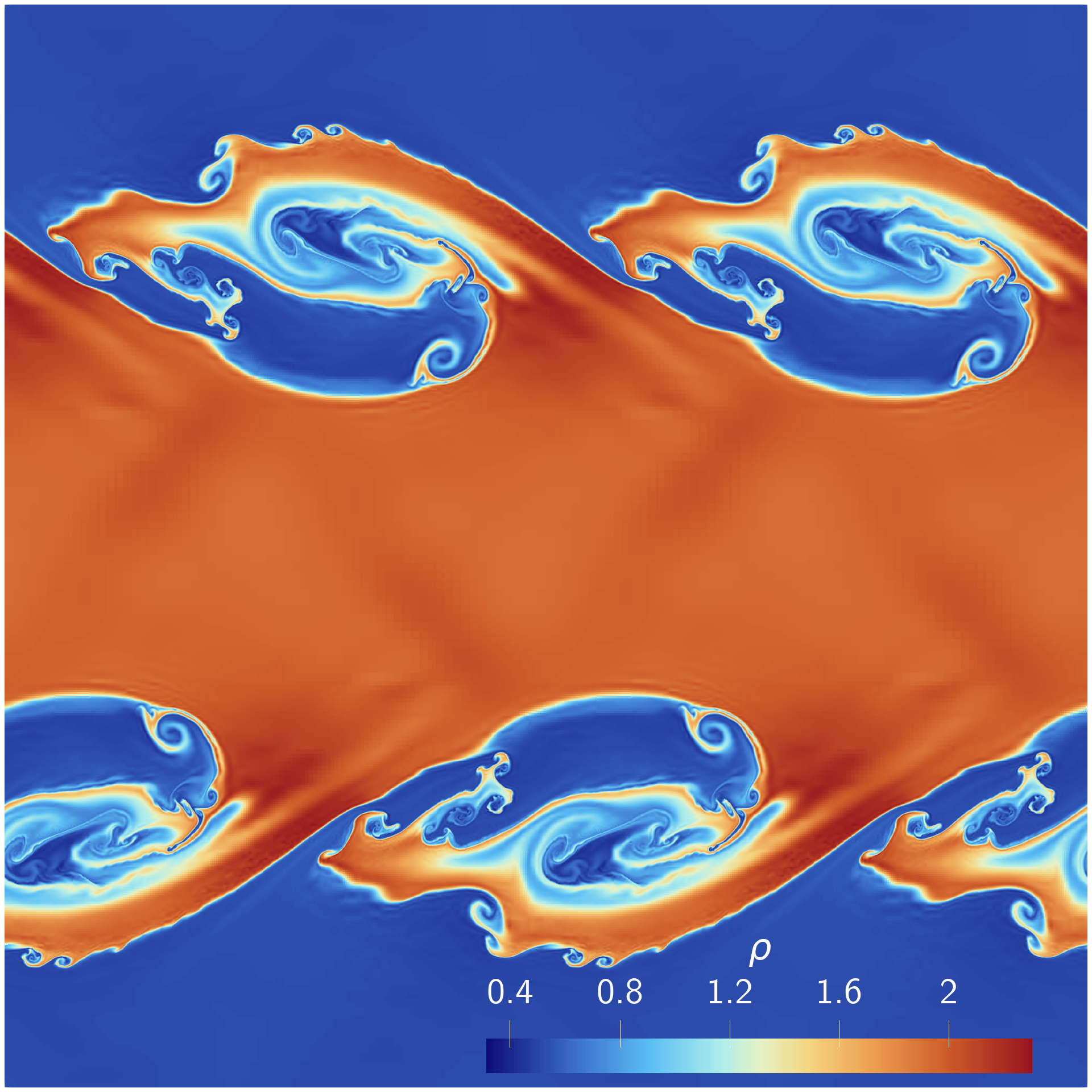}
		}
		\hfill
		\subfloat[{Adaptive mesh at final time $t_f = 3.2$ for the Kelvin-Helmholtz instability testcase.}]{
			\label{fig:Mesh_KelvinHelmholtzInstability}
			\centering
			\includegraphics[width=0.45\textwidth]{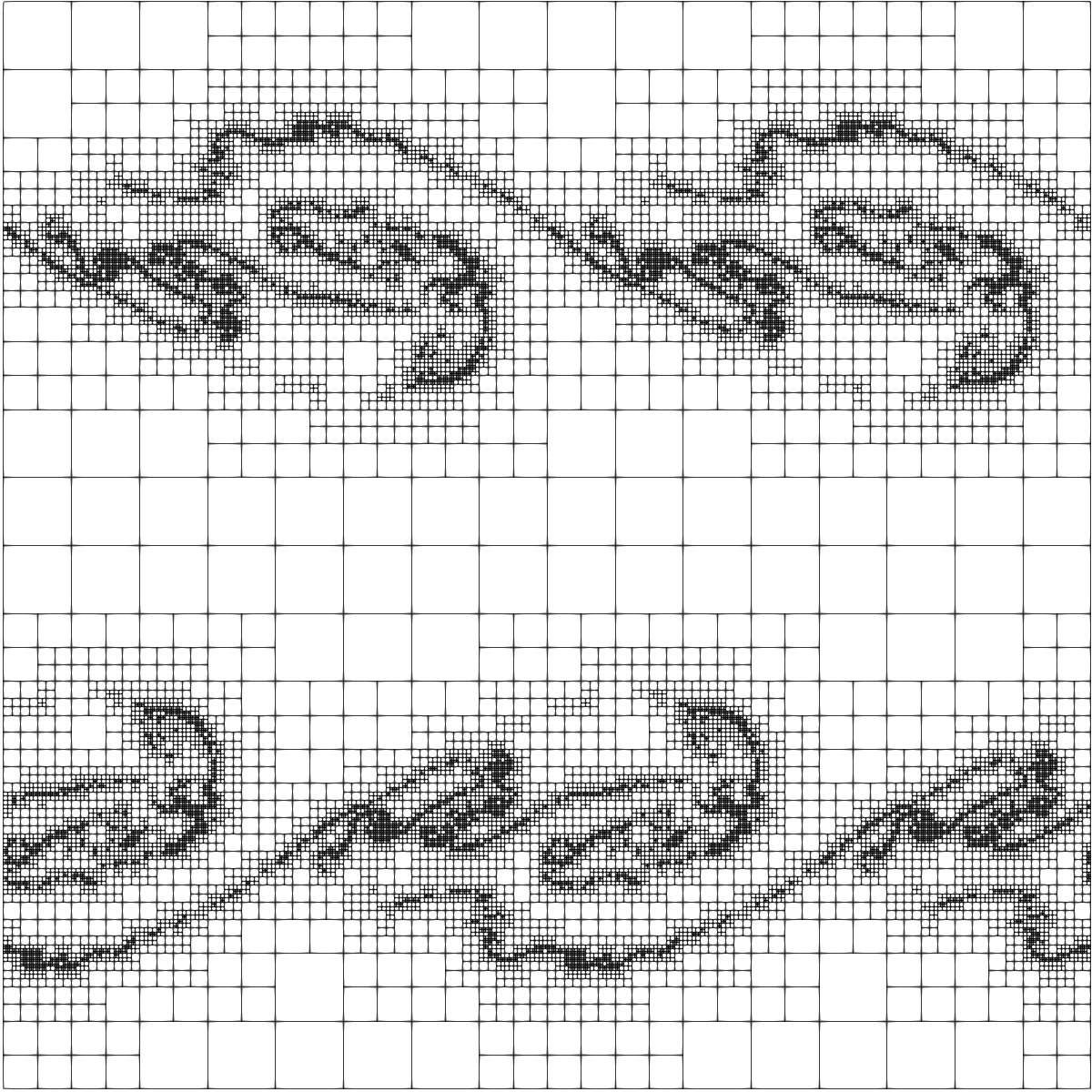}
		}
	\caption[Density $\rho$ and adaptive mesh (six levels) at final time $t_f = 3.2$ for the Kelvin-Helmholtz instability testcase obatined with the relaxed multirate P-ERRK scheme.]
	{\cref{subsec:KelvinHelmholtzInstability}: Density $\rho$ (\cref{fig:Density_KelvinHelmholtzInstability}) and adaptive mesh (\cref{fig:Mesh_KelvinHelmholtzInstability}) at final time $t_f = 3.2$ for the Kelvin-Helmholtz instability testcase obatined with the relaxed multirate \ac{PERRK} scheme.}
	\label{fig:DensityMesh_KelvinHelmholtzInstability}
	\end{figure}
	\begin{figure}[!ht]
		\centering
		\includegraphics[width=0.55\textwidth]{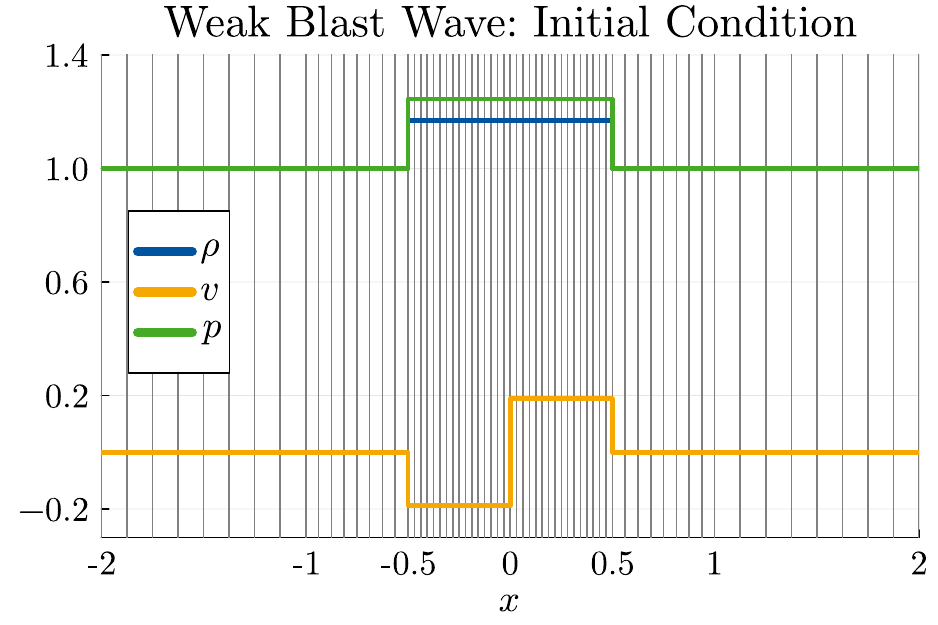}
		\caption[Initial hydrodynamic primal variables corresponding to the weak blast wave initial condition.]
		{\cref{subsec:EC_WeakBlastWaveEulerMHD}:
			Initial hydrodynamic primal variables corresponding to the weak blast wave initial condition \cref{eq:WeakBlastWave_IC}.
			The vertical lines depict the cell boundaries.
		}
		\label{fig:WeakBlastWave_hydrodynamic_variables_init}
	\end{figure}
	\begin{figure}[!ht]
		\centering
		\includegraphics[width=0.55\textwidth]{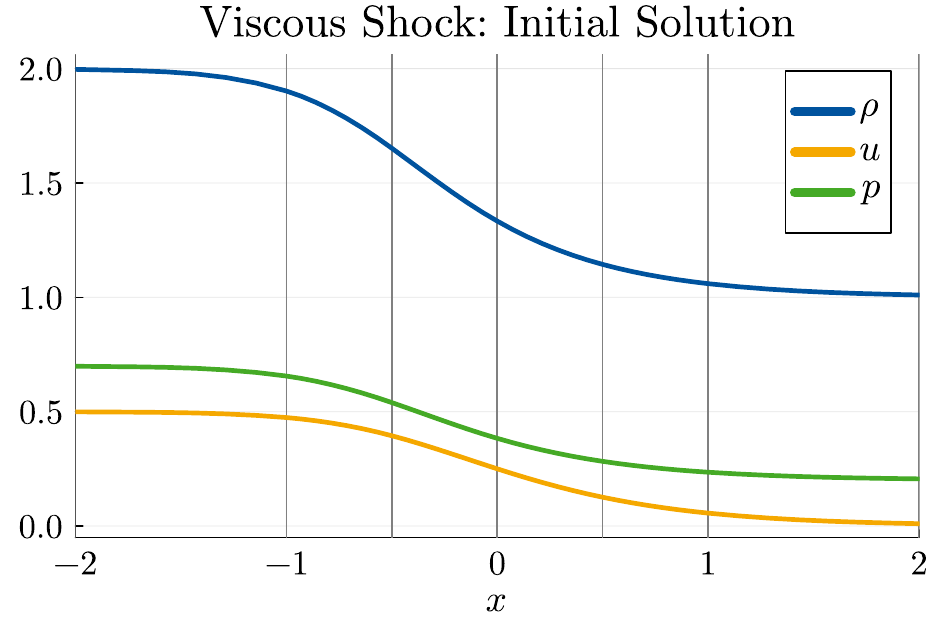}
		\caption[Initial primal variables corresponding to the viscous shock solution.]
		{\cref{subsec:ConvTest_ViscousShock}: Initial primal variables corresponding to the viscous shock solution \cref{eq:ViscousShockSolution}.
			The vertical lines depict the cell boundaries.
			Here, the coarsest mesh used in the convergence study with six cells only and $k=3$ solution polynomials is shown.
		}
		\label{fig:ViscousShock_Init}
	\end{figure}
	\begin{figure}[!ht]
		\centering
		\subfloat[{Relaxation parameter $\gamma_n$ over time for $p=2$.}]{
			\label{fig:gamma_VSP_p2}
			\centering
			\resizebox{.475\textwidth}{!}{\includegraphics{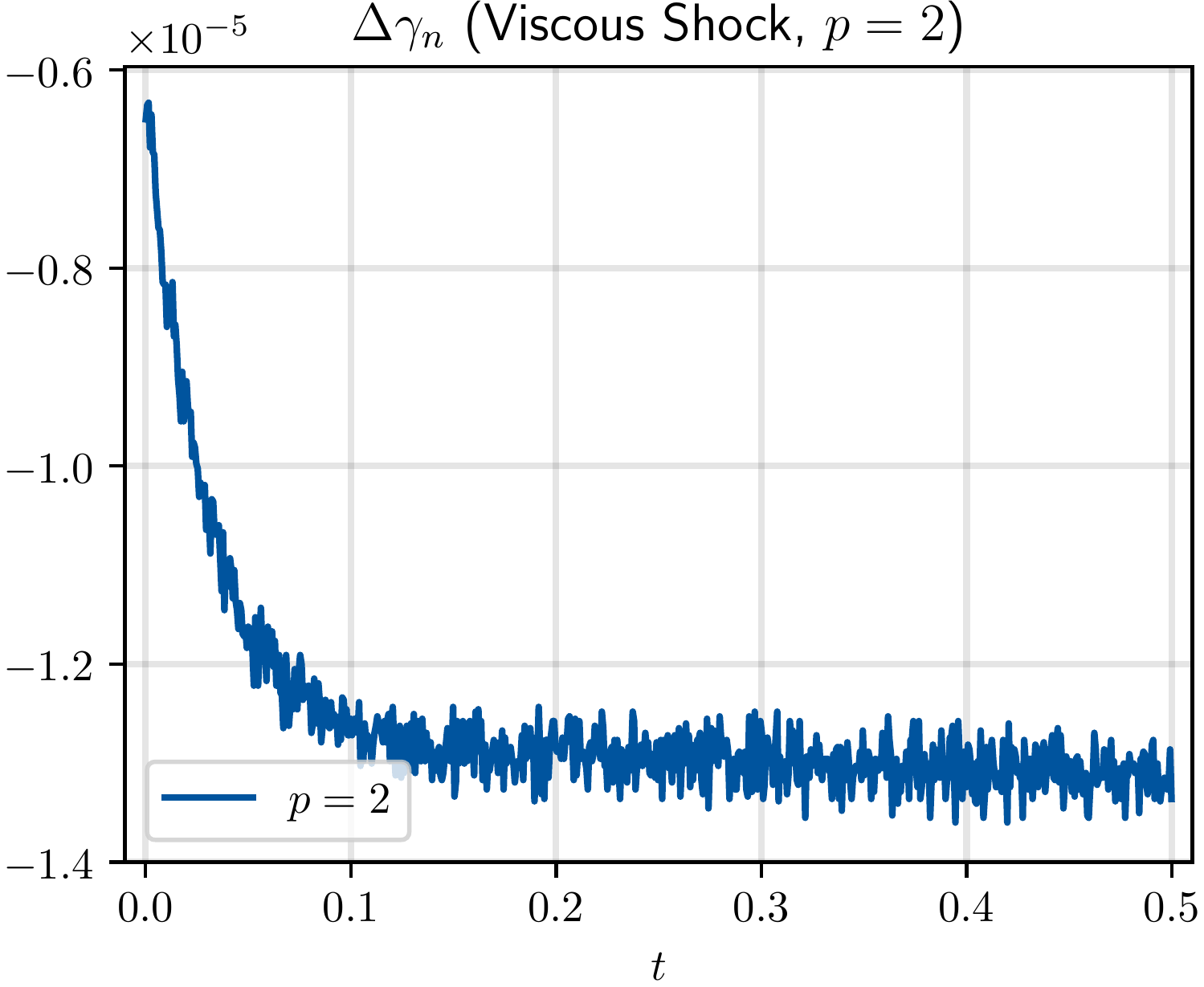}}
		}
		\hfill
		\centering
		\subfloat[{Relaxation parameter $\gamma_n$ over time for $p=3, 4$.}]{
			\label{fig:gamma_VSP_p3p4}
			\centering
			\resizebox{.475\textwidth}{!}{\includegraphics{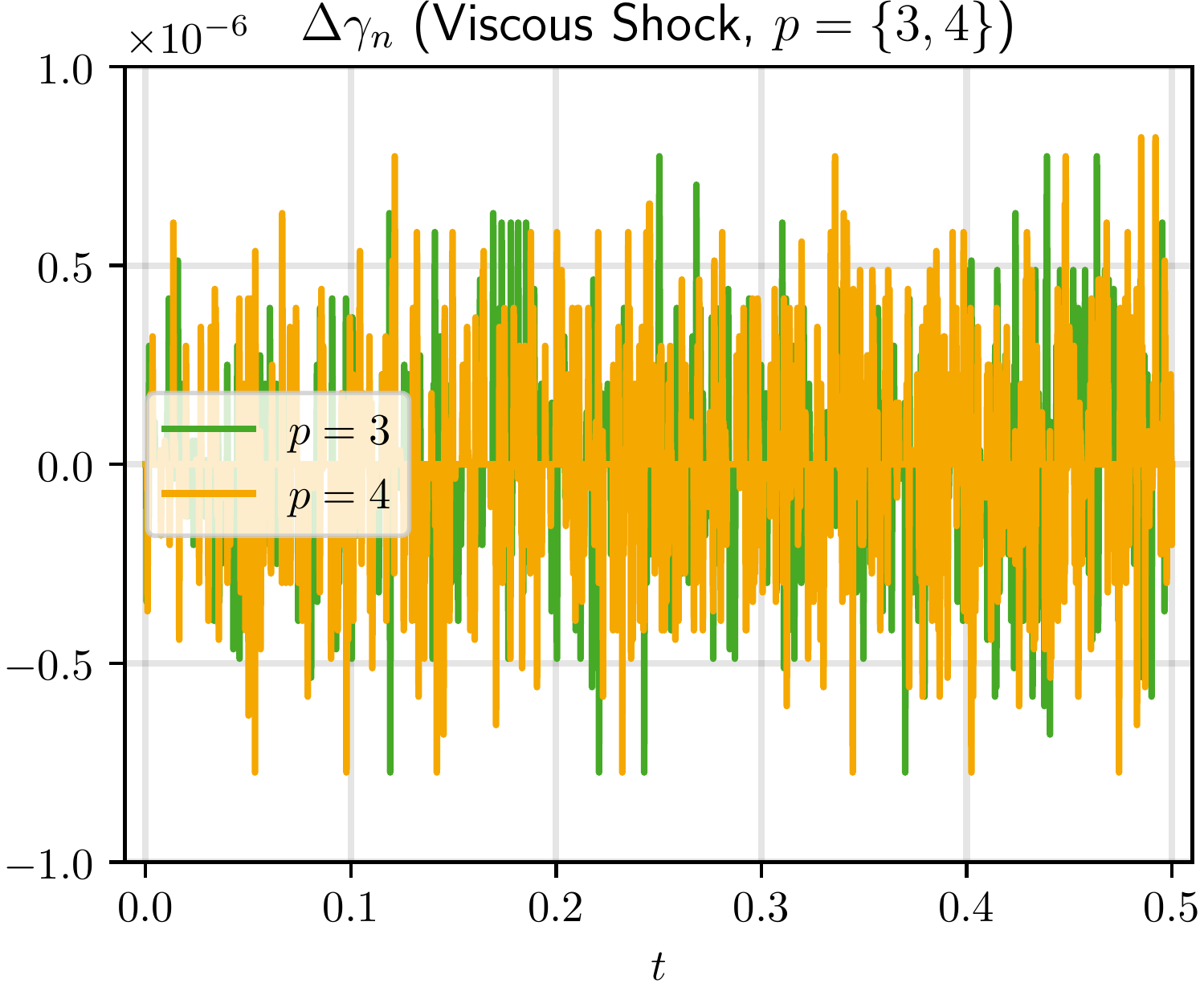}}
		}
		\caption[Evolution of the relaxation parameter $\gamma_n$ for the viscous shock convergence testcase with $N= 24$ grid cells.]
		{\cref{subsec:ConvTest_ViscousShock}:
		Evolution of the relaxation parameter $\gamma_n$ for the viscous shock convergence testcase with $N= 24$ grid cells for $p=2$ (\cref{fig:gamma_VSP_p2}) and $p=3, 4$ (\cref{fig:gamma_VSP_p3p4}) \ac{PERRK} schemes.
		In the plots above, the deviation from unity, i.e., $\Delta \gamma_n \coloneqq \gamma_n - 1$ is shown.}
		\label{fig:gamma_VSP_p2p3p4}
	\end{figure}
	\begin{figure}[!ht]
		\centering
		\subfloat[{Relaxation parameter $\gamma_n$ over time for $p=2$.}]{
			\label{fig:gamma_AW_p2}
			\centering
			\resizebox{.55\textwidth}{!}{\includegraphics{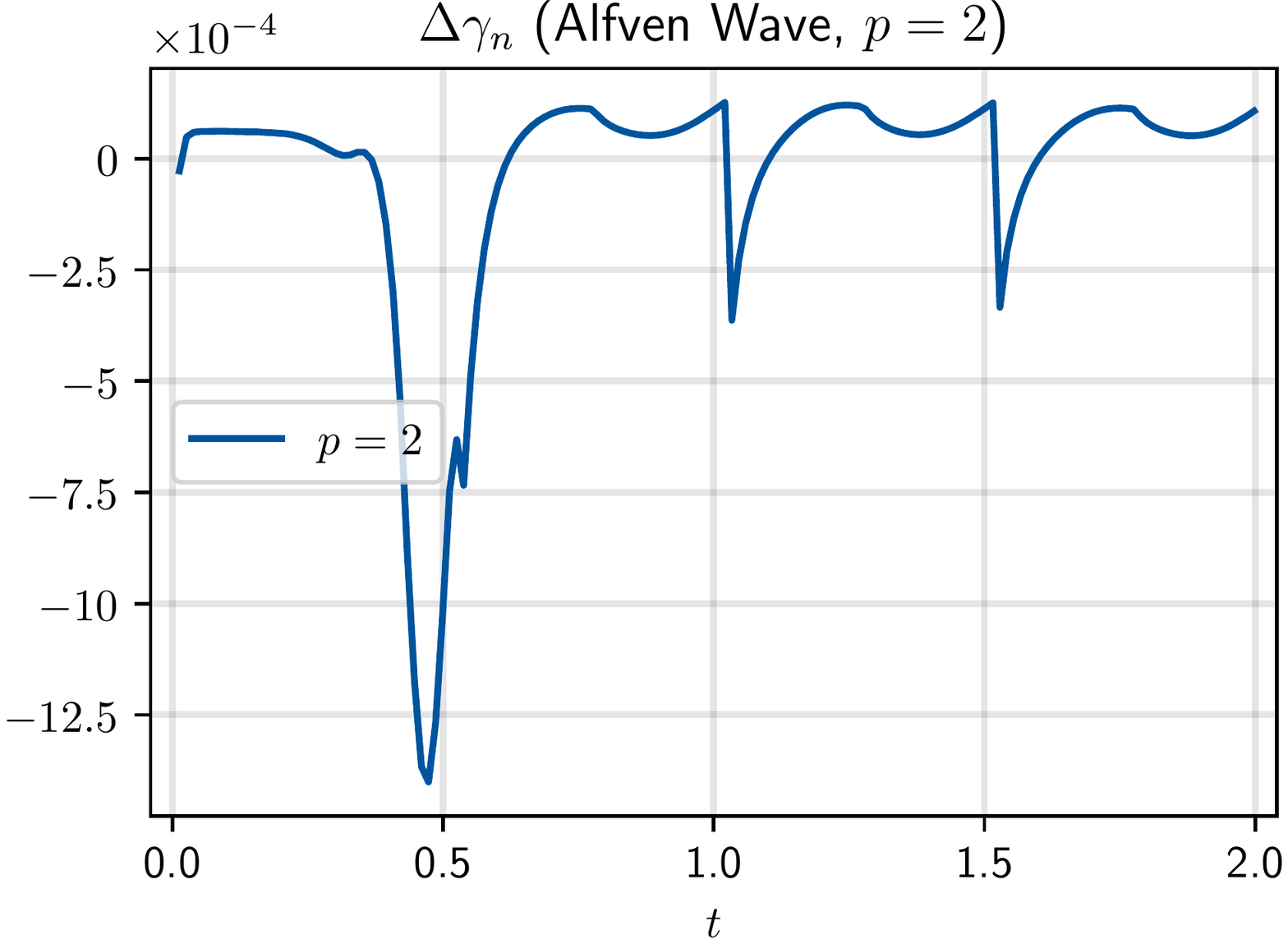}}
		}
		\\
		\subfloat[{Relaxation parameter $\gamma_n$ over time for $p=3$.}]{
			\label{fig:gamma_AW_p3}
			\centering
			\resizebox{.55\textwidth}{!}{\includegraphics{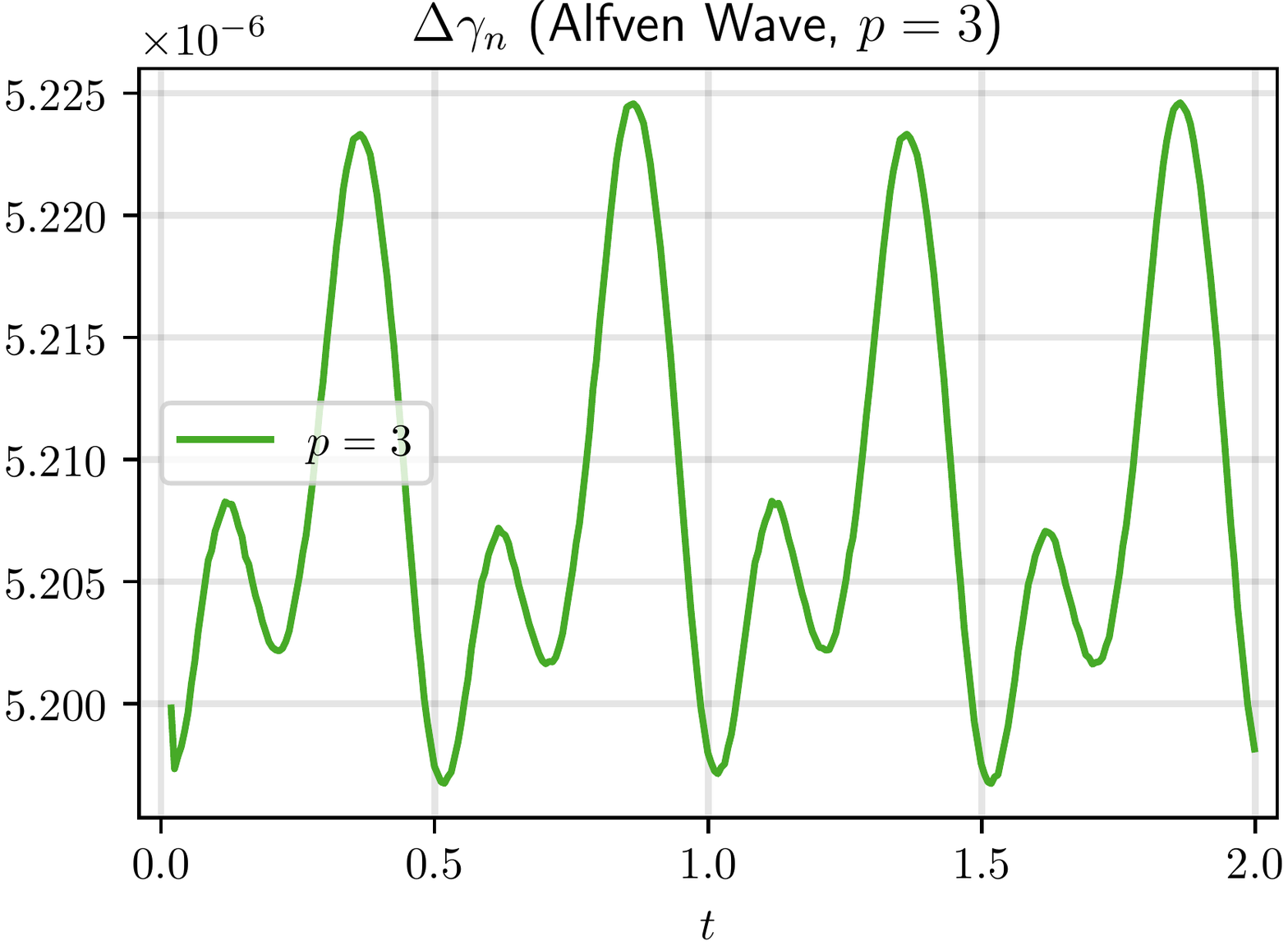}}
		}
		\\
		\subfloat[{Relaxation parameter $\gamma_n$ over time for $p=4$.}]{
			\label{fig:gamma_AW_p4}
			\centering
			\resizebox{.55\textwidth}{!}{\includegraphics{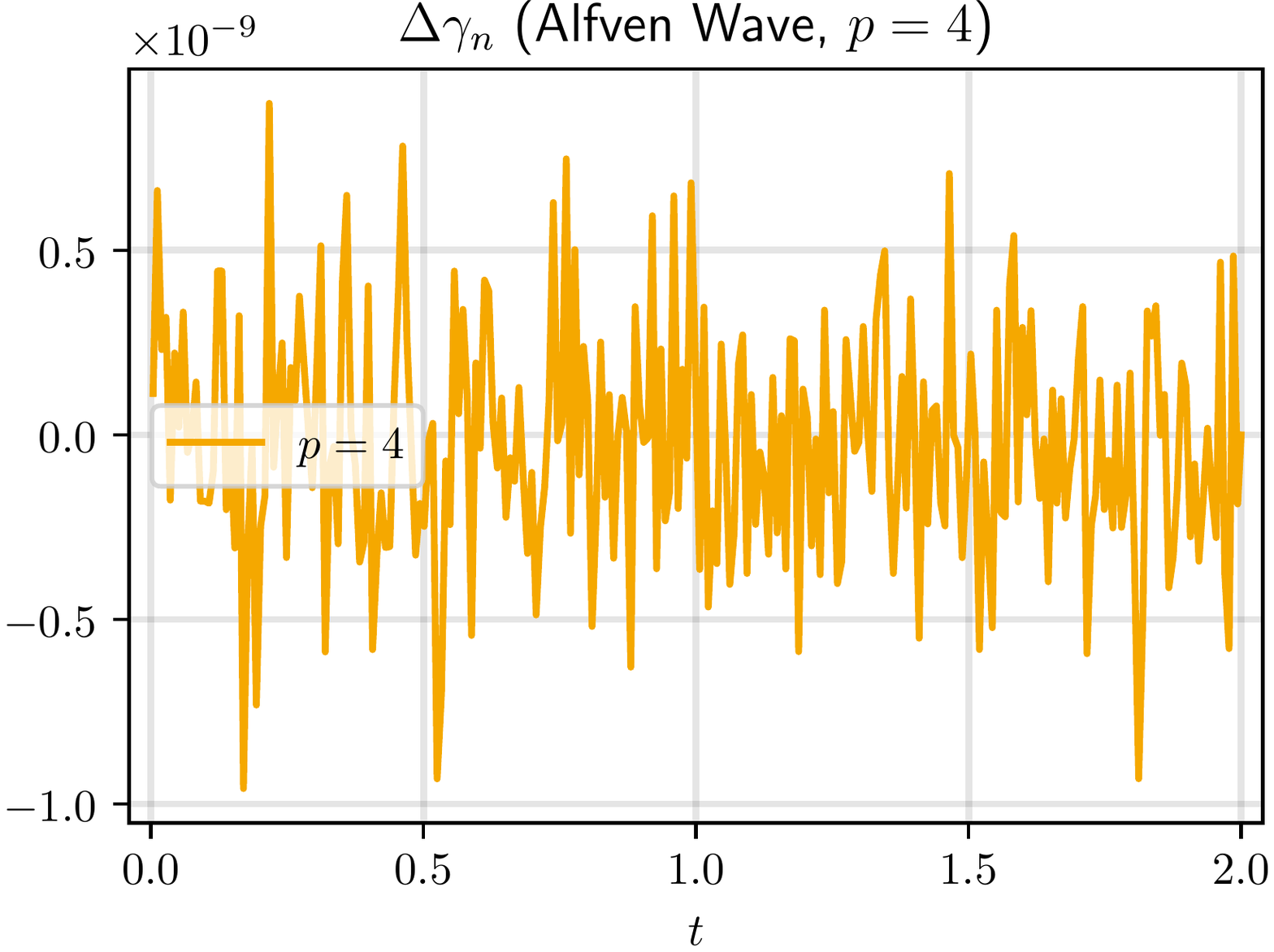}}
		}
		\caption[Evolution of the relaxation parameter $\gamma_n$ for the Alfvén wave convergence testcase with $N= 64$ grid cells.]
		{\cref{subsec:ConvTest_AlfvenWave}: 
		Evolution of the relaxation parameter $\gamma_n$ for the Alfvén wave convergence testcase with $N= 64$ grid cells for $p=2$ (\cref{fig:gamma_AW_p2}), $p=3$ (\cref{fig:gamma_AW_p3}), and $p=4$ (\cref{fig:gamma_AW_p4}) \ac{PERRK} schemes.
		In the plots above, the deviation from unity, i.e., $\Delta \gamma_n \coloneqq \gamma_n - 1$ is shown.}
		\label{fig:gamma_AW_p2p3p4}
	\end{figure}
	\begin{figure}[!ht]
		\centering
		\subfloat[{Domain-normalized $L^1$-error in $B_x$, cf. \cref{eq:L1ErrorDomainNormalized_Density}.}]{
			\label{fig:Convergence_AW_L1_RR}
			\centering
			\resizebox{.47\textwidth}{!}{\includegraphics{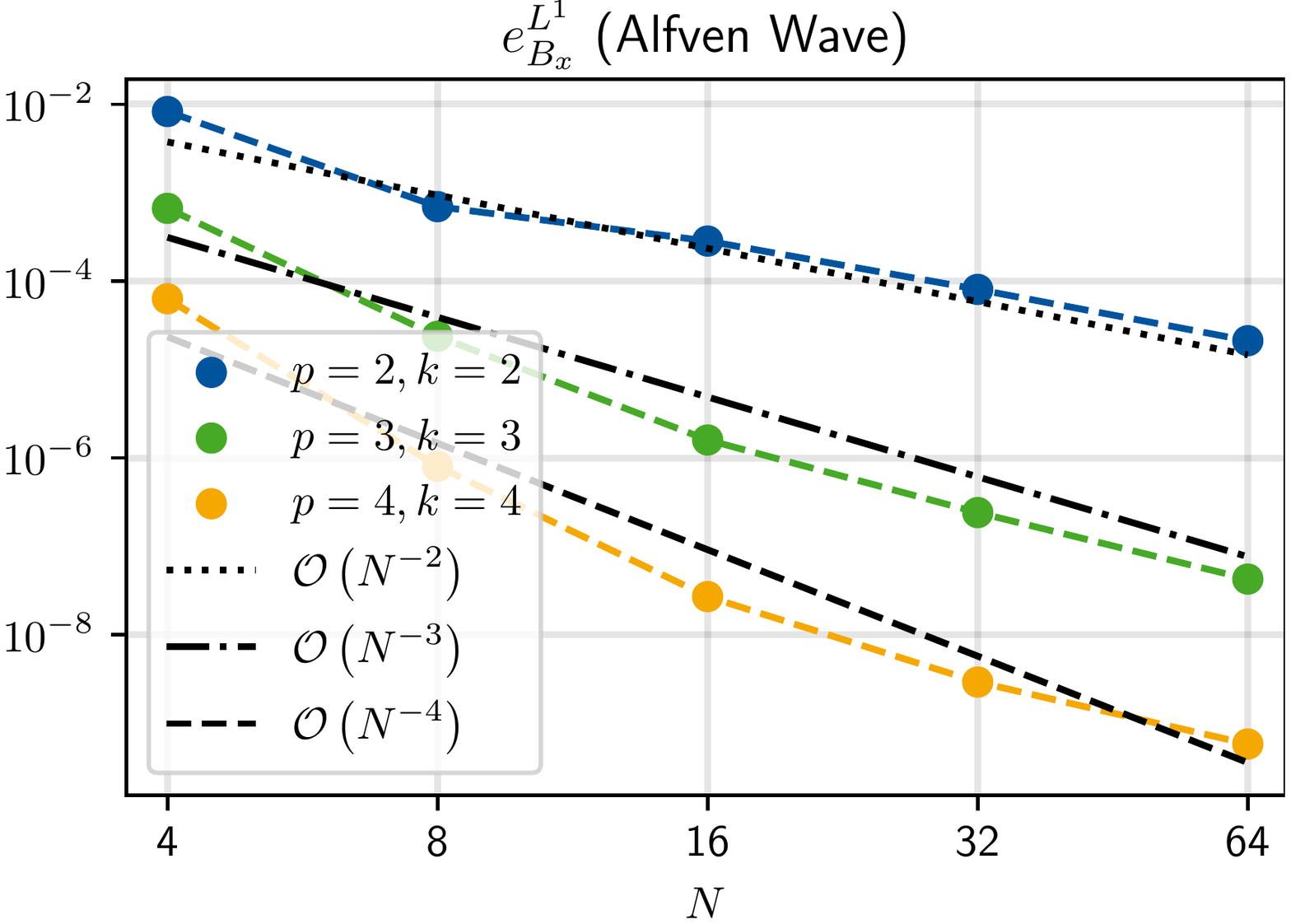}}
		}
		\hfill
		\subfloat[{$L^\infty$ error in $B_x$.}]{
			\label{fig:Convergence_AW_LInf_RR}
			\centering
			\resizebox{.47\textwidth}{!}{\includegraphics{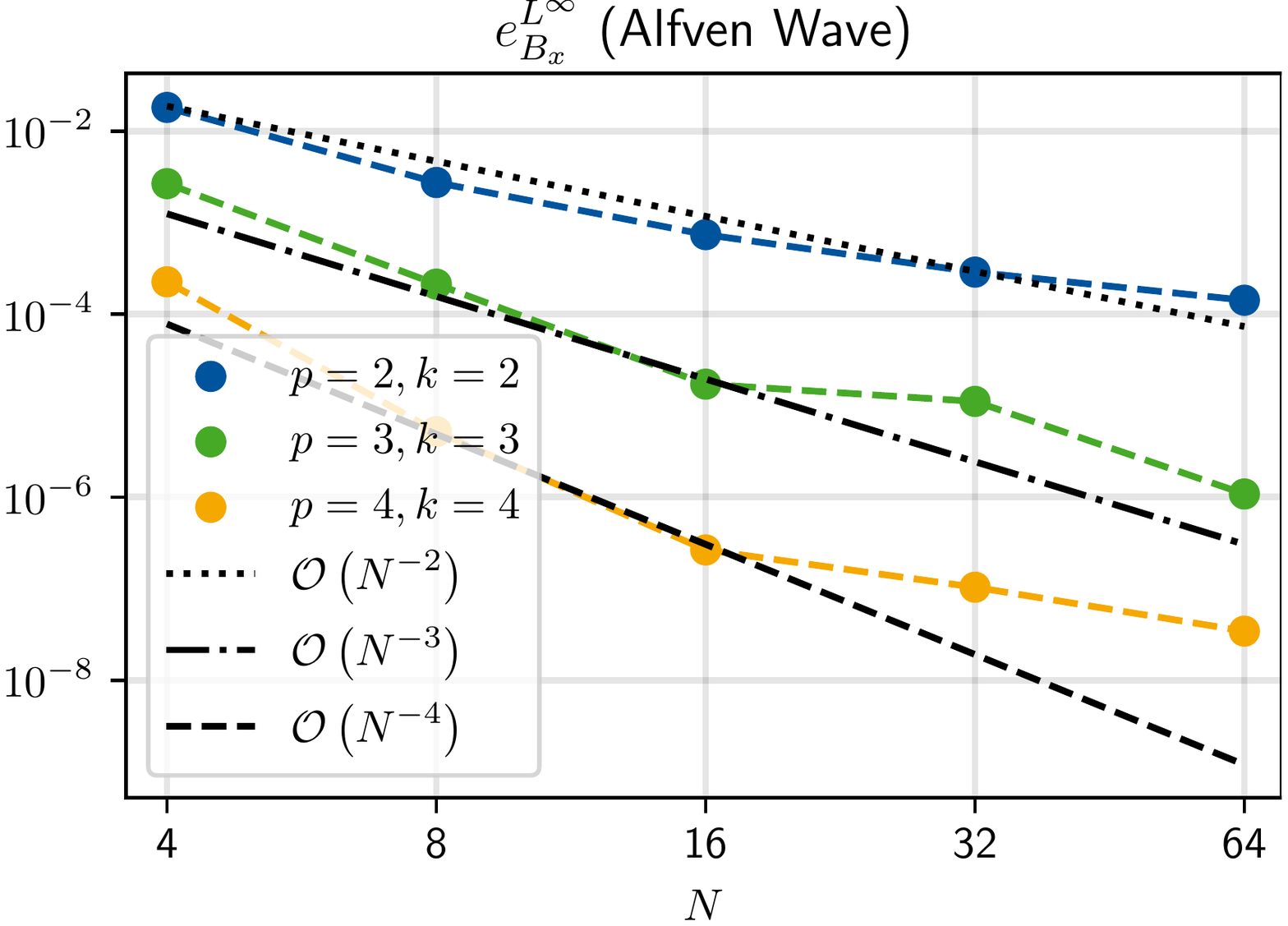}}
		}
		\caption[$L^1$ and $L^\infty$ errors of the $x$-component of the magnetic field for the Alfvén wave testcase.]
		{\cref{subsec:ConvTest_AlfvenWave}: $L^1$ (\cref{fig:Convergence_IVA_L1}) and $L^\infty$ (\cref{fig:Convergence_IVA_LInf}) errors of the $x$-component of the magnetic field $B_x$ for the Alfvén wave testcase at $t_f = 4.0$ for $p=2, 3, 4$ \ac{PERRK} schemes assigned in a round-robin fashion.}
		\label{fig:Convergence_AW_L1_LInf_RR}
	\end{figure}
	\begin{figure}[!ht]
		\centering
		\includegraphics[width=0.95\textwidth]{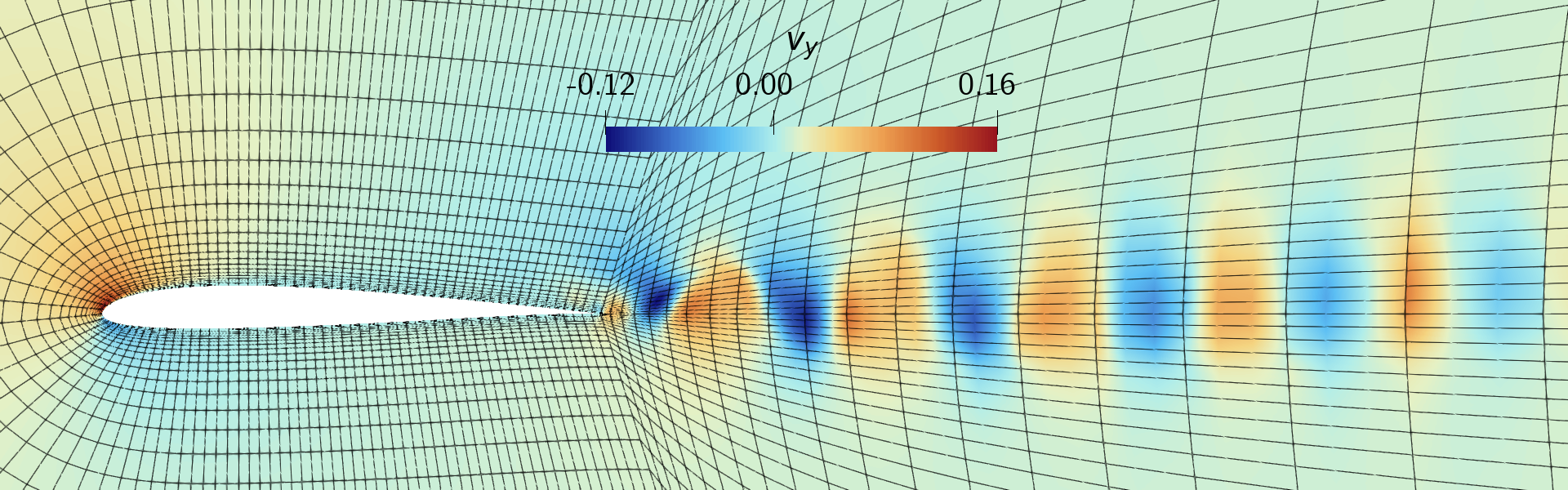}
		\caption[Vertical velocity $v_y$ around the SD7003 airfoil and wake at $t = 30 t_c\text{ s}$.]
		{\cref{subsec:SD7003}: 
		Vertical velocity $v_y$ around the SD7003 airfoil and wake at $t = 30 t_c\text{ s}$.
		}
		\label{fig:SD7003_vy}
	\end{figure}
	\begin{figure}[!ht]
		\centering
		\subfloat[{Lift coefficient $C_L = C_{L,p}$, cf. \cref{eq:LiftCoeff}.}]{
			\label{fig:LiftCoeff}
			\centering
			\resizebox{.47\textwidth}{!}{\includegraphics{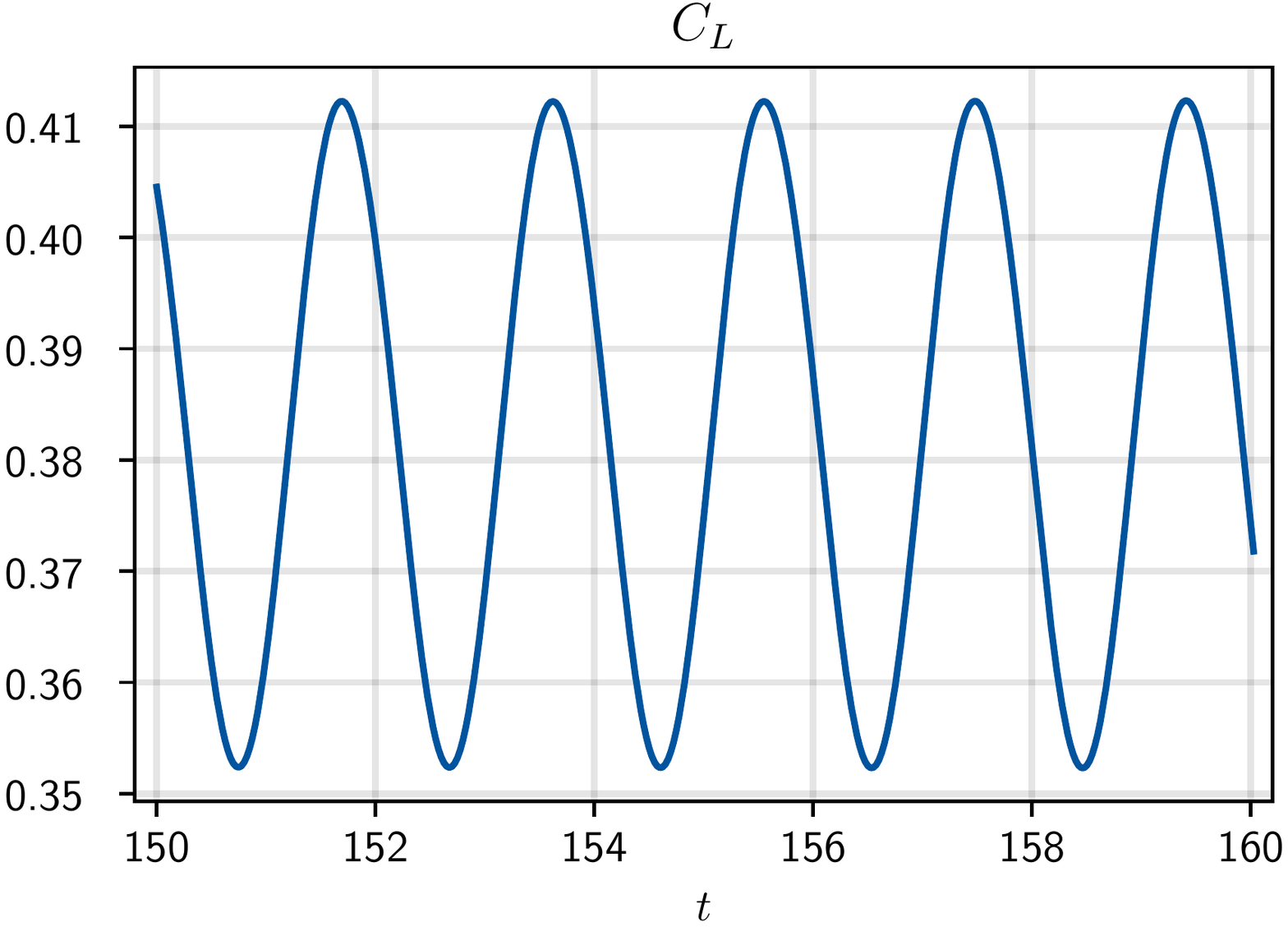}}
		}
		\hfill
		\subfloat[{Drag coefficient $C_D = C_{D,p} + C_{D, \mu}$, cf. \cref{eq:DragCoeffs}.}]{
			\label{fig:DragCoeff}
			\centering
			\resizebox{.47\textwidth}{!}{\includegraphics{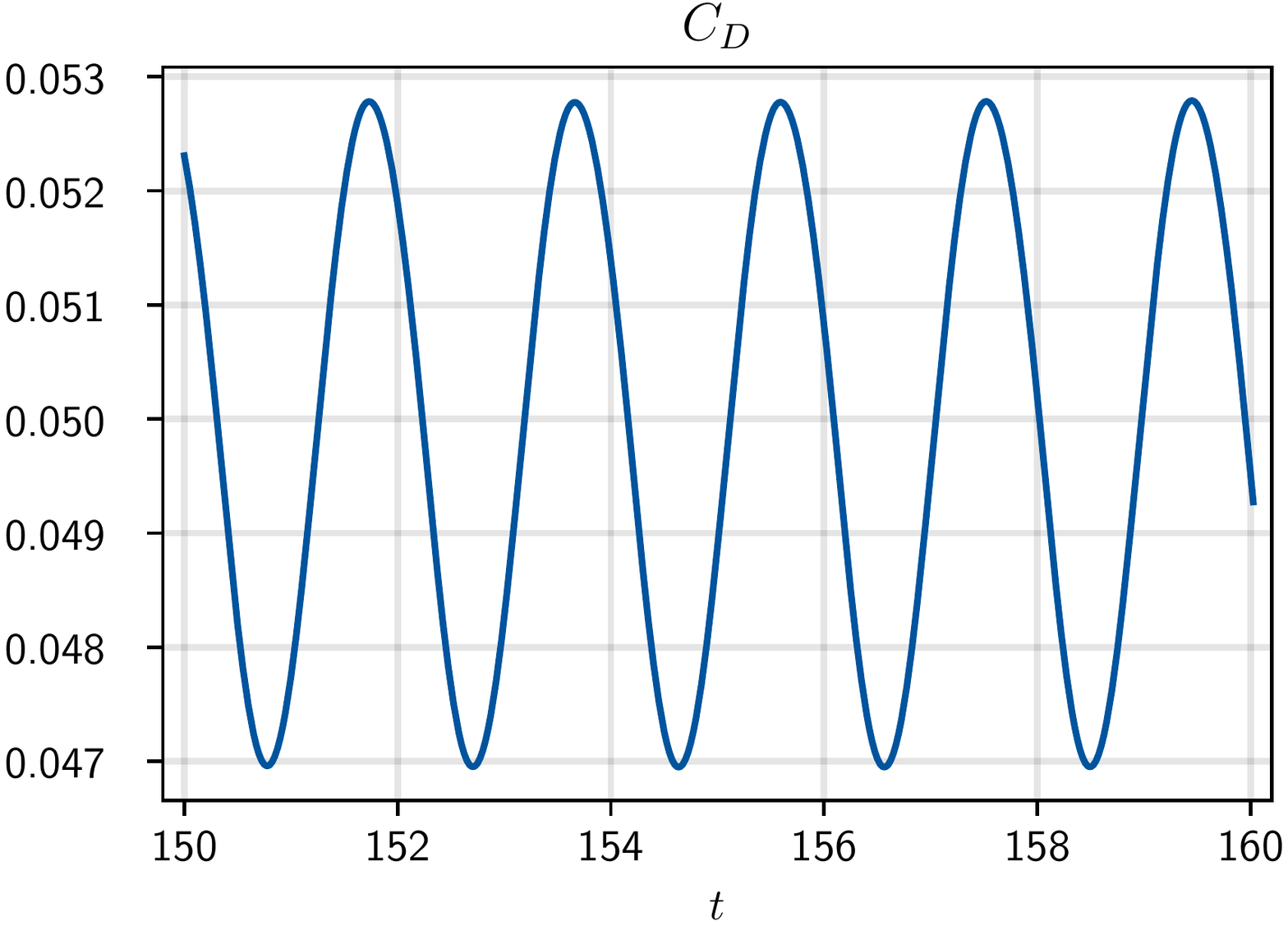}}
		}
	\caption[Unsteady lift and drag coefficients for the SD7003 airfoil.]
	{\cref{subsec:SD7003}: Unsteady lift (\ref{fig:LiftCoeff}) and drag coefficient (\ref{fig:DragCoeff}) for the SD7003 airfoil over the $[30 t_c, 32 t_c]$ time interval.}
	\label{fig:LiftDragCoeffs}
	\end{figure}
	\begin{figure}[!ht]
		\centering
		\includegraphics[width=0.75\textwidth]{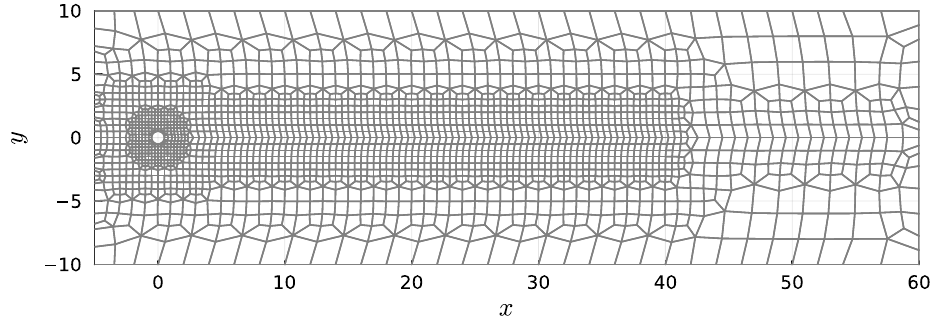}
		\caption[Computational mesh for the VRMHD flow past a cylinder.]
		{\cref{subsec:VRMHD_Cylinder}:
		Non-uniform, $y=0$ symmetric computational mesh composed of $2258$ quadrilateral elements employed for the \ac{vrMHD} flow past a cylinder.}
		\label{fig:Cylinder_Mesh}
	\end{figure}
	\begin{figure}[!ht]
		\centering
		\subfloat[{NACA0012 airfoil simulation on static mesh until $t_f = 100$.}]{
			\label{fig:NACA0012Static}
			\centering
			\includegraphics[width=0.75\textwidth]{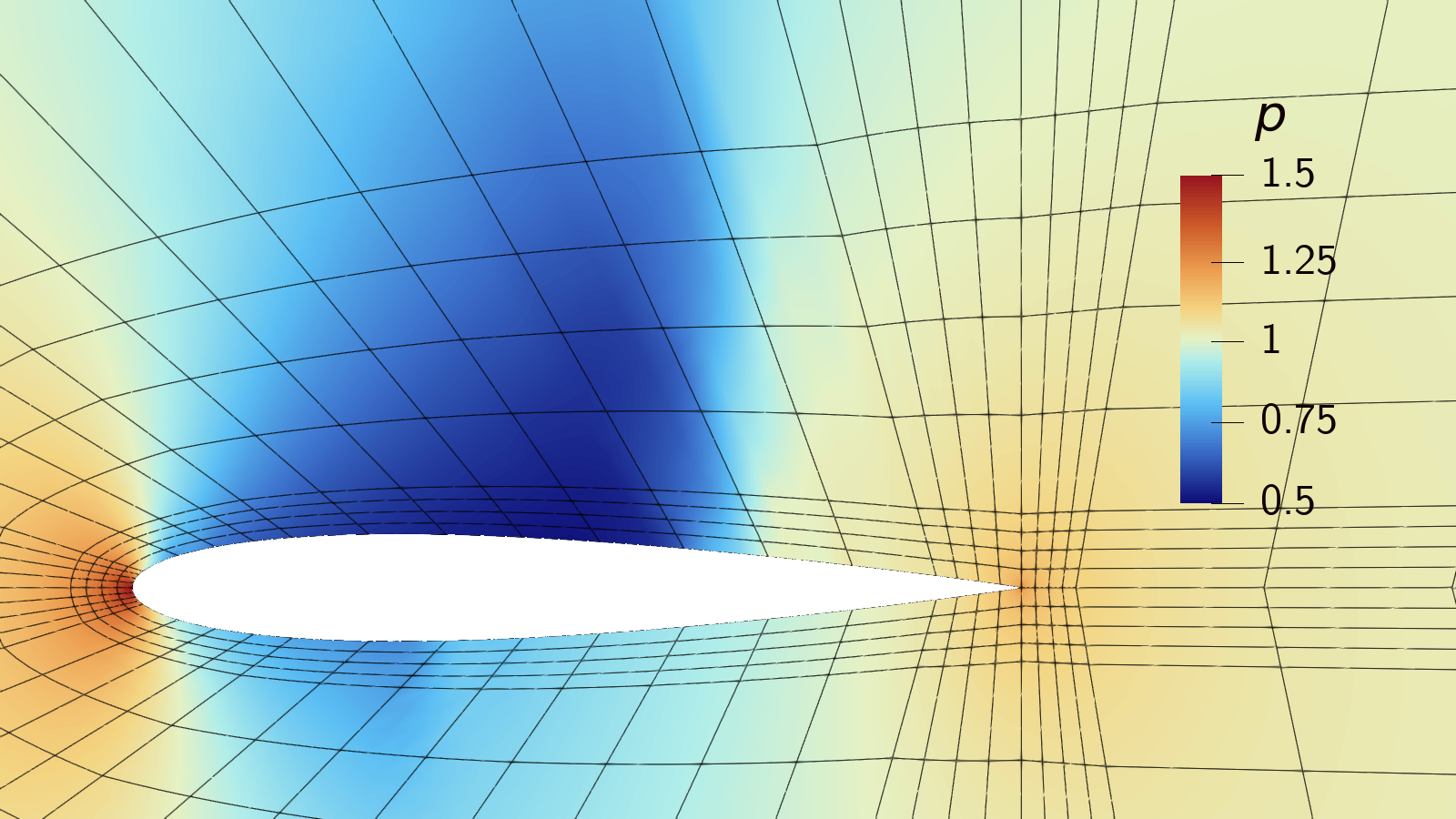}
		}
		\hfill
		\subfloat[{NACA0012 airfoil simulation on adaptively refined mesh until $t_f = 200$.}]{
			\label{fig:NACA0012AMR}
			\centering
			\includegraphics[width=0.75\textwidth]{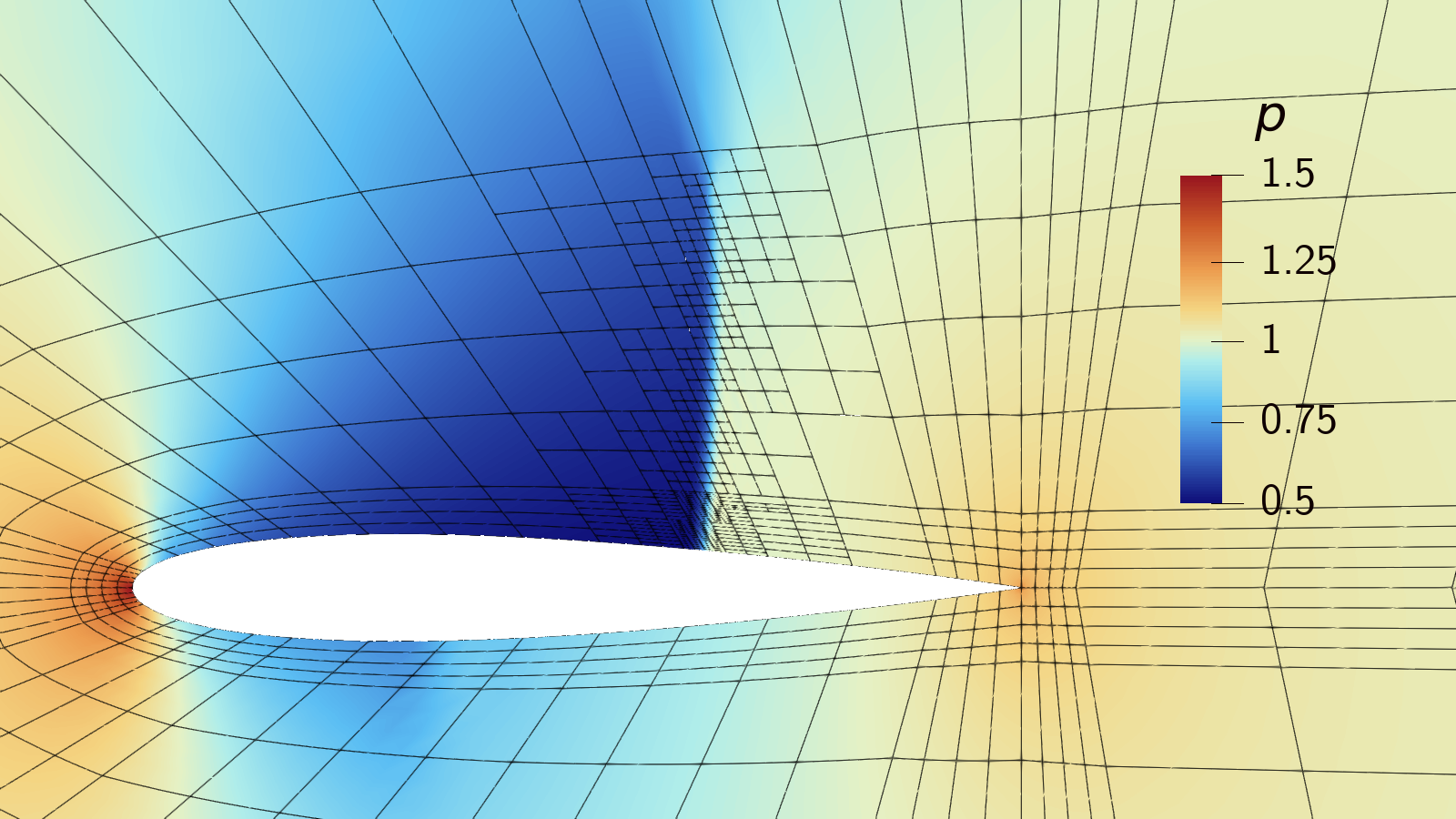}
		}
		\caption[Pressure $p$ at times $t_f = 100$ and $t_f = 200$ for the transonic inviscid flow over the NACA0012 airfoil.]
		{\cref{subsec:NACA0012}:
		Pressure $p$ at times $t_f = 100$ for the static mesh, cf. \cref{fig:NACA0012Static}, and $t_f = 200$ on dynamically adapted mesh, cf. \cref{fig:NACA0012AMR} for the transonic inviscid flow over the NACA0012 airfoil.
		The simulation with \ac{AMR} is restarted from the static simulation at $t_f = 100$ and then evolved until the final time $t_f = 200$.
		}
		\label{fig:NACA0012StaticAMR}
	\end{figure}
	\clearpage
	\bibliographystyle{elsarticle-num-names} 
	\bibliography{references.bib}
	\end{document}